%% file: Certainty_Equivalence_Control-Based_Heuristics_in_Multi-Stage_Convex_Stochastic_Optimization_Problems.tex
\documentclass[opre,nonblindrev]{informs3}

\OneAndAHalfSpacedXI % current default line spacing
\usepackage{endnotes}
\let\footnote=\endnote
\let\enotesize=\normalsize
\def\notesname{Endnotes}%
\def\makeenmark{\hbox to1.275em{\theenmark.\enskip\hss}}
\def\enoteformat{\rightskip0pt\leftskip0pt\parindent=1.275em
  \leavevmode\llap{\makeenmark}}

\usepackage[colorlinks=true,breaklinks=true,bookmarks=true,urlcolor=blue,
     citecolor=blue,linkcolor=blue,bookmarksopen=false,draft=false]{hyperref}

\def\EMAIL#1{\href{mailto:#1}{#1}}% When hyperref is used, otherwise outcomment
\usepackage{bold-extra}
\usepackage[utf8]{inputenc} % allow utf-8 input
\usepackage[T1]{fontenc}    % use 8-bit T1 fonts
\usepackage{hyperref}       % hyperlinks
\usepackage{url}            % simple URL typesetting
\usepackage{booktabs}       % professional-quality tables
\usepackage{nicefrac}       % compact symbols for 1/2, etc.
\usepackage{microtype}      % microtypography
\usepackage{xcolor}         % colors
\usepackage{dsfont,mathtools}
\usepackage[ruled,vlined,linesnumbered]{algorithm2e}
\usepackage{subcaption}
\usepackage{graphicx}
\usepackage{optidef}
\usepackage{enumitem}
\usepackage{diagbox}
\usepackage{tikz-cd}
\usepackage{multirow,arydshln} % for hdashline and multirow
\usepackage{tcolorbox}

\graphicspath{{}{figure/}}

\makeatletter
\newcommand{\mytag}[2]{%
  \text{#1}%
  \@bsphack
  \begingroup
    \@onelevel@sanitize\@currentlabelname
    \edef\@currentlabelname{%
      \expandafter\strip@period\@currentlabelname\relax.\relax\@@@%
    }%
    \protected@write\@auxout{}{%
      \string\newlabel{#2}{%
        {#1}%
        {\thepage}%
        {\@currentlabelname}%
        {\@currentHref}{}%
      }%
    }%
  \endgroup
  \@esphack
}
\makeatother

\usepackage[textwidth=\linewidth]{todonotes}		% for todo's

\newcommand{\bX}{\mathbf{X}}
\newcommand{\bp}{\mathbf{p}}
\newcommand{\bq}{\mathbf{q}}
\newcommand{\bW}{\mathbf{W}}
\newcommand{\bV}{\mathbf{V}}

\newcommand{\bI}{\mathbf{I}}

\newcommand{\bd}{\mathbf{d}}

\newcommand{\bx}{\mathbf{x}}

\newcommand{\calJ}{\mathcal{J}}
\newcommand{\calB}{\mathcal{B}}

\newcommand{\bY}{\mathbf{Y}}

\newcommand{\by}{\mathbf{y}}

\newcommand{\bU}{\mathbf{U}}

\newcommand{\bv}{\mathbf{v}}
\newcommand{\bw}{\mathbf{w}}
\newcommand{\barw}{\mathbf{\overline{w}}}

\newcommand\bC{\mathbf{C}}
\newcommand\bD{\mathbf{D}}

\newcommand\bu{\mathbf{u}}
\newcommand\bP{\mathbf{P}}
\newcommand\bM{\mathbf{M}}

\newcommand\bzero{\mathbf{0}}

\newcommand\bMN{\mathbf{M}^{(N)}}
\newcommand\bYN{\mathbf{Y}^{(N)}}

\newcommand\bEN{\mathbf{E}^{(N)}}

\newcommand{\R}{\mathbb{R}}

\newcommand\calO{\mathcal{O}}
\newcommand\calE{\mathcal{E}}
\newcommand\calI{\mathcal{I}}

\newcommand\calM{\mathcal{M}}
\newcommand\calN{\mathcal{N}}
\newcommand\calT{\mathcal{T}}

\newcommand\calD{\mathcal{D}}

\newcommand\calX{\mathcal{X}}

\newcommand\calF{\mathcal{F}}

\newcommand\calU{\mathcal{U}}

\newcommand\calP{\mathcal{P}}
\newcommand\calC{\mathcal{C}}
\newcommand\calL{\mathcal{L}}
\newcommand\calW{\mathcal{W}}
\newcommand\calTN{\mathcal{TN}}

\newcommand\bonet{\mathbf{1}^{\top}}
\newcommand\bone{\mathbf{1}}

\newcommand{\expect}[1]{\mathbb{E}\left[#1\right]}
\newcommand{\expectt}[1]{\mathbb{E}_{t-1}\left[#1\right]}
\newcommand{\var}[1]{\mathrm{var}\left[#1\right]}
\newcommand{\proba}[1]{\mathbb{P}\left(#1\right)}
\newcommand{\expectproj}[1]{\mathbb{E}_{\mathrm{proj}}\left[#1\right]}
\newcommand{\expectup}[1]{\mathbb{E}_{\mathrm{update}}\left[#1\right]}

\newcommand\cproj{c_{\mathrm{proj}}}

\newcommand{\abs}[1]{\left|#1\right|}

\newcommand{\norme}[1]{\left\| #1 \right\|}

\newcommand{\norminf}[1]{\left\| #1 \right\|_\infty}

\newcommand\bGPhi{\mathbf{G}^{\Phi}}
\newcommand\bGProj{\mathbf{G}^{\mathrm{Proj}}}
\newcommand\bGR{\mathbf{G}^{R_t}}
\newcommand\bHProj{\mathbf{H}^{\mathrm{Proj}}}
\newcommand\bHR{\mathbf{H}^{R_t}}
\newcommand\calHProj{\mathcal{H}^{\mathrm{Proj}}}
\newcommand\calHPhi{\mathcal{H}^{\Phi}}
\newcommand\calHR{\mathcal{H}^{R_t}}
\newcommand\slip{\mathfrak{S}_{\mathrm{lip}}}
\newcommand\slisse{\mathfrak{S}_{\mathcal{C}^2}}

\newcommand\proj{\mathrm{Proj}}
\newcommand\defeq{\stackrel{\triangle}{=}}

\usepackage{natbib}
 \bibpunct[, ]{(}{)}{,}{a}{}{,}%
 \def\bibfont{\small}%
\TheoremsNumberedThrough     % Preferred (Theorem 1, Lemma 1, Theorem 2)

\ECRepeatTheorems

\EquationsNumberedThrough    % Default: (1), (2), ...
\renewcommand{\theassumption}{\Roman{assumption}}

\begin{document}
%%%%%%%%%%%%%%%%

% Outcomment only when entries are known. Otherwise leave as is and
%   default values will be used.
%\setcounter{page}{1}
%\VOLUME{00}%
%\NO{0}%
%\MONTH{Xxxxx}% (month or a similar seasonal id)
%\YEAR{0000}% e.g., 2005
%\FIRSTPAGE{000}%
%\LASTPAGE{000}%
%\SHORTYEAR{00}% shortened year (two-digit)
%\ISSUE{0000} %
%\LONGFIRSTPAGE{0001} %
%\DOI{10.1287/xxxx.0000.0000}%

% Author's names for the running heads
% Sample depending on the number of authors;
% \RUNAUTHOR{Jones}
% \RUNAUTHOR{Jones and Wilson}
% \RUNAUTHOR{Jones, Miller, and Wilson}
% \RUNAUTHOR{Jones et al.} % for four or more authors
% Enter authors following the given pattern:
\RUNAUTHOR{Yan and Reiffers-Masson}

% Title or shortened title suitable for running heads. Sample:
% \RUNTITLE{Bundling Information Goods of Decreasing Value}
% Enter the (shortened) title:
\RUNTITLE{CEC-Based Heuristics in Stochastic Optimization Problems}

% Full title. Sample:
% \TITLE{Bundling Information Goods of Decreasing Value}
% Enter the full title:
\TITLE{Certainty Equivalence Control-Based Heuristics in Multi-Stage Convex Stochastic Optimization Problems}

% Block of authors and their affiliations starts here:
% NOTE: Authors with same affiliation, if the order of authors allows,
%   should be entered in ONE field, separated by a comma.
%   \EMAIL field can be repeated if more than one author
\ARTICLEAUTHORS{%
\AUTHOR{Chen Yan \footnote{Corresponding author}}
\AFF{ STATIFY, Inria, 38334 Saint Ismier, France; Biostatistics and Spatial Processes, INRAE, 84914 Avignon, France \EMAIL{chen.yan@inria.fr}} %, \URL{}}
\AUTHOR{Alexandre Reiffers-Masson}
\AFF{ IMT Atlantique, Lab-STICC, UMR CNRS 6285, Brest, France \EMAIL{alexandre.reiffers-masson@imt-atlantique.fr}}
% Enter all authors
} % end of the block

\ABSTRACT{%
We examine a multi-stage stochastic optimization problem characterized by stagewise-independent, decision-dependent noises with strict constraints. The problem assumes convexity in that, following a specific relaxation, it transforms into a deterministic convex program. The relaxation process is inspired by the principle of Certainty Equivalence Control, which substitutes uncertainties with their nominal values and requires the hard constraints to be satisfied only in an expected sense. Utilizing the solutions obtained from these convex programs, we propose two universal methodologies—re-solving-based and projection-based—to formulate feasible policies relevant to the original problem. These methodologies are subsequently amalgamated to develop a hybrid policy, equipped with a tuning parameter that manages the frequency of re-solving. We derive upper bounds on the gap between the performance of these heuristic policies and the optimal one. Under the Lipschitz-type regularity of the model, these bounds are proportional to the square root of the stochastic noise variance. Assuming additional $\calC^2$-smoothness regularity, an alternative bound, proportional to the variance of the stochastic noise, can be established—providing a refinement when variances are small. Our model provides a framework for dynamic decision-making under uncertainty, encompassing classic inventory and Markovian bandit problems while embracing a broader range of stochastic optimization challenges. We demonstrate our methods using numerical experiments on a network utility maximization problem.

}%

% Sample
%\KEYWORDS{deterministic inventory theory; infinite linear programming duality;
%  existence of optimal policies; semi-Markov decision process; cyclic schedule}

% Fill in data. If unknown, outcomment the field
% parametric optimization; constraint qualification} 

% \HISTORY{This paper was
% first submitted on April 12, 1922 and has been with the authors for
% 83 years for 65 revisions.}

\maketitle
%%%%%%%%%%%%%%%%%%%%%%%%%%%%%%%%%%%%%%%%%%%%%%%%%%%%%%%%%%%%%%%%%%%%%%

% Samples of sectioning (and labeling) in OPRE
% NOTE: (1) \section and \subsection do NOT end with a period
%       (2) \subsubsection and lower need end punctuation
%       (3) capitalization is as shown (title style).
%
%\section{Introduction.}\label{intro} %%1.
%\subsection{Duality and the Classical EOQ Problem.}\label{class-EOQ} %% 1.1.
%\subsection{Outline.}\label{outline1} %% 1.2.
%\subsubsection{Cyclic Schedules for the General Deterministic SMDP.}
%  \label{cyclic-schedules} %% 1.2.1
%\section{Problem Description.}\label{problemdescription} %% 2.

% Text of your paper here

\section{Introduction} 

Multi-stage stochastic optimization, with its intricate mix of uncertainty and decision-making across time, serves as a cornerstone in a multitude of real-world applications, including supply chain management, power systems operation, inventory control, and financial planning (\cite{shapiro2021lectures, pflug2014multistage, kuchler2009stability}). This paper focuses on a distinct variant of this type of problem, distinguished by convexity, stagewise-independent decision-dependent noises, and stringent constraints.

To tackle the inherent complexity of the problem, we utilize a specific relaxation process that transforms it into a deterministic convex program. This process is twofold. The first part, inspired by the well-established Certainty Equivalence Control (CEC) principle, replaces uncertainties with their nominal values—typically the mean—and acts as if these were the actual values. It has been demonstrated that, under certain conditions for the linear-quadratic stochastic control problem, this approach actually results in optimal control (\cite{simon1956dynamic, theil1957note, duchan1974clarification, runggaldier1981generalized}). The second part entails relaxing stringent constraints to be met only in an expected sense. This removes all uncertainties from the equation and offers a valuable approximation in handling intricate stochastic dynamics. 

The potency of the approach discussed above, within a multi-stage optimization framework, has been well demonstrated in inventory management (\cite{kunnumkal2008refined, jasin2012re, cooper2002asymptotic, secomandi2008analysis, bumpensanti2020re}) and Markovian bandits (\cite{hu2017asymptotically,ZayasCabn2017AnAO,Brown2020IndexPA,zhang2021restless,gast2021lp}) contexts, where the relaxed problems emerge as linear programs. Building on these successful applications, we broaden this methodology to encompass a more extensive range of problems, where the relaxed problem can now be a convex program. We propose two universal methodologies—re-solving and projection—to devise feasible policies applicable to the original problem. We then integrate these methodologies to derive a hybrid policy, featuring a tuning parameter to effectively regulate the frequency of re-solving.

In quantifying the effectiveness of these heuristic policies, we establish upper bounds on the gap between their performance and the optimal one. Establishing these upper bounds involves a meticulous manipulation of the propagation of stochastic errors, supplemented with concentration-type inequalities. Assuming Lipschitz-type regularity in the model, these bounds are proportional to the square root of the stochastic noise variance—a first-order estimation. An alternative bound, proportional to the variance of the stochastic noise, is established under the additional $\calC^2$-smoothness regularity. This second-order estimation provides a refinement when variances are small. Such regularity conditions are met on the model, given that the relaxed programs satisfy specific constraint qualifications. This analysis is facilitated by tools from parameterized optimization (\cite{bonnans2013perturbation} and \cite[Chapters 4,5]{facchinei2003finite}).

%We further delve into the intricacies of these bounds, focusing on the exponential growth of a multiplicative constant with respect to the stage number. From a computational complexity perspective, we argue that such exponential growth is generally unavoidable in a multi-stage optimization problem (\cite{dyer2006computational, shapiro2005complexity, reaiche2016note}).

\subsection*{Related Works and Our Contributions}

Our problem assumes an inherent convexity, demonstrated by the CEC-based relaxation process, transforming it into a deterministic convex program. As per our knowledge, previous applications of CEC in multi-stage stochastic optimization settings, such as inventory management or Markovian bandits mentioned previously, restrict themselves to linear programs. The inherent linearity allows for scaling within the model, with the scaling parameter governed by the system size and facilitating asymptotic limit analysis when the size grows large. In a more generalized convex setting, however, such scaling ceases to be applicable in general. Therefore, we shift the perspective to consider a single stochastic problem, presuming that the amplitude of the variance of stochasticities can be reduced—an approach similar to a density model in the linear case (see Remark \ref{rem:scaling-under-affinity} for details). Unquestionably, the absence of linearity also leads to more complex analysis and calculations.

We posit that this convexity requirement, to some extent, represents the most general framework to consider under this methodology. Several factors underpin this claim. Primarily, Jensen's inequality applies, ensuring the value of the relaxed convex program is \emph{larger} than the value of the optimal policy, which is typically elusive. Without this relationship, the upper bounds on the performance gap developed in this paper would become untenable (see Remark \ref{rem:necessitiy-of-the-assumtions} for details). The second factor is computational: convex programs straddle the boundary between efficiently solvable and intractable problems \cite[Lecture 5]{ben2001lectures}. The appeal of applying CEC diminishes if significant challenges already arise at the level of solving the relaxed problems.

In our model, the constraints are "hard": we cannot tolerate any violations, even minor ones. A major challenge we face is that the relaxed problem merely approximates the actual problem, and the decisions derived from it are generally not even feasible for the original problem. Thus, devising simple, efficient, and high-performing feasible policies for the original problem based on the solutions to the relaxed problems is a central theme in these CEC-based techniques. Examples include the LP-index policy (\cite{gast2021lp}) and randomized activation control policy (\cite{ZayasCabn2017AnAO}) for the Markovian bandit problem, as well as the booking limit policy and bid-price policy for inventory management (\cite{jasin2013analysis}). Another basic idea, based on re-solving and originating from Model Predictive Control (\cite{rakovic2018handbook}), is also well-documented (\cite{bumpensanti2020re, gast2022lp, brown2022fluid}). In this paper, we extend the re-solving concept, termed the "update policy," to a more general convex setting. Simultaneously, we also introduce, to the best of our knowledge, a new concept of projecting the relaxed solution onto the feasible set, which we term the "projection policy", thereby contributing another method for such constructions.

A noteworthy feature of our model, which sets it apart from classical multi-stage stochastic optimization problems as in \cite{shapiro2021lectures, pflug2014multistage, kuchler2009stability}, is the introduction of decision-dependent noise. This divergence is not purely academic; instead, it is motivated by a range of theoretical and practical applications we aim to address—from Markov decision processes (\cite{Puterman:1994:MDP:528623}), stochastic approximation (\cite{kushner2003stochastic}), computational complexity (\cite{papadimitriou1985games}), to reinforcement learning (\cite{sutton2018reinforcement}). The decision-dependent noise renders various techniques used to handle a standard multi-stage problem, such as quasi-Monte-Carlo approximations and scenario trees, not directly applicable. However, an advantage of the CEC-based heuristics lies in their ability to treat decision-dependent noise in the same manner as i.i.d. white noise, under the assumption of uniformly bounded variances, which we adopt in this paper (see \eqref{eq:law-of-stochastic-noise}).

Our second-order estimations hinge on local $\calC^2$-smoothness regularity. This characteristic has previously been noted in the simpler linear program case, as exemplified in \cite{jasin2012re, wu2015algorithms, bumpensanti2020re} for the inventory model, and in \cite{zhang2021restless, gast2021lp, brown2022fluid, gast2023exponential} for the Markovian bandit model. These studies introduce this additional regularity under various names such as non-degeneracy or non-singularity, with their precise definitions varying across different papers and contexts. We argue that this regularity property can be expressed as generally as in any finite-dimensional variational inequality (\cite{facchinei2003finite}), ultimately establishing a diffeomorphism that renders the constrained parameterized optimization problem locally akin to an unconstrained, parameter-free problem within a neighborhood of the origin in an appropriately dimensioned Euclidean space. We provide a unifying discussion in Appendix \hyperlink{appendix:policy-mapping}{EC.3}, building on existing works such as \cite{robinson1982generalized, robinson1987local, dunn1987convergence}. We anticipate this unified viewpoint may lend further insight into this property, especially considering its ability to accelerate convergence rates in the asymptotic regime (see Section \ref{subsec:interpretation-of-the-optimality-gap-bounds} for details).

\paragraph*{Outline}  The rest of the paper unfolds as follows. We introduce the general model in Section \ref{sec:general-model}, and subsequently present three motivating examples in Section \ref{sec:concrete-examples}: network utility maximization, inventory management, and Markovian bandit. The concept of CEC is articulated for the general model in Section \ref{sec:CEC}. Section \ref{sec:main-performance-bound} introduces heuristic policies based on CEC, encompassing the update policy, the projection policy, and a hybrid policy featuring a tuning parameter. We present upper bounds on the performance gap of the update and projection policies, predicated on specific regularity conditions, in Section \ref{sec:performance-results}. Numerical experiments focused on the network utility maximization example are offered in Section \ref{sec:numerical-experiments}. We conclude with additional comments and a discussion of future research directions in Section \ref{sec:conclusion}. Proofs, extensions, and additional discussions are in the electronic companion appendix.

\paragraph*{Notational Convention} 

We use bold letters to denote vectors and matrices, and vectors are represented as row vectors. Capital letters are used to denote random quantities and lowercase letters are used to denote deterministic quantities. The letters $\bx,\bu, \bW, \calE$ are reserved to represent system state, system control, exogenous and endogenous uncertainties, respectively. A norm $\norme{\cdot}$ without subscript is understood to be the $\calL^2$-norm for a vector, and the spectral norm for a matrix, so that it is compatible with the $\calL^2$-norm for vectors in the context of matrix-vector multiplication. We use $(\cdot)^\top$ to denote the transpose, so for instance $\bzero$ is a zero row vector, and $\bonet$ is a column vector of one's. We adopt the convention that time-step begins at $t=1$. For $1 \le t \le T$, denote by $\calF_t$ the $\sigma$-algebra  generated by the random quantities up to time-step $t$, prior to the system state transition, with the convention that $\calF_0 = \{ \bx \}$, where $\bx$ is the initial configuration of the stochastic system (the filtration $\calF_t$ is given more precisely in Assumption \ref{assumption-5} below). We write $\mathbb{E}_{\mathrm{update}} [\cdot]$ (resp. $\mathbb{E}_{\mathrm{proj}} [\cdot]$) to mean the expectation taken under the update policy (resp. the projection policy) that we shall analyse in this work. Suppose that $\bx$ takes values in a domain $\calX$, we denote by $\calB(\bx,\varepsilon)\defeq\left\{\bx' \in \calX \mid \norme{\bx' - \bx} \le \varepsilon \right\}$ the neighbourhood of $\bx$ of radius $\varepsilon$ in $\calX$. Denote by $(\bx,\bw,\bu)$ the row vector that concatenates three row vectors $\bx$, $\bw$ and $\bu$. We write $\bx[t,T]$ for the sequence of vectors $\bx(t')$ for $t \le t' \le T$, so that $\bx[t,T]$ is a short hand for the concatenation of vectors $\big( \bx(t), \bx(t+1), \dots, \bx(T)\big)$. For constraints that involve random quantities, we write \emph{a.s.} to mean that they are to be satisfied almost surely.

\section{The General Model}  \label{sec:general-model} 

Consider the following $T$-stage convex stochastic optimization problem with hard constraints:

\begin{tcolorbox}[colback=white,colframe=black,boxrule=1pt]

\emph{Known parameters}: initial condition of the system $\bx(1)$; horizon $T$; random vector $\bW$ with known distribution function $f(\bw)$; constraint functions $g_{t,i}(\cdot), h_{t,j}(\cdot)$ for $1 \le t \le T$ and $1 \le i \le I(t)$, $1 \le j \le J(t)$; system evolution with known Markovian laws $\phi(\cdot) + \calE(\cdot)$, where we have separated into two parts for later purpose: the deterministic part $\phi(\cdot)$ is \emph{affine}, the decision-dependent noise part $\calE(\cdot)$ is with zero mean; reward functions $R_t(\cdot)$ for $1 \le t \le T$.  

\

\emph{For each time-step $t=1, \dots, T$:}

\begin{enumerate}
  \item The decision-maker gets full knowledge of the current system state $\bX(t)$; 
  \item (The environment) independently draws $\bW(t) \sim f(\bw)$; 
  \item Once $\bW(t)$ has been observed, the decision-maker chooses a control $\bU(t)$ that satisfies the $I(t) + J(t)$ constrains $g_{t,i}(\bX(t), \bW(t), \bU(t)) \le 0 \ a.s.$ and $h_{t,j}(\bX(t), \bW(t), \bU(t)) = 0 \ a.s.$;
  \item The decision-maker collects a reward $R_t(\bX(t),\bW(t),\bU(t))$;
  \item The system evolves to the next state $\bX(t+1)$ such that $\bX(t+1) \sim \phi \left( \bX(t), \bW(t), \bU(t) \right) + \calE(\bX(t), \bW(t), \bU(t))$. 
\end{enumerate}

\

\emph{Objective}: Maximize the expected total sum of rewards over the $T$ time-steps.

\end{tcolorbox}

Mathematically, the problem can be formulated as follows:
\begin{maxi!}|s|{\bU[1,T]}{\mathbb{E} \left[ \sum_{t=1}^{T} R_t \left( \bX(t), \bW(t), \bU(t) \right)  \right] \label{eq:cost-overall} }{\label{eq:original-problem}}{V_{\mathrm{opt}} (\bx(1),T) =}
  \addConstraint{ \bX(1) = \bx(1) \ a.s.  \label{eq:initial-configuration} }
  \addConstraint{ g_{t,i}(\bX(t), \bW(t), \bU(t)) \le 0  \ a.s. \ \ \mbox{ for $ 1 \le t \le T$ and $ 1 \le i \le I(t)$} \label{eq:constraints} }{}
  \addConstraint{ h_{t,j}(\bX(t), \bW(t), \bU(t)) = 0  \ a.s. \ \ \mbox{ for $ 1 \le t \le T$ and $ 1 \le j \le J(t)$} \label{eq:constraints-eq} }{}
  \addConstraint{ \bX(t+1) = \phi \left( \bX(t), \bW(t), \bU(t) \right) + \calE(\bX(t), \bW(t), \bU(t)) \ a.s. \ \  \mbox{ for } 1 \le t \le T-1 \label{eq:evolution-of-the-system}}
\end{maxi!}
where for all $1 \le t \le T$ the model satisfies the following assumptions:
\begin{enumerate}
  \item   $\bx(t), \bX(t) \in \R^{n_x}$,  $\bU(t) \in \R^{n_u}$, $\bW(t) \in \R^{n_w}$ are \emph{continuous-valued} vectors, and are interpreted respectively as the system state (or configuration), system control (or decision, action), and exogenous uncertainties, with $\bx(1)$ being the \emph{deterministic} initial system configuration. We assume in addition that they all belong to \emph{convex} subsets of the corresponding Euclidean spaces.     \label{assumption-1}
  \item  $R_t: \R^{n_x} \times \R^{n_w} \times \R^{n_u} \rightarrow \R$ are real-valued \emph{concave} and $\calC^3$-smooth functions for all $t$, jointly for the three arguments. They are interpreted as reward or utility, and are \emph{additive} across times.  \label{assumption-2}
  \item  $\bW(t)$ are \emph{bounded independent and identically distributed}  (i.i.d.) random vectors with distribution function $f(\bw)$. We write
  \begin{equation*}
     \expect{\bW} \defeq \barw \in \R^{n_w} \mbox{ and } \var{\bW} \defeq \expect{(\bW - \barw) \cdot (\bW - \barw)^\top}
  \end{equation*}
 They are interpreted as exogenous uncertainties \emph{prior to} the decision-making $\bU(t)$ at time-step $t$.    \label{assumption-3}
  \item $g_{t,i}: \R^{n_x} \times \R^{n_w} \times \R^{n_u} \rightarrow \R$ are real-valued \emph{convex} and $\calC^2$-smooth functions for all $t$ and $1 \le i \le I(t)$, jointly for the three arguments, and $h_{t,j}: \R^{n_x} \times \R^{n_w} \times \R^{n_u} \rightarrow \R$ are real-valued \emph{affine} functions for all $t$ and $1 \le j \le J(t)$, jointly for the three arguments. They are interpreted as hard constraints to the system and $I(t), J(t) \in \mathbb{N}$ are the numbers of (in)equality constraints at time-step $t$.   \label{assumption-4}
  \item Given $\bX(t), \bW(t)$ and $\bU(t)$, the evolution from $\bX(t)$ to $\bX(t+1)$ is Markovian with an \emph{affine} behavior in expectation. Consequently, denote by $\calF_t$ the $\sigma$-algebra generated by the random vectors $\bX(t'), \bW(t')$, $\bU(t')$ for $1 \le t' \le t$, plus the random vectors $\calE(\bX(t'), \bW(t'), \bU(t'))$ for $1 \le t' \le t-1$ (note that $\calE(\bX(t), \bW(t), \bU(t))$ is not included in $\calF_t$, but rather in $\calF_{t+1}$),  then
      \begin{equation}\label{eq:law-of-Markovian-evolution}
        \expect{\bX(t+1) \mid \calF_t} = \expect{\bX(t+1) \mid \bX(t), \bW(t), \bU(t)} \defeq \phi \left( \bX(t), \bW(t), \bU(t) \right)
      \end{equation}
  where $\phi: \R^{n_x} \times \R^{n_w} \times \R^{n_u} \rightarrow \R^{n_x}$ is an affine function: there exists $\bC \in \R^{(n_x+n_w+n_u) \times n_x}$ and $\bD \in \R^{n_x}$ such that 
  \begin{equation}\label{eq:phi-explicit}
    \phi \left( \bX(t), \bW(t), \bU(t) \right) = (\bX(t),\bW(t),\bU(t)) \cdot \bC + \bD
  \end{equation}
 A priori $\bC$ and $\bD$ may also depend on $t$. The random decision vector $\bU(t)$ is $\calF_t$-measurable, and therefore cannot depend on any information beyond that available at time-step $t$, adhering to the principle of \emph{non-anticipative constraints}. \label{assumption-5}
  \item  Write $\calE(\bX(t), \bW(t), \bU(t)) \defeq \bX(t+1) - \phi \left( \bX(t), \bW(t), \bU(t) \right)$, we deduce from Assumption \ref{assumption-5} that 
      \begin{align*}  
       & \expect{\calE(\bX(t), \bW(t), \bU(t)) \mid \calF_t} = \expect{\calE(\bX(t), \bW(t), \bU(t)) \mid \bX(t), \bW(t), \bU(t)} \\
       & = \int_{\bx \in \R^{n_x}} \bx \ d \nu( \bx \mid \bX(t), \bW(t), \bU(t)) = \mathbf{0}   
      \end{align*}
     where $\nu( \cdot \mid \bX(t), \bW(t), \bU(t))$ is a probability distribution on $\R^{n_x}$ parameterized by $(\bX(t), \bW(t), \bU(t))$. Denote furthermore by
      \begin{align*}
       & \var{\calE(\bX(t), \bW(t), \bU(t))} = \var{\calE(\bX(t), \bW(t), \bU(t)) \mid \bX(t), \bW(t), \bU(t)} \\
       & = \int_{\bx \in \R^{n_x}} \bx \cdot \bx^\top \ d \nu(\bx \mid \bX(t), \bW(t), \bU(t))
      \end{align*}
    We suppose that 
      \begin{itemize}
        \item The support of $\nu(\cdot \mid \bx,\bw,\bu)$ is \emph{bounded} uniformly on $(\bx,\bw,\bu)$ 
      \end{itemize}
      which implies in particular that 
      \begin{equation} \label{eq:law-of-stochastic-noise}
        \var{\calE} \defeq \sup_{\bx, \bw, \bu} \var{\calE(\bx, \bw, \bu)} < \infty
      \end{equation}
      $\calE(\bX(t), \bW(t), \bU(t))$ are interpreted as endogenous uncertainties \emph{posterior to} the decision-making $\bU(t)$ at time-step $t$. They are decision-dependent noises. 
      \label{assumption-6}
  \item The problem has \emph{relatively complete recourse} \cite[Chapter 3]{shapiro2021lectures}. More precisely, for any system configuration $\bX(t)$, any realisation $\bW(t)$, there exists at least one decision $\bU(t)$ that satisfies the hard constraints \eqref{eq:constraints} and \eqref{eq:constraints-eq} at time-step $t$. This assumption is to ensure that the constrained problem \eqref{eq:original-problem} is feasible at any moment, so that the feasible region $\calU_t(\bx, \bw)$ given in \eqref{eq:feasible-region-time-t} are always non-empty.  
  \label{assumption-7}
\end{enumerate}
In the above formulation, for each time-step $t$, just before taking our decision, we observe a realization of the random vector $\bW(t)$. The vectors $\bX(t),\bW(t),\bU(t)$ are constrained by $I(t)$ inequality convex functions given by \eqref{eq:constraints}, plus $J(t)$ equality affine functions given by \eqref{eq:constraints-eq}. The system then evolves to the next state $\bX(t+1)$ according to \eqref{eq:evolution-of-the-system} in a Markovian way, which is an affine function of $(\bX(t),\bW(t),\bU(t))$ for the deterministic part, plus the zero-mean stochastic part $\calE(\bX(t), \bW(t), \bU(t))$ governed by the probability measure $\nu(\cdot \mid \bX(t), \bW(t), \bU(t))$. Our goal is to maximize the sum of rewards over the whole horizon, where the instantaneous reward $R_t(\bX(t),\bW(t),\bU(t))$ being a utility function is concave.

Denote by
\begin{equation}\label{eq:feasible-region-time-t}
  \calU_t(\bx,\bw) \defeq \left\{ \bu \in \R^{n_u} \ \big| \ g_{t,i} (\bx,\bw,\bu) \le 0, \ \ \mbox{ for }  1 \le i \le I(t); \ \ h_{t,j} (\bx,\bw,\bu) = 0, \ \ \mbox{ for } 1 \le j \le J(t) \right\}
\end{equation} 
which is the set of feasible controls at time-step $t$, given $\bx = \bX(t)$ and $\bw = \bW(t)$. By our previous assumptions on the model, $\calU_t(\bx,\bw)$ is a non-empty (Assumption \ref{assumption-7}), convex (Assumption \ref{assumption-4}), and compact (Assumptions \ref{assumption-3} and \ref{assumption-6}) set in $\R^{n_u}$. We point out that it is possible to relax the boundedness condition on the distributions of $\bW$ and $\calE$ by requiring that they have a tail distribution that converges to zero exponentially fast at infinity, as e.g. a Gaussian distribution, by combining with a concentration inequality as in Lemma \ref{lem:multivariate-hoeffding} on all our subsequent analysis. We impose the stronger uniformly bounded condition to simplify the matter, since then we can suppose that all the vectors $(\bX(t),\bW(t),\bU(t))$ for $1 \le t \le T$ take values on a common bounded and closed set, and invoke results using compactness.

In order to justify that the optimization problem \eqref{eq:original-problem} is well formulated under such generality, we need to impose two additional technical assumptions on the model.

\begin{enumerate}[resume]
  \item \begin{enumerate}
          \item (Slater CQ for feasible sets) For all $1 \le t \le T$ and all possible realization $(\bx,\bw)$ of $(\bX(t),\bW(t))$, there exists $\bu' \in \calU_t(\bx,\bw)$ such that $g_{t,i}(\bx,\bw,\bu') < 0$ for all $1 \le i \le I(t)$. 
          \item (Lipschitz-continuity on $\nu(\cdot \mid \bx,\bw,\bu)$) There exists $c_{\nu} > 0$ such that for all $(\bx_1,\bw_1,\bu_1)$ and $(\bx_2,\bw_2,\bu_2)$, we have 
              \begin{equation*}
                K(\nu(\cdot \mid \bx_1,\bw_1,\bu_1), \nu(\cdot \mid \bx_2,\bw_2,\bu_2)) \le c_{\nu} \cdot \norme{(\bx_1,\bw_1,\bu_1) - (\bx_2,\bw_2,\bu_2)}
              \end{equation*}
              where $K(\nu_1,\nu_2)$ is the Kantorovich distance between two probability distributions, or equivalently, the Wasserstain $1$-distance via the Kantorovich duality theorem \cite[Chapter 1]{villani2021topics}.
        \end{enumerate} 
  \label{assumption-8}
\end{enumerate}

We remark that the Slater CQ is a common assumption in convex programming, and the Lipschitz property on the probability measures is a control on the speed of change of the Markovian transition laws as a function of the input $(\bx,\bw,\bu)$. Under this additional Assumption \ref{assumption-8}, we can justify that the optimization problem \eqref{eq:original-problem} is well-formulated, the proof is provided in Appendix \hyperlink{appendix:proof-well-defined}{EC.1}.

\begin{proposition} [Existence of Optimal Solutions of \eqref{eq:original-problem}]  \label{prop:existence-of-optimal-solution}
  Under Assumptions \ref{assumption-1}-\ref{assumption-8}, the optimization problem \eqref{eq:original-problem} is well-defined and an optimal solution exists. Moreover, the mapping $\bx \mapsto V_{\mathrm{opt}} (\bx,T)$ is a continuous function of $\bx$. 
\end{proposition}

Note that an optimal solution to \eqref{eq:original-problem}, denoted as $U_{\mathrm{opt}}[1,T]$, is a series of random variables, where $U_{\mathrm{opt}}(t) \in \calF_t$, for $1 \le t \le T$. 

\begin{remark}[Discussion on the Model Assumptions]  \label{rem:model-assumptions}
  The convexity assumption, which plays an essential role in our method, will be justified later in Remark \ref{rem:necessitiy-of-the-assumtions}, after introducing the certainty equivalence control on the problem. Formulated as a multistage stochastic optimization problem, the model has \emph{stagewise-independent} noise \cite[Chapter 3]{shapiro2021lectures}. It incorporates two sources of stochasticity, $\bW$ and $\calE$, that have distinct nature and are motivated from different application scenarios, as we shall see in the examples displayed in Section \ref{sec:concrete-examples}. Clearly the i.i.d. noise $\bW$ is a special and simplified case from the decision-dependent noise $\calE$. Throughout the analysis of the paper we are only concerned with their first and second moments, and it turns out that this subtlety does not play an essential role in our method. Additinally, if the model lacks decision-dependent noise \(\calE\), the Sample Average Approximation, a Monte Carlo method, can be utilized to approximate the optimal solution of \eqref{eq:original-problem} \cite[Chapter 5]{shapiro2021lectures}. Let us also note that state augmentation can be applied by defining \(\widetilde{\bX} \defeq (\bX,\bW)\). Under this setup, the \(\bW\)-part acts as an exogenous and uncontrollable segment of the system state \(\widetilde{\bX}\), as seen in, for example, \cite{brown2022strength}. While this alternative modeling does not introduce technical changes to our subsequent discussions, separating \(\bX\) and \(\bW\) allows algorithms to operate over a reduced state space \cite[Section 1.4]{bertsekas2012dynamic}.
   
  In terms of computational complexity, already for linear two-stage stochastic programs with fixed recourse, it is shown in \cite{hanasusanto2016comment} that it is \#P-hard to find an approximate solution with sufficiently high accuracy. In the same paper, it is argued that problems with non-relatively complete recourse is even more challenging to solve \cite[Theorem 4]{hanasusanto2016comment}. Intuitively, this is caused by additional implicit constraints at each time-step to ensure feasibility in the future. One common approach to dealt with this issue is to incorporate a high cost for any violation of the hard constraints. Quite often this can transform the program into an equivalent problem with relatively complete recourse, but instead the reward function will have a huge Lipschitz constant \cite[Section 3]{shapiro2005complexity}. To avoid technical difficulties we hence opt to assume that our problem has relatively complete recourse. 
  
  For a general $T$-stage program with stagewise-independent and decision-dependent noise, it is shown in \cite{dyer2006computational} that the problem is PSPACE-hard, by treating $T$ as an input parameter. The statue for the complexity of the simpler problem with decision-independent noise is still open (e.g. with only the $\bW$-part in our model \eqref{eq:original-problem}). It is conjectured in \cite{dyer2006computational} that it is PSPACE-hard as well.  \Halmos
  
\end{remark}

\section{Concrete Examples}  \label{sec:concrete-examples}

The problem \eqref{eq:original-problem} covers several classes of models widely studied in the literature as particular cases. We list three below: The first requires the full modeling generality of \eqref{eq:original-problem} and will be studied numerically later in Section \ref{sec:numerical-experiments}; the two later are with affine rewards and polyhedron constraint sets, originate respectively from inventory management and weakly coupled Markov decision processes. 

\subsection{Network Resource Allocation and Utility Maximization}  \label{example-3}

In a network utility maximization model (\cite{shakkottai2008network, srikant2004mathematics, palomar2006tutorial}), $\bX(t)$ represents the bandwidth occupation of $n_x$ routing paths in a communication network with a known topology, $\bW(t)$ represents the arrival of service demands among the $n_x$ paths, and $\bU(t)$ represents the amount of bandwidth that we allocate to each path. Each of the $n_x$ paths use a certain collection of links in the network, and \eqref{eq:constraints}-\eqref{eq:constraints-eq} refers to the constraints on the link capacities, as well as the delays suffered on each routing path, which are typically non-linear convex functions. \eqref{eq:evolution-of-the-system} describes the stochastic evolution of the bandwidth dynamics. Our goal is to maximize a reward function of transferring data flows via the network. One typical choice is $R_t(\bX) = \sum_{i=1}^{n_x} X_i^{1-\alpha}/(1-\alpha)$ for $\alpha > 0$, called the $\alpha$-fairness utility (\cite{mo2000fair}). This is a dynamic extension of the classical network utility maximization problem. 

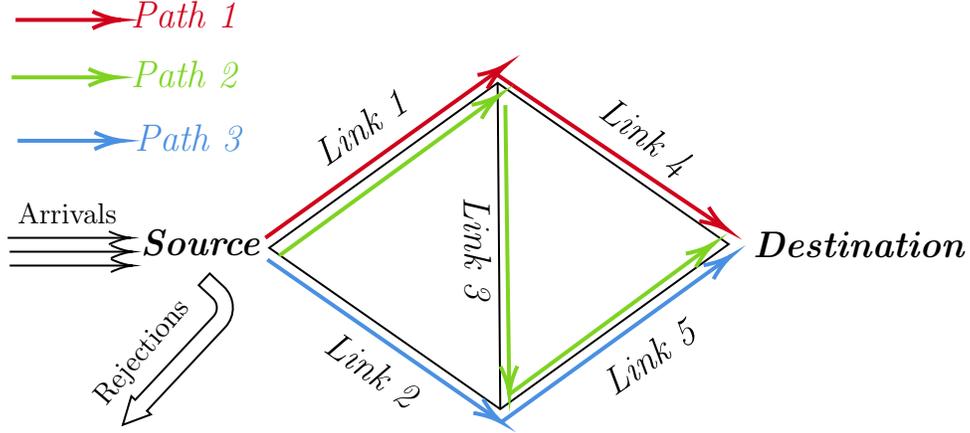
\begin{figure}[htbp]
  \centering
  \input{figure/network4.tex}
  \caption{An example of diamond-shaped network with $3$ paths and $5$ links.}\label{fig-network-example}
\end{figure}

Typically, consider a diamond-shaped routing network that can be represented as a directed graph in Figure \ref{fig-network-example}. Each directed edge of the graph is called a \emph{link}, enumerated by $1 \le l \le 5$. Each link has a maximum bandwidth capacity, denoted as $c_l > 0$. At each discrete time-step $1 \le t \le T$, a certain quantity of demand of bandwidth arrives at the source, which will occupy a specific \emph{path} to reach its destination. A path is a sequence of consecutive directed links that connect the source to the destination, enumerated by $1 \le p \le 3$. We assume that the allocation of bandwidth is instantaneous, meaning that it immediately occupies all the subsequent links of this path. Denote by $(X_1(t),X_2(t),X_3(t))$ the bandwidth occupation of the three routing paths just arriving at time-step $t$, by $(W_1(t),W_2(t),W_3(t))$ the arrivals of new demands of bandwidth on each path at time-step $t$, and the decision is the allocation $(U_1(t),U_2(t),U_3(t))$ of bandwidth among each path. We suppose that the demand is \emph{elastic}, so that the non-satisfied demand incurs no cost. A first set of constraints on the model can then be expressed as 
\begin{equation}\label{eq:network-model-0}
   0 \le U_1(t) \le W_1(t), \ \ \  0 \le U_2(t) \le W_2(t), \ \ \ 0 \le U_3(t) \le W_3(t) \ \ \mbox{ for $1 \le t \le T$} 
\end{equation}
and
\begin{align}
 & Y_1(t)\defeq U_1(t) + X_1(t) + U_2(t) + X_2(t) \le c_1, \ \ \  Y_2(t) \defeq U_3(t) + X_3(t) \le c_2  \ \ Y_3(t) \defeq U_2(t) + X_2(t) \le c_3 \nonumber  \\
 & Y_4(t) \defeq U_1(t) + X_1(t) \le c_4, \ \ \ Y_5(t) \defeq U_2(t) + X_2(t) + U_3(t) + X_3(t) \le c_5  \mbox{ for $1 \le t \le T$}  \label{eq:network-model-1}
\end{align}
To model the stochasticity, we introduce the notation $\calTN(\mu,\sigma^2,a,b)$ to denote a \emph{truncated normal distribution}, which represents a normal distribution $\calN(\mu,\sigma^2)$ truncated at the interval $[a,b]$. We assume that for $1 \le p \le 3$ and $1 \le t \le T$, the occupation of bandwidth of path $p$ follows the dynamic 
\begin{align} \label{eq:network-model-2}
 & X_p(t+1) = (X_p(t) + U_p(t)) \cdot q_p  \nonumber \\ 
 & + \calTN \left( 0, \sigma^2, -(X_p(t) + U_p(t))\cdot \min \{q_p, 1-q_p \}, (X_p(t) + U_p(t)) \cdot \min \{q_p, 1-q_p \} \right) 
\end{align}
with $0 < q_p < 1$ and $\sigma > 0$. This modeling approach seeks to replicate the scenario in which the bandwidth of a path $p$ is consumed by numerous discrete flows. Each of these flows has a lifespan, which we model using a geometric distribution with parameter $q_p$. After its lifespan ends, the flow releases its occupation of the links corresponding to its path. These flows operate independently of one another. Equation \eqref{eq:network-model-2} can be interpreted as a continuous approximation of this dynamic. For simplicity, we have chosen to truncate on a symmetric interval around $0$, ensuring that the truncated normal distribution has a mean of zero. In a similar vein, we suppose that the arrivals of new demands follow the distribution
\begin{equation} \label{eq:network-model-3}
  W_p(t) = \bar{w}_p + \calTN \left(0,\sigma^2, -\bar{w}_p, \bar{w}_p \right) \mbox{ for $1 \le p \le 3$ and $1 \le t \le T$}
\end{equation}
with $\bar{w}_p > 0$. Note that all variances of the stochastic uncertainties in the model are governed by a single parameter $\sigma$, which provides convenience for our later numerical study.

Following \cite{parag2011value}, we next introduce the \emph{link quality degradation function} $d_l(y)$ for a link $l$ with a occupation bandwidth $0 \le y < c_l$ as
\begin{equation} \label{eq:network-model-4}
  d_l(y) \defeq \frac{y}{c_l(c_l - y)} = \frac{1}{c_l - y} - \frac{1}{c_l} \ \ \mbox{ for $1 \le l \le 5$ and $0 \le y < c_l$}
\end{equation}
Note that $d_l(\cdot)$ is non-negative, convex, and increases from $0$ to $\infty$ when $y$ ranges in $[0, c_l)$. Other choices of $d_l(y)$, as well as its practical implication are discussed in detail in \cite{parag2011value}. Note that \eqref{eq:network-model-4} is motivated by seeing link $l$ as a $M/M/1$ queue with arrival rate $y$ and service rate $c_l$, since then $d_l(y)$ is the expected waiting time in the queue. 

We assume that the total degradation experienced along path $p$ is the sum of degradation of links that it traverses, and there is a pre-specified positive value $D_p$ being the maximum degradation that can be tolerated on path $p$ at each time-step. This gives rise to a second set of constraints on the model as
\begin{align}  \label{eq:network-model-5}
  & d_1(Y_1(t)) + d_4(Y_4(t)) \le D_1, \ \ \ d_1(Y_1(t)) + d_3(Y_3(t)) + d_5(Y_5(t)) \le D_2, \ \ \ d_2(Y_2(t)) + D_5(Y_5(t)) \le D_3  \nonumber \\
  & \mbox{ for $1 \le t \le T$} 
\end{align}
The decision-maker aims to maximize the following $\alpha$-fairness utility (with $\alpha > 0$)
\begin{equation}  \label{eq:network-model-6}
  \expect{\sum_{t=1}^{T} \sum_{p=1}^{3} \frac{(X_p(t)+U_p(t))^{1 - \alpha}}{1 - \alpha} }
\end{equation}
gained by allocating the bandwidth demands over a finite horizon $T$, while respecting the dynamics and constraints described in \eqref{eq:network-model-0}-\eqref{eq:network-model-5}. 

We suppose that initially the network is in a state $\bX(1) \defeq \bx$ such that if we take $\bU(1) \equiv \bzero$, all constraints in \eqref{eq:network-model-0}, \eqref{eq:network-model-1} and \eqref{eq:network-model-5} are satisfied with strict inequality. We then remark that by always taking $U_p(t) \equiv \epsilon$ with a fixed $\epsilon > 0$ small enough for all subsequent time-steps $t$, all inequality constraints in \eqref{eq:network-model-0}-\eqref{eq:network-model-5} can be satisfied with strict inequality, so Assumption \ref{assumption-7} and Slater CQ of Assumption \ref{assumption-8} are met for this model. Also by our choice of the truncated normal distributions, other conditions in Assumptions \ref{assumption-1}-\ref{assumption-8} concerning the stochastic part of the model are also satisfied. 

We shall study numerically one such example using the heuristics developed in this paper in Section \ref{sec:numerical-experiments}.   

\subsection{Network Revenue and Inventory Management}   \label{example-1}

In an inventory management scenario, at each time-step $t$, $\bX(t)$ represents the amount of $n_x$ types of resource available in storage, $\bW(t)$ represents the quantity of arrivals of $n_w$ types of (real-valued) customer demands, and $\bU(t)$ consists of two parts: the first part $\bU_1(t)$ represents the amount of additional resource $j$, for $1 \le j \le n_x$, ordered to replenish the storage, at a cost of $c_j$; the second part $\bU_2(t)$ represents the amount of each type $i$, for $1 \le i \le n_w$, of customer demands that we actually serve at time-step $t$, hence $\bU_2(t) \le \bW(t)$. Each type $i$ of customer demand, for $1 \le i \le n_w$, requires a certain combination $a_{ij}$ units of resources $j$ to be served, for $1 \le j \le n_x$, and produce the amount $r_i$ of profit to the decision-maker. While the customer demands that are not served leave the system and are considered as lost with no extra cost. In addition, the resources that are left in the storage incur a holding cost $h_j$ per unit of time, for $1 \le j \le n_x$. So to summarize the reward can be written as
  \begin{equation*}
    R_t \left( \bX(t), \bW(t), \bU(t) \right) = - \bU_1(t) \cdot \mathbf{c}^\top + \bU_2(t) \cdot \mathbf{r}^\top - \left( \bX(t) + \bU_1(t) - \bU_2(t) \cdot \mathbf{A} \right) \cdot \mathbf{h}^\top 
  \end{equation*}  
where $\mathbf{A}$ is the matrix with entries $a_{ij}$. There is clearly also the constraint $\bX(t) + \bU_1(t) - \bU_2(t) \cdot \mathbf{A} \ge \bzero$. Hence the optimization problem \eqref{eq:original-problem} for this inventory management problem can be formulated as
  \begin{maxi!}|s|{\bU[1,T]}{\expect{\sum_{t=1}^{T} - \bU_1(t) \cdot \mathbf{c}^\top + \bU_2(t) \cdot \mathbf{r}^\top - \left( \bX(t) + \bU_1(t) - \bU_2(t) \cdot \mathbf{A} \right) \cdot \mathbf{h}^\top } }{ \label{eq:inventory} }{V_{\mathrm{opt}} (\bx(1),T) =}
  \addConstraint{ \bX(1) = \bx(1)  \ a.s. \nonumber }
  \addConstraint{ \bU_2(t) - \bW(t) \le \bzero  \ a.s. \ \ \mbox{ for $ 1 \le t \le T$ } \nonumber }{}
  \addConstraint{ - \bX(t)  \le \bzero  \ a.s. \ \ \mbox{ for $ 2 \le t \le T$ } \nonumber }{}
  \addConstraint{ - \bU_1(t) \le \bzero, \ - \bU_2(t) \le \bzero  \ a.s. \ \ \mbox{ for $1 \le t \le T$ } \nonumber }{}
  \addConstraint{ \bX(t+1) = \bX(t) + \bU_1(t) - \bU_2(t) \cdot \mathbf{A} \ a.s. \ \ \mbox{ for } 1 \le t \le T-1 \nonumber }
\end{maxi!}
which is a particular case with linear rewards, linear constraints and without the $\calE$ stochastic part. This is the model that has been considered in \cite[Section 1.3.3]{shapiro2021lectures} as a $T$-stage stochastic optimization problem. 
  
Many of its variants have been studied extensively in the literature. Notably, the extension of the model with customer-choice has been considered in \cite{kunnumkal2008refined, bront2009column, jasin2012re}, which incorporate randomness from customer demands' resource consumptions, and can be modeled using the $\calE$-part: $\bU_2(t)$ is hence replaced everywhere by $\bU_2(t) + \calE(\bU_2(t))$, where $\calE(\bU_2(t))$ are noises with distributions depend on $\bU_2(t)$. A simpler model without the possibility of replenishing the storage has been considered in \cite{cooper2002asymptotic, secomandi2008analysis, bumpensanti2020re}, using a re-solving policy similar to the update policy considered in the current work. A thorough discussion for inventory models can be found in \cite{talluri2004theory,zipkin2000foundations}.  

\subsection{Markovian Bandits and Weakly Coupled Markov Decision Processes}   \label{example-2}
 
A bandit consists of an integer number of symmetric arms, where each arm has $n_x$ states, that evolves as a Markov decision process (MDP). The whole bandit itself can be seen as a large MDP that consists of $N$ statistically identical sub-MDPs (arms), with $N \in \mathbb{N}$ being an integer. The vector $\bX(t)$ represents the number of arms of the bandit being in each of the $n_x$ states, that sums to $N$. An arm can undertake one among a certain number $n_u$ of actions, where each action $j$, for $1 \le j \le n_u$, consumes a certain combination $c_{jk}$, for $1 \le k \le I(t)$ of resources, with a total availability of $\mathbf{b}(t) \in \R^{I(t)}$ resources. In the meantime, an arm in state $i$, for $1 \le i \le n_x$, undertaking action $j$, for $1 \le j \le n_u$, earns a reward $r_{ij}$ for the decision-maker. The vector $\bU(t)$ represents the number of arms undertaking each action, hence it will be more convenient to treat it as a matrix of dimension $n_x \times n_u$, and naturally $\bU(t) \cdot \bonet = \bX(t)$. The resource constraints can be written compactly as $ \bone \cdot \bU(t) \cdot \mathbf{c} \le \mathbf{b}(t)$, with $\mathbf{c}$ being the matrix of entries $c_{jk}$; while the utility collected at time-step $t$ is $\bone \cdot \bU(t) \cdot \mathbf{r}^\top \cdot \bonet$, with $\mathbf{r}$ being the matrix of entries $r_{ij}$. Denote by "$\mathrm{Multinomial}(n,\bp)$" the multinomial distribution with $n$ trials and probability vector of success $\bp$.  For $1 \le i \le n_x$ and $1 \le j \le n_u$, the Markov evolution of $\bX(t)$ can be written as 
\begin{align*}
  \bX(t+1) & = \sum_{i=1}^{n_x} \sum_{j=1}^{n_u} \mathrm{Multinomial}(U_{ij}(t), \bP^{(j)}_i ) \\
  & =  \underbrace{\sum_{j=1}^{n_u} \sum_{j=1}^{n_u} U_{ij}(t) \cdot \bP^{(j)}_i}_{\phi(\bU(t))} + \underbrace{\sum_{i=1}^{n_x} \sum_{j=1}^{n_u} \mathrm{Multinomial}(U_{ij}(t), \bP^{(j)}_i ) - \sum_{j=1}^{n_u} \sum_{j=1}^{n_u} U_{ij}(t) \cdot \bP^{(j)}_i}_{\calE(\bU(t))}   
\end{align*}
where $\bP^{(j)}$ for $1 \le j \le n_u$, are in total $n_u$ transition probability matrices of size $n_x \times n_x$, and $\calE(\bU(t))$ are decision-dependent noise obtained by substituting from the multinomial distributions their mean values $\phi(\bU(t))$, the latter being an affine function. So in summary the optimization problem for the Markovian bandit model is 
\begin{maxi!}|s|{\bU[1,T]}{\expect{\sum_{t=1}^{T} \bone \cdot \bU(t) \cdot \mathbf{r}^\top \cdot \bonet }}{ \label{eq:WCMDP-1} }{V_{\mathrm{opt}} (\bx(1),T) =}
  \addConstraint{ \bX(1) = \bx(1) \ a.s. \nonumber }
  \addConstraint{ \bone \cdot \bU(t) \cdot \mathbf{c} - \mathbf{b}(t) \le \bzero  \ a.s. \ \ \mbox{ for $ 1 \le t \le T$ } \nonumber }{}
  \addConstraint{ - \bU(t) \le \bzero  \ a.s. \ \ \mbox{ for $1 \le t \le T$ } \nonumber }{}
  \addConstraint{ - \bU(t) \cdot \bonet + \bX(t) = \bzero  \ a.s. \ \ \mbox{ for $ 1 \le t \le T$ } \nonumber }{}
  \addConstraint{ \bX(t+1) = \sum_{i=1}^{n_x} \sum_{j=1}^{n_u} \mathrm{Multinomial}(U_{ij}(t), \bP^{(j)}_i ) \ \ a.s. \ \mbox{ for $1 \le t \le T-1$}  \label{eq:WCMDP-2} }
\end{maxi!}
which has linear rewards, linear constraints and without the $\bW$ stochastic part. Note that in the current model $\bX(t)$ and $\bU(t)$ are discrete and integer-valued, and does not fit into the continuous-valued setting we supposed in problem \eqref{eq:original-problem}. We refer to Remark \ref{rem:scaling-under-affinity} for a further discussion on this issue.   

This finite horizon Markovian bandit model has been widely studied in the literature, see e.g. \cite{hu2017asymptotically,ZayasCabn2017AnAO,Brown2020IndexPA,zhang2021restless,gast2021lp}, as well as its generalization to weakly coupled Markov decision processes \cite{adelman2008relaxations, carpentier2020mixed, gast2022lp, brown2022strength, brown2022fluid}. In \cite{Papadimitriou99thecomplexity}, a problem of routing and scheduling in closed queueing networks, called "Network of Queues" has been studied. This problem can be reformulated into the form of \eqref{eq:WCMDP-1}, by considering each class of jobs as a state, and each server as an action, with an additional action of being idle. The constraints correspond to the classes of jobs that a particular server can serve. By treating $N$ as an input parameter, the problem is shown to be EXP-complete, provided that the horizon $T$ is exponential in $N$, see \cite[Corollary 1]{Papadimitriou99thecomplexity} and \cite[Section 5.2]{blondel2000survey}. The simpler "Restless Bandits" problem that has a single server is proven to be PSPACE-complete \cite[Theorem 4]{Papadimitriou99thecomplexity}.

\section{The Certainty Equivalent Control (CEC)} \label{sec:CEC}

The heuristic policies to problem \eqref{eq:original-problem} that we shall discuss in this work, i.e. the update policy in Section \ref{subsec:re-solving-heuristic}, the projection policy in Section \ref{subsec:policy-projection}, are inspired from the \emph{certainty equivalent control} (CEC). The CEC is in general a sub-optimal control that applies at each stage the control that would be optimal if some or all of the uncertain quantities were fixed at their expected values \cite[Chapter 6]{bertsekas2012dynamic}. There is an additional difficulty in the current situation to apply CEC, however, as we are facing a problem with \emph{hard constraints} that depend on both the current system state $\bX(t)$ and disturbance $\bW(t)$, the \emph{feasibility} of an action is hence of major concern. A key point in the development of the subsequent sections is centered around how to design feasible controls after taking the expectation.

Based on CEC, we apply the following relaxation to the original problem \eqref{eq:original-problem}: define $\expect{\bX(t)} \defeq \bx(t)$ and $\expect{\bU(t)} \defeq \bu(t)$ where the expectation is taken with the whole trajectory. From the convexity requirement made in Assumption \ref{assumption-1}, these are well-defined system states and controls. By Assumptions \ref{assumption-2}, \ref{assumption-3}, \ref{assumption-4} and Jensen's inequality, we have 
\begin{equation} \label{eq:jensen-concave}
  \expect{\sum_{t=1}^{T} R_t \left( \bX(t), \bW(t), \bU(t) \right)} \le \sum_{t=1}^{T} R_t \left( \bx(t), \barw, \bu(t) \right).
\end{equation}
\begin{equation}\label{eq:jensen-convex}
  \expect{g_{t,i}(\bX(t), \bW(t), \bU(t))} \ge g_{t,i}(\bx(t), \barw, \bu(t)) \ \ \mbox{ for $ 1 \le t \le T$ and $ 1 \le i \le I(t)$}
\end{equation}
\begin{equation}\label{eq:jensen-affine}
  \expect{h_{t,j}(\bX(t), \bW(t), \bU(t))} = h_{t,j}(\bx(t), \barw, \bu(t)) \ \ \mbox{ for $ 1 \le t \le T$ and $ 1 \le j \le J(t)$}
\end{equation}
By Assumption \ref{assumption-5}, the expectations can be interchanged with the affine function $\phi(\cdot)$ in \eqref{eq:evolution-of-the-system}. By Assumption \ref{assumption-6}, the stochastic parts $\calE(\bX(t), \bW(t), \bU(t))$ are with zero means. All these considerations lead to the following relaxed mathematical program with decision variables $\bu[1,T]$:

\begin{maxi!}|s|{\bu[1,T]}{\sum_{t=1}^{T} R_t \left( \bx(t), \barw, \bu(t) \right)   \label{eq:cost-overall-rel} }{\label{eq:original-problem-rel}}{V_{\mathrm{rel}-} (\bx,T) =}
  \addConstraint{ \bx(1) = \bx  \label{eq:initial-configuration-rel} }
  \addConstraint{g_{t,i}(\bx(t), \barw, \bu(t)) \le 0 \ \ \mbox{ for $ 1 \le t \le T$ and $ 1 \le i \le I(t)$} \label{eq:constraints-rel} }{}
  \addConstraint{h_{t,j}(\bx(t), \barw, \bu(t)) = 0 \ \ \mbox{ for $ 1 \le t \le T$ and $ 1 \le j \le J(t)$} \label{eq:constraints-rel-eq} }{}
  \addConstraint{\bx(t+1) = \phi \left( \bx(t), \barw, \bu(t) \right)  \ \ \ \mbox{ for } 1 \le t \le T-1  \label{eq:evolution-of-the-system-rel}}
\end{maxi!}
The subscript "$_{\mathrm{vel}-}$" in the notation $V_{\mathrm{rel}-}$ is to make contrast with the later notation $V_{\mathrm{rel}+}$, emphasizing the fact that the relaxed program \eqref{eq:original-problem-rel} is solved without the knowledge of $\bW(1)$ at time-step $1$. More generally, for each decision epoch $1 \le t \le T$, we write $V_{\mathrm{rel}-} (\bx,T+1-t)$ to refer to the relaxed CEC problem at time-step $t$ with current system state $\bx(t) = \bx$ and time-span $[t,T]$. Remark that \eqref{eq:original-problem-rel} is a deterministic convex program. From \eqref{eq:jensen-concave} and \eqref{eq:jensen-convex}, we see that in \eqref{eq:original-problem-rel} we are maximizing a larger objective function over a less restrictive feasible region. We hence obtain the key inequality $V_{\mathrm{opt}} (\bx,T) \le V_{\mathrm{rel}-} (\bx,T)$.

\begin{remark} [Necessity of the Convexity Assumption]  \label{rem:necessitiy-of-the-assumtions}

The convexity assumptions made in model \eqref{eq:original-problem} merit further justification. Specifically, the concavity of the reward functions $R_t(\cdot)$, the convexity of the inequality constraint functions $g_{t,i}(\cdot)$, and the affinity of the equality constraint functions $h_{t,j}(\cdot)$ as well as the system evolution function $\phi(\cdot)$, together ensure the relaxed problem \eqref{eq:original-problem-rel} following the CEC remains a convex program. From a computational point of view, there exist solution methods that efficiently solve every convex optimization program satisfying very mild computability restrictions; in contrast, no efficient universal solution methods for non-convex programs are known, and there are strong reasons to expect that no such methods exist, see \cite[Lecture 5]{ben2001lectures} for a thorough discussion. So it is reasonable to remain in the convex paradigm, as otherwise even the relaxed problem may not be tractable.  

More importantly, these properties are pivotal as they guarantee that in \eqref{eq:original-problem-rel} we are maximizing a larger objective function over a superset of the original feasible region, securing the vital inequality $V_{\mathrm{opt}} (\bx,T) \le V_{\mathrm{rel}-} (\bx,T)$. The significance of maintaining $V_{\mathrm{opt}} (\bx,T) \le V_{\mathrm{rel}-} (\bx,T)$ lies in providing an upper-bound estimate for the \emph{sub-optimality gap} of any heuristic policy under consideration. This relationship is encapsulated in the inequality:
\begin{align}  \label{eq:pivotal-relation}
 \mbox{sub-optimality gap} & \defeq V_{\mathrm{opt}} (\bx,T) - \mbox{value of an heuristic policy} \nonumber \\
 & \ \le V_{\mathrm{rel}-} (\bx,T) - \mbox{value of an heuristic policy}
\end{align}
As previously noted in Remark \ref{rem:model-assumptions}, determining an exact value for $V_{\mathrm{opt}} (\bx,T)$ to evaluate the sub-optimality gap is generally a complex task. However, $V_{\mathrm{rel}-}(\bx,T)$ can be acquired by solving a convex program, and the value of an efficient heuristic policy can be estimated via Monte-Carlo simulation. Thus, for the remainder of this paper, we will use the right-hand side of \eqref{eq:pivotal-relation} as an upper-bound estimate of the sub-optimality gap, termed an \emph{optimality gap bound}.  \Halmos

\end{remark}

\section{Heuristic Policies based on CEC}  \label{sec:main-performance-bound}

This section contains the algorithmic results of this paper. Having introduced the CEC in Section \ref{sec:CEC}, we construct two heuristic policies based on this principle, the first given in Section \ref{subsec:re-solving-heuristic} is called the update policy, which requires re-solving a new relaxed mathematical program at each time-step. The second given in Section \ref{subsec:policy-projection} is called the projection policy, which solves a single program at the start, and simpler Euclidean projection problems at each later time-step. A hybrid policy that combines the advantages of both of these two policies is then introduced in Section \ref{subsec:hybrid-policy}.

\subsection{The Update Policy with Re-Solving}   \label{subsec:re-solving-heuristic}

Denote by $\bu^*_{\bx,T}$ an optimal solution of \eqref{eq:original-problem-rel}, which exists by Proposition \ref{prop:existence-of-optimal-solution}. We use the subscript "$_{\bx,T}$" to keep track that the solution is with respect to initial system configuration $\bx$ and horizon length $T$. The observation is that the first control $\bu^*_{\bx,T}(1)$ from this solution is in general not feasible to the original problem \eqref{eq:original-problem}, due to the fact that the realization of $\bW(1)$ is still unknown at that moment, and the planning is only guaranteed to be feasible were it be that $\bW(1) = \barw$. The key is that we should apply the re-solving \emph{after} the realization of $\bW(1)$ is known. 

More generally, for each decision epoch $1 \le t \le T$, just before taking a decision, the decision-maker observes the realization $\bw$ of $\bW(t)$ as well as the system state $\bx$, and then solves the following relaxed problem:

\begin{maxi!}|s|{\bu[t,T]}{R_{t} \left( \bx(t), \bw, \bu(t) \right) + \sum_{t'= t+1}^{T} R_{t'} \left( \bx(t'), \barw, \bu(t') \right)   \label{eq:cost-overall-res} }{\label{eq:original-problem-res}}{\hat{V}_{\mathrm{rel}+} (\bx,T+1-t,\bw) =}
  \addConstraint{ \bx(t) = \bx  \label{eq:initial-configuration-res} }
  \addConstraint{g_{t,i}(\bx(t), \bw, \bu(t)) \le 0 \ \ \mbox{ for $ 1 \le i \le I(t)$}  \label{eq:constraints-res-1} }{}
  \addConstraint{h_{t,j}(\bx(t), \bw, \bu(t)) = 0 \ \ \mbox{ for $ 1 \le j \le J(t)$}  \label{eq:constraints-res-1-eq} }{}
  \addConstraint{\bx(t+1) = \phi \left( \bx(t), \bw, \bu(t) \right)  \label{eq:evolution-of-the-system-res-1}}
  \addConstraint{g_{t',i}(\bx(t'), \barw, \bu(t')) \le 0 \ \ \mbox{ for $ t+1 \le t' \le T$ and $ 1 \le i \le I(t')$} \label{eq:constraints-res-2} }{}
  \addConstraint{h_{t',j}(\bx(t'), \barw, \bu(t')) \le 0 \ \ \mbox{ for $ t+1 \le t' \le T$ and $1 \le j \le J(t')$} \label{eq:constraints-res-2-eq} }{}
  \addConstraint{\bx(t'+1) = \phi \left( \bx(t'), \barw, \bu(t') \right)  \ \ \ \mbox{ for } t+1 \le t' \le T-1  \label{eq:evolution-of-the-system-res-2}}
\end{maxi!}
Note that the only difference between $V_{\mathrm{rel}-} (\bx,T+1-t)$ and $\hat{V}_{\mathrm{rel}+} (\bx,T+1-t,\bw)$ is that the latter has taken the information $\bW(t) = \bw$ into account, and we actually have 
\begin{equation*}
  V_{\mathrm{rel}-} (\bx,T+1-t) = \hat{V}_{\mathrm{rel}+} (\bx,T+1-t,\barw)
\end{equation*}
Denote by $\bu^*_{\bx,T+1-t,\bw}[t,T]$ an optimal solution of \eqref{eq:original-problem-res}, which exists by Proposition \ref{prop:existence-of-optimal-solution}. We use the subscript "$_{\bx,T+1-t,\bw}$" to keep track that the solution is with respect to system configuration $\bx$, horizon $T+1-t$ and upon observation of $\bW(t)=\bw$. The first control $\bu^*_{\bx,T+1-t,\bw}(t)$ from this control sequence is by construction feasible to \eqref{eq:original-problem}. Note that we write $\bu^*_{\bx,T+1-t,\bw}(t)$ rather than $\bu^*_{\bx,T+1-t,\bw}(1)$, since although it is the first control to problem \eqref{eq:original-problem-res}, it is the $t$-th control to the original problem \eqref{eq:original-problem}. We repeat this process at each decision epoch $1 \le t \le T$, and this gives rise to the update policy described in Algorithm \ref{algo:update}. 

%Also see Figure \ref{fig:evolution-of-system} for an illustration.

%\begin{figure}[htbp]
%  \centering
%  \input{figure/evolution.tex}
%  \caption{Illustration of the evolution of the process under the update policy. The "cloud" shape is used to represent the stochasticity. Note that there is no "cloud" near $\bx^*(1) = \bX(1)$. Sample paths reaching outside the $\varepsilon$-neighbourhood of the deterministic trajectory are considered "rare", and their analysis are disregarded, with a value estimation given by the trivial bound $\bar{V}$. For instance, one sample path with the realization of $\bW'(2)$, and another sample path with the realization of $\bX'(3)$ (represented in orange color) are considered rare events.  }\label{fig:evolution-of-system}
%\end{figure}

\begin{algorithm}[h]
  \SetAlgoLined
  \SetKwInput{KwInput}{Input}
  \KwInput{Initial system configuration vector $\bx(1)$ and time horizon $T$.}
  Set $\bx \defeq$  current system configuration vector \;
  \For{$t = 1,2,\dots,T$}{
    Observe the realization $\bw$ of $\bW(t)$  \; 
    Solve $\hat{V}_{\mathrm{rel}+} (\bx,T+1-t,\bw)$ and obtain the control sequence $\bu^*_{\bx,T+1-t,\bw}[t,T]$  \;
    Use control $\bu^*_{\bx,T+1-t,\bw}(t)$ to advance to the next time-step \;
    Set $\bx \defeq$  current system configuration vector  \;
  }
\caption{The Update Policy.}
  \label{algo:update}
\end{algorithm}
We use "\emph{update}" as a short hand for the update policy defined in Algorithm \ref{algo:update}, and denote by $V_{\mathrm{update}} (\bx,T)$ the value of the update policy with initial condition $\bx$ and horizon $T$.

\subsection{The Projection Policy without Re-Solving}  \label{subsec:policy-projection}

The update policy in Algorithm \ref{algo:update} requires re-solving a convex program that spans the horizon $[t,T]$ at each time-step $t$. This re-solving procedure adjusts the decision based on current available information, and may help to avoid the accumulation of stochastic noises. Yet, from a computational point of view it is not efficient. In this section we propose the idea of policy projection, that replaces the task of re-solving a complex convex program by computing a simpler Euclidean projection.   

As before, by solving \eqref{eq:original-problem-rel} we obtain an optimal solution $\bu^*[1,T]$ with the corresponding $\bx^*[1,T]$. For each time-step $t$, after the realization of $\bX(t) = \bx(t)$ and $\bW(t) = \bw$, the control $\bu^*(t)$ is in general not feasible, as it was planned with respect to $\bW(t) = \barw$ and $\bX(t) = \bx^*(t)$. In other words, recall the following notation for the set of feasible controls at time-step $t$:
\begin{equation}  \label{eq:feasible-set-for-projection}
  \calU_t(\bx,\bw) = \left\{ \bu \in \R^{n_u} \ \big| \ g_{t,i} (\bx,\bw,\bu) \le 0, \ \ \mbox{ for }  1 \le i \le I(t); \ \ h_{t,j} (\bx,\bw,\bu) = 0, \ \ \mbox{ for } 1 \le j \le J(t) \right\}
\end{equation} 
then $\bu^*(t) \in \calU_t (\bx^*(t), \barw)$ but in general $\bu^*(t) \notin \calU_t (\bx(t), \bw)$. The idea is that since $\calU_t (\bx(t), \bw)$ is the set of all feasible actions at time-step $t$, which is a non-empty closed set by our assumption, we may apply the Euclidean projection of $\bu^*(t)$ onto $\calU_t (\bx(t), \bw)$ to obtain a feasible action, which is the closest feasible action to $\bu^*(t)$ measured by Euclidean distance. Denote by $\Pi_{\calU_t (\bx(t), \bw)}(\bu^*(t))$ this Euclidean projection. It is the unique solution to the following convex program \emph{parameterized} by $(\bx(t),\bw)$:
\begin{equation}  \label{eq:convex-program-for-projection}
  \min_{\bu} \frac{1}{2} (\bu - \bu^*(t)) \cdot (\bu - \bu^*(t))^\top \mbox{ subject to } \bu \in \calU_t(\bx(t),\bw) \mbox{ given by \eqref{eq:feasible-set-for-projection} }
\end{equation} 
We apply this projected action to the system, and repeat this procedure at each time-step $1 \le t \le T$. This gives rise to the projection policy summarized in Algorithm \ref{algo:projection}.

\begin{algorithm}[h]
  \SetAlgoLined
  \SetKwInput{KwInput}{Input}
  \KwInput{Initial system configuration vector $\bx(1)$ and time horizon $T$.}
  Solve \eqref{eq:original-problem-rel} for time-span $[1,T]$ to obtain an optimal solution $\bu^*[1,T]$ \;
  Set $\bx \defeq $  current system configuration vector \;
  \For{$t = 1,2,\dots,T$}{
    Observe the realization $\bw$ of $\bW(t)$  \; 
    Compute the Euclidean projection $\bu_{\pi}(t) = \Pi_{\calU_t(\bx, \bw)}(\bu^*(t))$  \;
    Use control $\bu_{\pi}(t)$ to advance to the next time-step \;
    Set $\bx \defeq$  current system configuration vector  \;
  }
\caption{The Projection Policy.}
  \label{algo:projection}
\end{algorithm}
We use "\emph{proj}" as a short hand for the projection policy defined in Algorithm \ref{algo:projection}, and denote by $V_{\mathrm{proj}} (\bx,T)$ the value of the projection policy with initial condition $\bx$ and horizon $T$.

\subsection{A Hybrid Policy}  \label{subsec:hybrid-policy} 

From our previous discussion, we see that the update policy and the projection policy each has its advantage and disadvantage, which are complementary to each other: the update policy is more robust against stochastic uncertainties while requires significantly more computational resource; the projection policy on the other hand, is time-efficient but may not prevent the propagation of estimation errors. 

\begin{algorithm}[h]
  \SetAlgoLined
  \SetKwInput{KwInput}{Input}
  \KwInput{Initial system configuration vector $\bx(1)$ and time horizon $T$. A threshold parameter $\Theta > 0$. }
  Solve \eqref{eq:original-problem-rel} for time-span $[1,T]$ to obtain an optimal solution $\bu^*[1,T]$ \;
  Set the current deterministically optimal control as $\mathfrak{U} \defeq \bu^*[1,T]$, and the deterministically optimal system trajectory as $\mathfrak{X} \defeq \bx^*[1,T]$ \;
  Set $\bx \defeq$  current system configuration vector \;
  \For{$t = 1,2,\dots,T$}{
    Observe the realization $\bw$ of $\bW(t)$  \; 
    Compute the Euclidean projection $\bu_{\pi} \defeq \Pi_{\calU_t(\bx, \bw)}(\mathfrak{U}(t))$ \;
    Compute the Euclidean distance $\theta \defeq \norme{(\mathfrak{X}(t), \barw, \mathfrak{U}(t)) - (\bx,\bw,\bu_{\pi}) }$ \;
    \If{$\theta < \Theta$}{
    Use control $\bu_{\pi}$ to advance to the next time-step \;
    Set $\bx \defeq$ current system configuration vector  \;
    }
    \Else{
    Solve $\hat{V}_{\mathrm{rel}+} (\bx,T+1-t,\bw)$ in \eqref{eq:original-problem-res} for time-span $[t,T]$ \;
    Update $\mathfrak{U} \defeq \hat{\bu}^*[t,T]$ and $\mathfrak{X} \defeq \hat{\bx}^*[t,T]$ from its solution  \;
    Use control $\mathfrak{U}(t)$, which is feasible by construction, to advance to the next time-step \;
    Set $\bx \defeq$ current system configuration vector  \;
    }
  }
\caption{A hybrid policy that combines the projection policy and the update policy.}
  \label{algo:hybrid}
\end{algorithm}

Motivated by this observation, we propose a hybrid algorithm that combines the strength of both policies, which works as follows: we fix a tuning parameter $\Theta > 0$ and keep in the memory a deterministically optimal control $\mathfrak{U}$ and system trajectory $\mathfrak{X}$ from the lastly-solved relaxed convex program. We compute the Euclidean projection as we do in the projection policy. If the "deviation", denoted as $\theta$, of the system to its deterministic counter-part is greater than $\Theta$, then we apply a re-solving as in the update policy. In the meantime we also update $\mathfrak{U}$ and  $\mathfrak{X}$ in the memory. This is summarized in Algorithm \ref{algo:hybrid}. We remark that as compared to the projection policy, the additional computation of $\theta$ in the hybrid policy is almost free. 

Note that in Line 7 of Algorithm \ref{algo:hybrid}, we have defined $\theta = \norme{(\mathfrak{X}(t), \barw, \mathfrak{U}(t)) - (\bx,\bw,\bu_{\pi})}$ as the deviation at each time-step $t$. We believe that this may capture more information than just using $\norme{\mathfrak{U}(t) - \bu_{\pi}}$. The latter quantity plays an essential role in the theory of \emph{error bounds}, see e.g. \cite{luo1996mathematical} and \cite[Chapter 6]{facchinei2003finite}, by considering a (non-negative) residual function such as
\begin{equation*}
  r_t(\bx, \bw, \bu) \defeq \sum_{i=1}^{I(t)} \left( g_{t,i}(\bx, \bw, \bu)\right)_+ + \sum_{j=1}^{J(t)} \abs{ h_{t,j} (\bx, \bw, \bu) }
\end{equation*}
where $(\cdot)_+ = \max \{ \cdot, 0 \}$. This theory aims at establishing bounds of type
\begin{equation}  \label{eq:resdual-function}
  \tau' \cdot r_t(\bx, \bw, \mathfrak{U}(t))^{\gamma'} \le \norme{\mathfrak{U}(t) - \bu_{\pi}} = \norme{\mathfrak{U}(t) - \Pi_{\calU_t(\bx, \bw)}(\mathfrak{U}(t)) } \le \tau \cdot r_t(\bx, \bw, \mathfrak{U}(t))^{\gamma}
\end{equation}
with some positive constants $\tau, \tau'$ and $\gamma, \gamma'$. An immediate application of these error bounds in the implementation of the hybrid policy is the following: the exact value of $\norme{\mathfrak{U}(t) - \bu_{\pi}}$ is not essential, since we only need a comparison of $\theta$ with an empirically chosen threshold $\Theta$. We may very well use the easily computable residual functions in the left and right hand sides of \eqref{eq:resdual-function} to approximate the quantity $\norme{\mathfrak{U}(t) - \bu_{\pi}}$ for this purpose, provided that these bounds can be proven valid. This saves us from the computation of an Euclidean projection, rendering the hybrid policy even more efficient than the update policy.

\section{Performance Results}  \label{sec:performance-results}

This section contains the major theoretical results of this paper. We state first-order and second-order theorems on the optimality gap bounds for both the update and the projection policies respectively in Sections \ref{subsec:update-policy-perf} and \ref{subsec:projection-policy-perf}, and mention briefly the proof ideas. The technical details are given in Appendix \hyperlink{appendix:proof-of-main-results}{EC.2}. We then discuss interpretation of theses bounds in Section \ref{subsec:interpretation-of-the-optimality-gap-bounds}. Sufficient conditions to satisfy the conditions of these theorems are discussed in length in Appendix \hyperlink{appendix:policy-mapping}{EC.3}.

\subsection{The Update Policy}   \label{subsec:update-policy-perf}

Following Remark \ref{rem:necessitiy-of-the-assumtions} and in particular \eqref{eq:pivotal-relation}, our goal is to compute an upper bound on the difference between $V_{\mathrm{update}} (\bx,T)$ and $V_{\mathrm{rel}-} (\bx,T)$. Our analysis relies heavily on sensitivity analysis of the convex program \eqref{eq:original-problem-rel} locally around an optimal solution, which exists by our assumptions. Hence:
\begin{equation}\label{eq:fix-optimal-solution}
  \mbox{we fix once and for all an optimal solution $\bu^*[1,T]$ of \eqref{eq:original-problem-rel} with the corresponding $\bx^*[1,T]$}
\end{equation}
For each $1 \le t \le T$, consider the following \emph{policy mapping}:
\begin{align} \label{eq:policy-mapping}
S^*_t \colon \R^{n_x} \times \R^{n_w} & \rightrightarrows \R^{n_u} \nonumber \\
(\bx,\bw) & \mapsto \left\{ \bu(t) \ \Big| \ \bu[t,T] \mbox{ is an optimal solution to $\hat{V}_{\mathrm{rel}+} (\bx,T+1-t,\bw)$ in \eqref{eq:original-problem-res}}  \right\}
\end{align}
The mapping $S^*_t(\cdot)$ defined above is in general a \emph{multi-function}, or a \emph{set-valued mapping} (where "S" stands for \emph{solution}). We mention that $S_t^*(\cdot)$ are implicit (multi)-functions defined from optimal solutions to a class of parameterized mathematical programs, and in general are not possible to obtain explicit formulas (but see \cite{gast2022lp} and the discussion in Section \hyperlink{subsec:non-degeneracy-and-strict-complementarity}{EC.3.2} for the linear program case). We define the following two well-behaved assumptions on the policy mappings. 

\begin{tcolorbox}[colback=white,colframe=black,boxrule=1pt]

\begin{assumption}[Local Lipschitz-Continuity for $S^*_t(\cdot)$] \label{assumption:lipschitz-continuity} 
For all $1 \le t \le T$, there exists $\varepsilon_t > 0$ such that for all $(\bx,\bw) \in \calB((\bx^*(t),\barw), \varepsilon_t)$, the set $S^*_t(\bx,\bw)$ is single-valued. Moreover, the (single-valued) function $S^*_t(\bx,\bw)$ defined in \eqref{eq:policy-mapping} is locally Lipschitz-continuous in $\calB((\bx^*(t),\barw), \varepsilon_t)$.  
\end{assumption}

\begin{assumption}[Local $\calC^2$-Smoothness for $S^*_t(\cdot)$] \label{assumption:C2-smoothness} 
For all $1 \le t \le T$, there exists $\varepsilon_t > 0$ such that for all $(\bx,\bw) \in \calB((\bx^*(t),\barw), \varepsilon_t)$, the set $S^*_t(\bx,\bw)$ is single-valued. Moreover, the (single-valued) function $S^*_t(\bx,\bw)$ defined in \eqref{eq:policy-mapping} is locally $\calC^2$-smooth in $\calB((\bx^*(t),\barw), \varepsilon_t)$.  
\end{assumption}

\end{tcolorbox}
  
We defer the verification and justification of these assumptions to Appendix \hyperlink{appendix:policy-mapping}{EC.3}. Specifically, Theorems \hyperlink{thm:sufficient-conditions-C2}{EC.1} and \hyperlink{thm:sufficient-conditions-lipschitz}{EC.2} provide sufficient conditions for satisfying these assumptions. Generally, we anticipate that Assumption \ref{assumption:lipschitz-continuity} is applicable across a wide variety of contexts. Meanwhile, the more stringent Assumption \ref{assumption:C2-smoothness} depends on additional regularities necessitated by a certain type of Implicit Function Theorem.

We state and prove two performance bound results under each of these assumptions. The key difference between these two results is that in Theorem \ref{thm:general-theory-of-resolving-Lipschitz}, the constant $C_1$ is expressed by the square root of the variances, while in the constant $C_2$ of Theorem \ref{thm:general-theory-of-resolving-C2}, the square roots have been removed. 

\begin{theorem}[Optimality Gap Bound with Lipschitz-Continuity in Update Policy]  \label{thm:general-theory-of-resolving-Lipschitz}
  Let $V_{\mathrm{opt}} (\bx,T)$ be the value of the stochastic optimization problem \eqref{eq:original-problem} that satisfies Assumptions \ref{assumption-1}-\ref{assumption-8}, and let $V_{\mathrm{update}} (\bx,T)$ be the value of the update policy defined in Algorithm \ref{algo:update}. Under the additional Assumption \ref{assumption:lipschitz-continuity}, there exists constants $C_1, \mathfrak{P}, \bar{V} > 0$ such that 
  \begin{equation*}
    V_{\mathrm{opt}} (\bx,T) - V_{\mathrm{update}} (\bx,T) \le \mathfrak{P} C_1 +  (1 - \mathfrak{P}) \bar{V}
  \end{equation*}
 The constant $\mathfrak{P}$ converges to $1$ exponentially fast as both $\var{\bW}$ and $\var{\calE}$ converge to $0$ (an explicit expression is given in \eqref{eq:event-proba-formula}). The constant $\bar{V}$ is a finite upper bound of $V_{\mathrm{opt}} (\bx,T)$. 
  The constant $C_1 = \calO(\sqrt{\var{\bW}} + \sqrt{\var{\calE}})$, with an explicit expression given in \eqref{eq:choice-of-constant}.
  
\end{theorem}

\begin{theorem}[Optimality Gap Bound with $\calC^2$-Smoothness in Update Policy]  \label{thm:general-theory-of-resolving-C2}
  Let $V_{\mathrm{opt}} (\bx,T)$ be the value of the stochastic optimization problem \eqref{eq:original-problem} that satisfies Assumptions \ref{assumption-1}-\ref{assumption-8}, and let $V_{\mathrm{update}} (\bx,T)$ be the value of the update policy defined in Algorithm \ref{algo:update}. Under the additional Assumption \ref{assumption:C2-smoothness}, there exists constants $C_2, \mathfrak{P}, \bar{V} > 0$ such that 
  \begin{equation*}
    V_{\mathrm{opt}} (\bx,T) - V_{\mathrm{update}} (\bx,T) \le \mathfrak{P} C_2 +  (1 - \mathfrak{P}) \bar{V}
  \end{equation*}
  The constant $\mathfrak{P}$ converges to $1$ exponentially fast as both $\var{\bW}$ and $\var{\calE}$ converge to $0$ (an explicit expression is given in \eqref{eq:event-proba-formula}). The constant $\bar{V}$ is a finite upper bound of $V_{\mathrm{opt}} (\bx,T)$. 
  The constant $C_2 = \calO(\var{\bW} + \var{\calE})$, with an explicit expression given in \eqref{eq:choice-of-constant-refined}.
  
\end{theorem}

Here is the main ingredient of the proofs for these two theorems, details are given in Appendix \hyperlink{appendix:proof-of-main-results}{EC.2}: we first apply the concentration inequality in Lemma \ref{lem:multivariate-hoeffding} to bound the probability of the stochastic trajectory remains inside the $\varepsilon_t$-neighbourhood for all $t$, as required in Assumption \ref{assumption:lipschitz-continuity} (resp. Assumption \ref{assumption:C2-smoothness}). Next, assuming Assumption \ref{assumption:lipschitz-continuity} (resp. Assumption \ref{assumption:C2-smoothness}), we use the Lipschitz-continuity (resp. $C^2$-smoothness) properties to control the deviations by the variance of the stochastic noises. % See Figure \ref{fig:evolution-of-system} for an illustration. 
The $C^2$-smoothness provides a refinement, since the first-order approximation is linear and the stochastic errors cancel out upon taking expectation, hence only second-order terms remain. See Section \ref{subsec:interpretation-of-the-optimality-gap-bounds} for more discussions.

\subsection{The Projection Policy}  \label{subsec:projection-policy-perf}

The analysis of the projection policy relies on understanding how $\Pi_{\calU_t(\bx(t), \bw)}(\bu^*(t))$ behaves for $(\bx(t),\bw)$ in a neighbourhood of $(\bx^*(t),\barw)$, which is a projection onto a perturbed set. Since the Euclidean projection is itself an optimization problem (with a quadratic objective function), we use the same toolkit as in Section \ref{subsec:re-solving-heuristic}. In abbreviation, we write
\begin{equation}\label{eq:proj}
  \proj_t(\bx,\bw) \defeq \Pi_{\calU_t(\bx, \bw)}(\bu^*(t))
\end{equation}
to emphasize the dependence of the projection mapping on $\bx,\bw$ and $t$. Note on the other hand that the vector for applying the projection, $\bu^*(t)$, is \emph{fixed} for all time-step $t$. Also be definition, $\proj_t(\bx^*(t),\barw) = \bu^*(t)$. 

Much like Assumptions \ref{assumption:lipschitz-continuity} and \ref{assumption:C2-smoothness}, we define the Lipschitz-continuity and the $C^2$-smoothness assumptions on the functions $\proj_t(\cdot)$. A discussion on sufficient conditions to satisfy these assumptions are given in Theorem \hyperlink{thm:sufficient-conditions-projector}{EC.3}, which are specified for the Euclidean projector.  

\begin{tcolorbox}[colback=white,colframe=black,boxrule=1pt]

\begin{assumption}[Lipschitz-Continuity for $\proj_t(\cdot)$] \label{assumption:lipschitz-continuity-projection} 
For all $1 \le t \le T$, there exists $\varepsilon_t > 0$ such that $\proj_t(\bx,\bw)$ as defined in \eqref{eq:proj} is locally Lipschitz-continuous in $\calB((\bx^*(t),\barw), \varepsilon_t)$.
\end{assumption}

\begin{assumption}[$\calC^2$-Smoothness for $\proj_t(\cdot)$] \label{assumption:C2-smoothness-projection} 
For all $1 \le t \le T$, there exists $\varepsilon_t > 0$ such that $\proj_t(\bx,\bw)$ as defined in \eqref{eq:proj} is locally $\calC^2$-smooth in $\calB((\bx^*(t),\barw), \varepsilon_t)$.
\end{assumption}

\end{tcolorbox}

The next two theorems are the analogues of Theorem \ref{thm:general-theory-of-resolving-Lipschitz} and Theorem \ref{thm:general-theory-of-resolving-C2} for the projection policy. Proofs are detailed in Appendix \hyperlink{subsec:proof-of-projection}{EC.2.2}.

\begin{theorem}[Optimality Gap Bound with Lipschitz-Continuity in Projection Policy]  \label{thm:general-theory-of-projection-Lipschitz}
Let $V_{\mathrm{opt}} (\bx,T)$ be the value of the stochastic optimization problem \eqref{eq:original-problem} that satisfies Assumptions \ref{assumption-1}-\ref{assumption-8}, and let $V_{\mathrm{proj}} (\bx,T)$ be the value of the projection policy defined in Algorithm \ref{algo:projection}. Under the additional Assumption \ref{assumption:lipschitz-continuity-projection}, there exists constants $C_3, \mathfrak{P}', \bar{V} > 0$ such that 
    \begin{equation*}
    V_{\mathrm{opt}} (\bx,T) - V_{\mathrm{proj}} (\bx,T) \le \mathfrak{P}' C_3 +  (1 - \mathfrak{P}') \bar{V}
  \end{equation*}
The constant $\mathfrak{P}'$ converges to $1$ exponentially fast as both $\var{\bW}$ and $\var{\calE}$ converge to $0$ (an explicit expression is given in \eqref{eq:event-proba-formula-1}). The constant $\bar{V}$ is a finite upper bound of $V_{\mathrm{opt}} (\bx,T)$. The constant $C_3 = \calO(\sqrt{\var{\bW}} + \sqrt{\var{\calE}})$, with an explicit expression given in \eqref{eq:choice-of-constant-projection-Lipschitz}.

\end{theorem}

\begin{theorem}[Optimality Gap Bound with $\calC^2$-Smoothness in Projection Policy]  \label{thm:general-theory-of-projection-C2}
Let $V_{\mathrm{opt}} (\bx,T)$ be the value of the stochastic optimization problem \eqref{eq:original-problem} that satisfies Assumptions \ref{assumption-1}-\ref{assumption-8}, and let $V_{\mathrm{proj}} (\bx,T)$ be the value of the projection policy defined in Algorithm \ref{algo:projection}. Under the additional Assumption \ref{assumption:C2-smoothness-projection}, there exists constants $C_4, \mathfrak{P}', \bar{V} > 0$ such that 
    \begin{equation*}
    V_{\mathrm{opt}} (\bx,T) - V_{\mathrm{proj}} (\bx,T) \le \mathfrak{P}' C_4 +  (1 - \mathfrak{P}') \bar{V}
  \end{equation*}
The constant $\mathfrak{P}'$ converges to $1$ exponentially fast as both $\var{\bW}$ and $\var{\calE}$ converge to $0$ (an explicit expression is given in \eqref{eq:event-proba-formula-1}). The constant $\bar{V}$ is a finite upper bound of $V_{\mathrm{opt}} (\bx,T)$. The constant $C_4 = \calO(\var{\bW} + \var{\calE})$, with an explicit expression given in \eqref{eq:choice-of-constant-projection-C2}.

\end{theorem}

\subsection{Interpretation of the Optimality Gap Bounds} \label{subsec:interpretation-of-the-optimality-gap-bounds}

In this subsection, we illustrate via en elementary observation as how the additional $\calC^2$-smoothness can achieve for a refinement, and justify how and when this can be called a refinement via the previously studied examples in Section \ref{sec:concrete-examples}. 

For this purpose, let $Y$ be a real-valued random variable with mean value $\expect{Y}=y$ and a small variance, with the meaning of "small" that will be made precise in the analysis. Suppose we are interested in finding an upper bound for the quantity $\abs{\expect{v(Y)} - v(y)}$, with a certain real-valued $\calC^2$-smooth function $v(\cdot)$. Let $\bar{v}$ be an upper bound of the function $v(\cdot)$ over the support of $Y$, and we suppose that $\bar{v} < \infty$. Let us choose a small $\varepsilon$-neighbourhood of $y$.  If we only rely on the Lipschitz-continuity of the function $v(\cdot)$, we obtain 
\begin{align}  
  \abs{\expect{v(Y)} - v(y)} & \le L_v \cdot \expect{\abs{Y - y} \cdot \mathds{1}_{Y \in [y-\varepsilon, y+\varepsilon]}} + 2 \bar{v} \cdot \proba{Y \notin [y-\varepsilon, y+\varepsilon]}  \nonumber \\
  & \ \le L_v \cdot \varepsilon + C_1 \cdot \exp{ \left( - \frac{C_2 \cdot \varepsilon^2 }{\var{Y}} \right) }  \label{eq:illustration-1}
\end{align}
with $L_v$ being the Lipschitz-constant of $v(\cdot)$ in the interval $[y-\varepsilon, y+\varepsilon]$, and $C_1,C_2 > 0$ are constants from Lemma \ref{lem:multivariate-hoeffding}. On the other hand, if we use the additional $\calC^2$-smoothness, we have
\begin{equation*}
  v(Y) = v(y) + v'(y) \cdot (Y-y) + \frac{1}{2} v^{''}(\tilde{y}_Y) \cdot (Y-y)^2  
\end{equation*}
with $\tilde{y}_Y$ being some value depending on $Y$. So we deduce that
\begin{align}  
  \abs{\expect{v(Y)} - v(y)} & \le \frac{H_v}{2} \cdot \expect{(Y-y)^2 \cdot \mathds{1}_{Y \in [y-\varepsilon, y+\varepsilon]}} + 2 \bar{v} \cdot \proba{Y \notin [y-\varepsilon, y+\varepsilon]} \nonumber \\
  & \ \le \frac{H_v}{2} \cdot \varepsilon^2 +  C_1 \cdot \exp{\left( - \frac{C_2 \cdot \varepsilon^2 }{\var{Y}}\right)} \label{eq:illustration-2}
\end{align}
where $H_v$ is some upper found for the function $v''(\cdot)$ in $[y-\varepsilon, y+\varepsilon]$. If $\var{Y}$ is much smaller than $\varepsilon^2$, say $\var{Y} = \varepsilon^3$, then the second term in \eqref{eq:illustration-1} and \eqref{eq:illustration-2} is negligible as compared to their respective first terms. If in the meantime, the value of $\varepsilon$ can also chosen to be small enough, the bound in \eqref{eq:illustration-1} is of first-order $\calO(\varepsilon)$, while the bound in \eqref{eq:illustration-2} is of second-order $\calO(\varepsilon^2)$. So in summary, a sufficient condition on $\var{Y}$ and $\varepsilon$ to ensure that the upper bound \eqref{eq:illustration-2} is a refinement to the upper bound \eqref{eq:illustration-1} is $\var{Y} \ll \varepsilon^2 \ll 1$. Note that this is \emph{not} a necessary condition and in any case, a second-order bound, if exists, always provides an alternative estimation to the optimality gap. We refer to Figure \ref{fig:fit-rate-with-sigma} for an illustration of this effect in a numerical study of the network utility maximization example.  

\begin{remark}  [Scaling under Affinity]  \label{rem:scaling-under-affinity}

In Examples \ref{example-1} and \ref{example-2}, both relaxed problems corresponding to \eqref{eq:inventory} and \eqref{eq:WCMDP-1} are linear programs, enabling scaling of the model size, which in turn gives rise to a reduction in variances commensurate with the scaling. For instance, in Example \ref{example-2}, the vectors $\bX$ and $\bU$ denote the number of bandit arms in each of the $n_x$ states and taking each of the $n_u$ actions respectively, totaling $N$. Consequently, the vectors $\bMN \defeq \bX/N$ and $\bYN \defeq \bU/N$ represent the corresponding \emph{proportions}. Given that every (in)equality defining the optimization problem \eqref{eq:WCMDP-1} is affine, dividing each formula by $N$, while allowing \emph{fractional} numbers of arms, transforms it into an equivalent problem involving the proportionality real-valued vectors $\bMN$ and $\bYN$ within some simplex with appropriate dimensions. Specifically, \eqref{eq:WCMDP-2} transforms into $\bMN(t+1) = \phi(\bYN(t)) + \bEN(t)$, where $\phi(\cdot)$ represents an affine function and $\bEN(t)$ serves as a random vector constituting the stochastic component, $\calE$. Here, $\bEN(t)$ adheres to $\var{\calE} = \calO(1/N)$, as proven in \cite[Lemma 1]{gast2022lp}. Thus, as the number of symmetric arms $N$ increases, a second-order result as found in Theorems \ref{thm:general-theory-of-resolving-C2} and \ref{thm:general-theory-of-projection-C2} delivers a faster $\calO(1/N)$ asymptotic convergence rate for the optimality gap bounds, contrasting with the $\calO(1/ \sqrt{N})$ rate seen in Theorem \ref{thm:general-theory-of-resolving-Lipschitz} and Theorem \ref{thm:general-theory-of-projection-Lipschitz}. In fact, \cite{gast2022lp} has demonstrated that with a perfect rounding scheme to fit the original integer-valued problem \eqref{eq:WCMDP-1}, the convergence rate can be enhanced to an even faster rate of $\exp(-\calO(N))$. 

Similarly, in Example \ref{example-1} for the inventory model without the option of replenishing storage, a common scaling method is to multiply both the initial resource and the time horizon by a factor $N$ (not necessarily integer). As shown in \cite{jasin2012re}, a faster $\calO(1/N)$ asymptotic convergence rate can be achieved if the relaxed linear program satisfies a non-degenerate condition. We emphasize that these efficient approximations do not contradict the computational hardness results, e.g. those mentioned in Section \ref{example-2}, since the scaling regimes under the two situations are different. For example, the "Restless Bandits" problem in \cite{Papadimitriou99thecomplexity} under this scaling involves identifying the best $5$, $50$, or $500$ arms out of $10$, $100$, or $1000$ arms, respectively. By contrast, identifying a fixed quantity of $5$ arms, regardless of the total number of arms, is PSPACE-complete. Also, note that in the setting of a convex program studied in this paper, such scaling is in general not possible. For situations where the amplitude of the variances are kept constant, thereby making a consideration of robust optimality more suitable, we refer to Remark \hyperlink{rem:asymptotic-vs-robustness}{EC.1} for further discussion.  \Halmos

\end{remark}

\section{Numerical Results on Network Utility Maximization}  \label{sec:numerical-experiments}

\begin{table}[h]
\centering
\begin{tabular}{|c|c|c|c|c|c|c|}
\hline
\textbf{Path $p$/Link $l$} & \textbf{x} & $\barw$ & \textbf{q} & \textbf{$\alpha$} & \textbf{D} & \textbf{c} \\
\hline
1 & 1 & 2 & 0.6 & 0.5 & 100 & 6 \\
\hline
2 & 1 & 2 & 0.7 & 0.5 & 100 & 4 \\
\hline
3 & 1 & 2 & 0.5 & 0.5 & 100 & 3 \\
\hline
4 & - & - & - & - & - & 4 \\
\hline
5 & - & - & - & - & - & 6 \\
\hline
\end{tabular}

\caption{Explicit parameters used in the numerical example of Figure \ref{fig-network-example}. }

%For $T=3$, this model satisfies the local $\calC^2$-smoothness Assumption \ref{assumption:C2-smoothness} for the update policy, while it does not satisfy the local $\calC^2$-smoothness Assumption \ref{assumption:C2-smoothness-projection} for the projection policy.  

\label{tab:paras-for-network}
\end{table}

\begin{figure}[h!]
    \centering
    \begin{subfigure}{.45\textwidth}
        \centering
        \includegraphics[width=.98\linewidth]{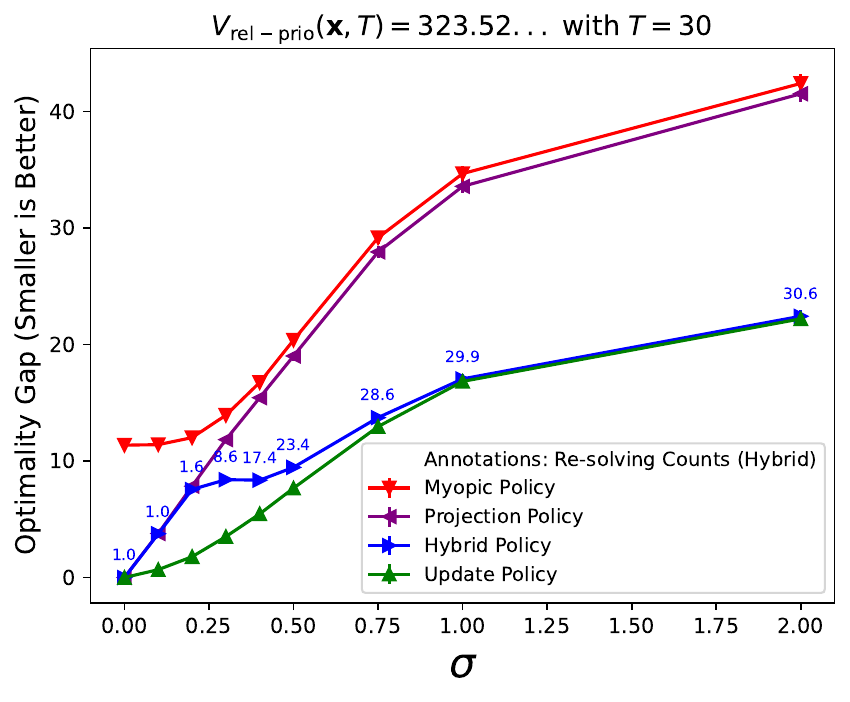}
        \caption{$\Theta \equiv 1.5$ }
        \label{fig:perf-compare-1}
    \end{subfigure}%
    \hfill
    \begin{subfigure}{.45\textwidth}
        \centering
        \includegraphics[width=.98\linewidth]{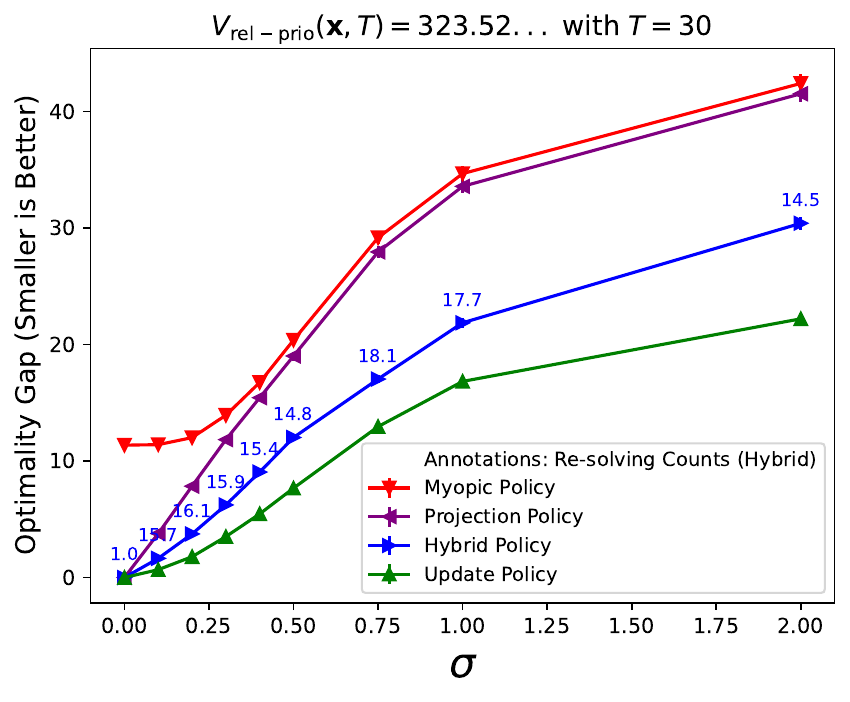}
        \caption{$\Theta = [0.4, 0.8, 1.2, 1.6, 2.0, 2.5, 3.0, 4.0]$}
        \label{fig:perf-compare-3}
    \end{subfigure}
    \begin{subfigure}{\textwidth}
        \centering
        \includegraphics[width=.65\linewidth]{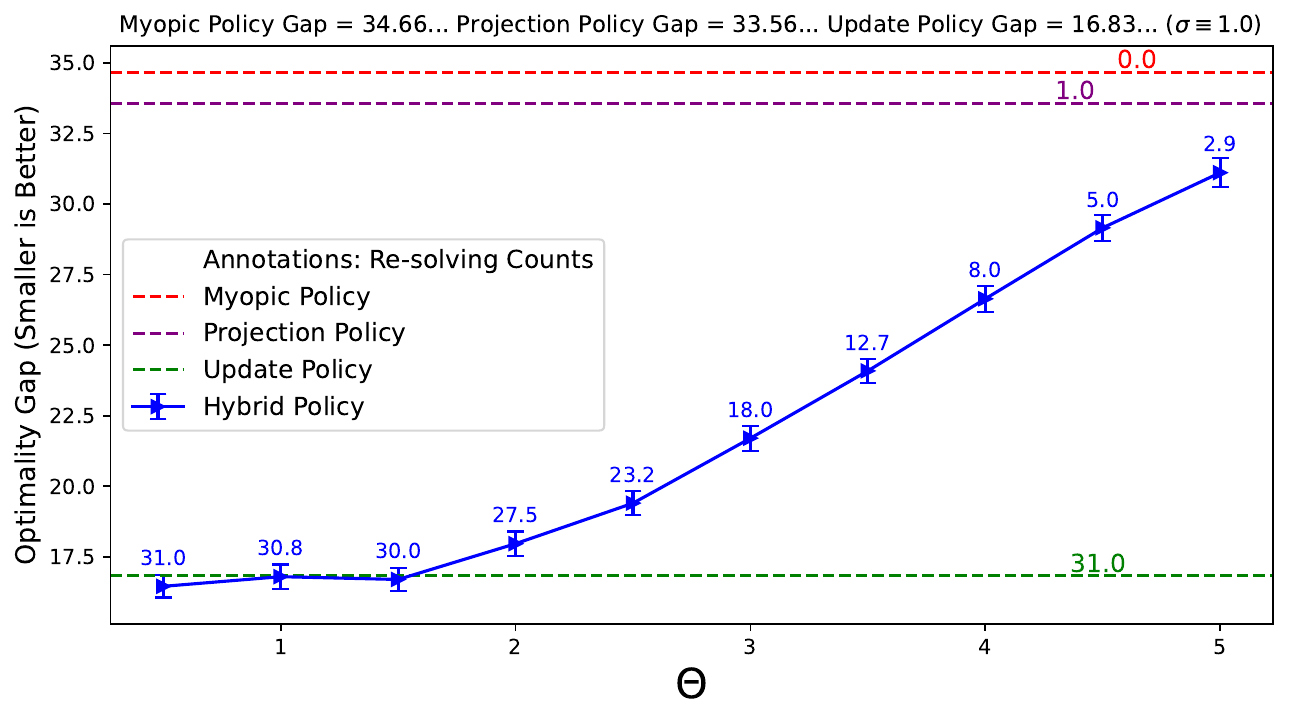}
        \caption{Fix $\sigma \equiv 1.0$ and vary $\Theta$: the trade-off between policy performance and computation efficiency} 
        \label{fig:perf-compare-4}
    \end{subfigure}
    \caption{The model in Figure \ref{fig-network-example} with parameters given in Table \ref{tab:paras-for-network}, horizon is $T=30$. The annotated numbers are the counts of re-solving used in the hybrid policy, averaged over $400$ simulations (the initial solving in Line $1$ of Algorithm \ref{algo:hybrid} before observing $\bW(1)$ is also \emph{included}). Since $T=30$, this count is upper bounded by $31$.}
    \label{fig:perf-compare-all}
\end{figure}

In this section we provide a numerical study on Example \ref{example-3}. The parameters employed in our experiments are detailed in Table \ref{tab:paras-for-network}. Our code, written in Python, can be accessed from our GitHub repository (available via \url{https://gitlab.inria.fr/phdchenyan/network_scheduling.git}). We utilize the CVXPY package for solving convex programs. In addition to the update, projection, and hybrid policies, we also introduce the myopic policy as a benchmark; this policy undertakes a random feasible action at each time-step. This random selection is subject to the specific convex program solver, but is generally considered the least computationally demanding of all feasible policies.

In the first set of experiments, we set the horizon $T=30$ with the aim of comparing the performance of various policies studied in this paper. The simulation results are presented in Figure \ref{fig:perf-compare-all}. We evaluate the variance variable $\sigma$ within the range of $[0, 2]$. In Figure \ref{fig:perf-compare-1}, the tuning parameter $\Theta$ for the hybrid policy remains fixed at $\Theta=1.5$, while in Figure \ref{fig:perf-compare-3}, a different and increasing $\Theta$ is selected for each ascending value of $\sigma$. In Figure \ref{fig:perf-compare-4}, we hold $\sigma = 1$ constant and adjust $\Theta$ within the scope of [0.5,5]. The annotated numbers indicate the counts of re-solving utilized in each policy, inclusive of the initial solving of the convex program to maintain a uniform comparison. This count is capped at $31$. By definition, it stands at $0$ for the myopic policy, at $1$ for the projection policy, and at $31$ for the update policy. For the hybrid policy, it assumes a value between $1$ and $31$ depending on $\Theta$, with a larger $\Theta$ resulting in fewer re-solving counts.

For the hybrid policy, the re-solving count does not account for \emph{when} the re-solving occurs, as this is entirely stochastic and reliant on the realization of each sample run. However, it can serve as an indicator of computational resource consumption. From Figure \ref{fig:perf-compare-all}, we deduce that more frequent re-solving can improve performance, and depending on the variance $\sigma$, a model-driven selection of $\Theta$ can help strike a balance between performance assurance and computational resource consumption. For instance, when $\sigma = 1$, choosing $\Theta = 3$ appears to be a judicious decision for this purpose. Meanwhile, when $\sigma = 0.1$, the projection policy proves ideal: Not only is it substantially more efficient than the update policy, but it also achieves a noteworthy improvement over the myopic policy, merely by solving the initial convex program once. Note that in the asymptotic limit $\sigma = 0$, all policies, except the myopic, achieve performance that coincides with the value of the relaxed convex program.

\begin{figure}[h!]
    \centering
    \begin{subfigure}{0.495\textwidth}
        \centering
        \includegraphics[width=\textwidth]{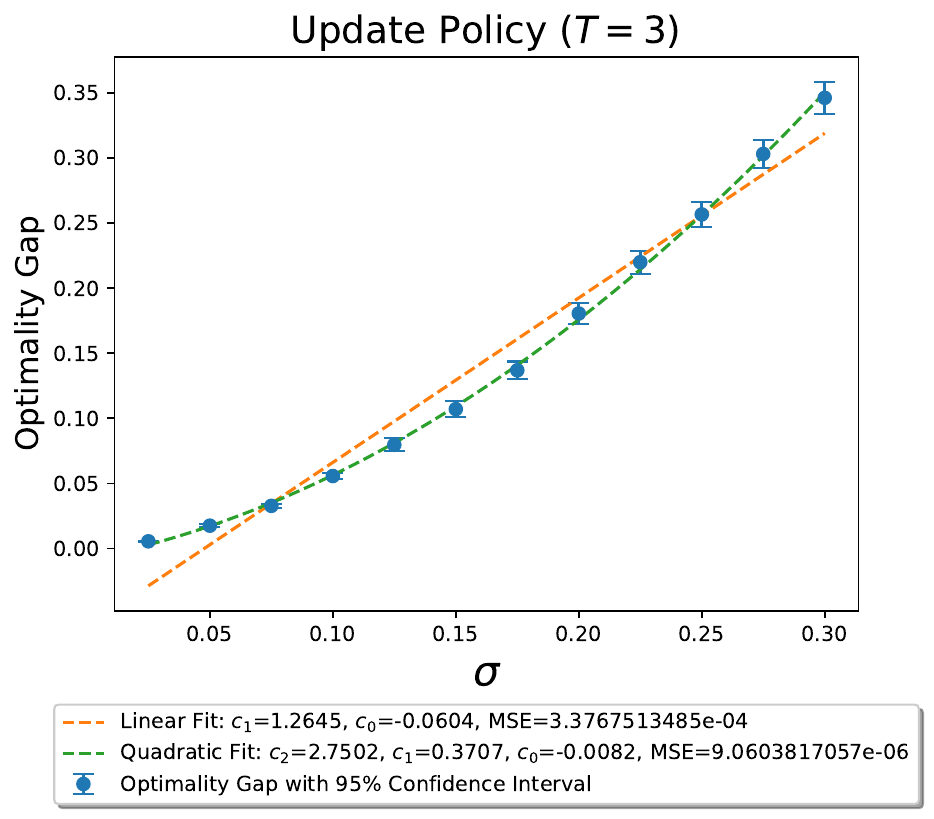}
        %\caption{$\sigma$ in the range of $[0.025,0.3]$ with a step-size of $0.025$}
         \label{fig:update-fit-1}
    \end{subfigure}
    \hfill
    \begin{subfigure}{0.495\textwidth}
        \centering
        \includegraphics[width=\textwidth]{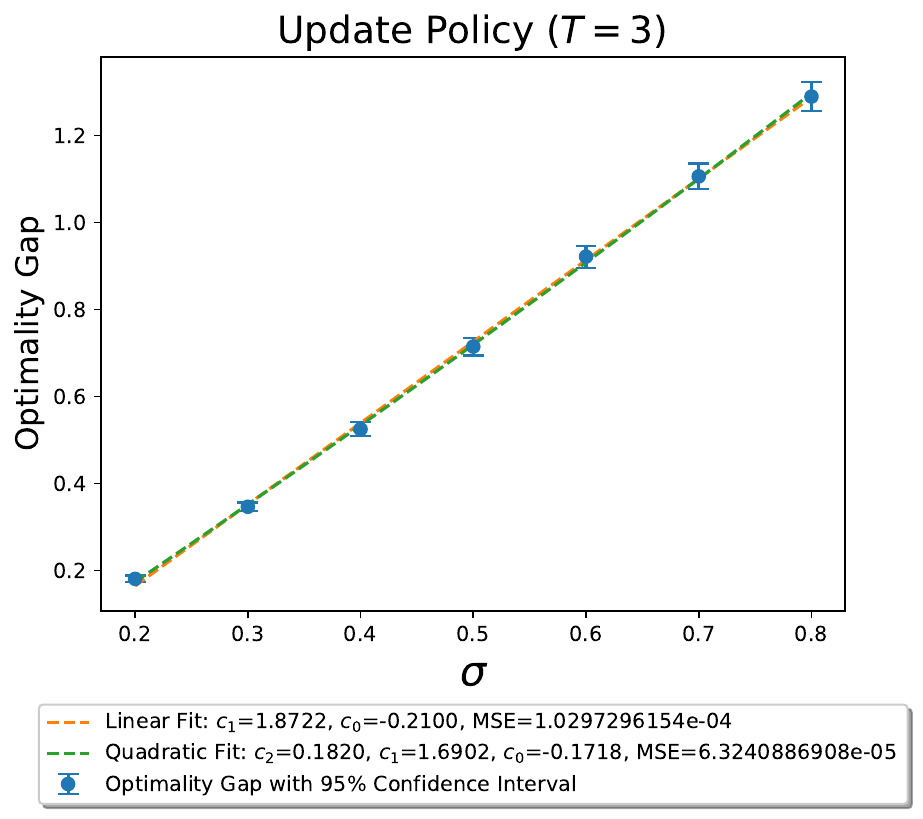}
        %\caption{$\sigma$ in the range of $[0.2,0.8]$ with a step-size of $0.1$}
           \label{fig:update-fit-2}
    \end{subfigure}
    \caption{The model in Figure \ref{fig-network-example} with parameters given in Table \ref{tab:paras-for-network}, horizon is $T=3$. We plot the optimality gap as function of a \emph{single} parameter $\sigma$, which quantifies the strength of the variances, from both $\bW$ and $\calE$, for the update policy. We then apply a linear and quadratic fit. The plots are done in two scale, one for $\sigma$ in the range $[0.025, 0.3]$, where a quadratic fit is more suitable; another in the range $[0.2,0.8]$, where a linear fit appears just as adequate as a quadratic fit. }
    \label{fig:fit-rate-with-sigma}
\end{figure}

For a second set of experiments, we set the horizon $T=3$ and aim to validate the convergence rate backed by our theoretical results. The results are depicted in Figure \ref{fig:fit-rate-with-sigma}. We begin by checking the sufficient conditions in Theorem \ref{thm:general-theory-of-resolving-C2} for the update policy. As per Theorem \hyperlink{thm:sufficient-conditions-C2}{EC.1}, this boils down to validating LICQ and the strict complementarity to the initial convex program \eqref{eq:original-problem-rel}. It should be noted that for the current model, the constraints in \eqref{eq:network-model-1} are actually superfluous and implied by the constraints in \eqref{eq:network-model-5}. Numerically, we discover that for each of the three time-steps, only one constraint in \eqref{eq:network-model-5} is saturated, with the LICQ being upheld. The strict complementarity is also satisfied, after examining the Lagrange multipliers from the solution to \eqref{eq:original-problem-rel}. On the other hand, the projection mapping \eqref{eq:proj} for time-step $t=1$ is degenerate (see Appendix \hyperlink{appendix:policy-mapping}{EC.3} and in particular Figure \ref{fig:degeneracy2} for an illustration), suggesting that the rate claimed in Theorem \ref{thm:general-theory-of-projection-C2} does not apply to the projection policy.

As a consequence of the preceding analysis, we observe in the left panel of Figure \ref{fig:fit-rate-with-sigma} that for the update policy with $\sigma$ in the range of $[0.025,0.3]$, a quadratic fit of the optimality gap is more suitable than a linear fit. However, in the right panel of Figure \ref{fig:fit-rate-with-sigma}, when we extend the range of $\sigma$ to $[0.2,0.8]$, a linear fit appears just as adequate as a quadratic fit for $\sigma$ in this scope. This aligns with the analysis in Section \ref{subsec:interpretation-of-the-optimality-gap-bounds} and the first step of the proof of Theorem \ref{thm:general-theory-of-resolving-Lipschitz} via concentration-type inequalities: when the variance is substantial, the concentration bound of Lemma \ref{lem:multivariate-hoeffding} becomes too broad, and the bound in \eqref{eq:illustration-1} (a first-order linear bound) proves more accurate than the one in \eqref{eq:illustration-2} (a second-order quadratic bound). We also mention that for the projection policy, the fits appear to be linear for $\sigma$ in both the ranges $[0.025,0.3]$ and $[0.2,0.8]$ (not shown in the figure). This supports the observation that the projection mapping is degenerate at $t=1$, and hence the quadratic convergence rate of Theorem \ref{thm:general-theory-of-projection-C2} is not applicable.

\section{Conclusion and Future Works}  \label{sec:conclusion}

We addressed multi-stage stochastic optimization problems that exhibit a convex structure upon applying a CEC-based relaxation. Our framework not only embraces established models from inventory management and Markovian bandits, but also adeptly manages more intricate, non-linear problems like network utility maximization. The heuristic policies we introduced offer a quantifiable performance gap in relation to the optimal ones, while also balance between performance and efficiency. As we conclude, we chart a path for several avenues of future research.

\begin{enumerate} [label=(\roman*)]

\item(Infinite Horizon) In situations where all the data of the model (i.e., rewards, noises, and constraints) are time-independent, we can formulate the problem under the infinite horizon with a time-averaged reward criterion. Specific models of \eqref{eq:original-problem}, such as the "Restless Bandits", have already been analysed within this framework, as seen with the well-known Whittle index policy (\cite{whittle-restless}) and the LP-priority policy (\cite{Ve2016.6}). A distinguishing feature of the infinite horizon case is that, when appropriately formulated, we only have a single relaxed program that describes the system's stationary behavior - significantly simpler than the finite horizon case, which is size $T$ times larger. In principle, the idea from the projection policy in this paper still applies. However, it now necessitates an additional verification of asymptotic stability to ensure optimality. In a recent paper (\cite{hong2023restless}), which studied the "Restless Bandits" model, the authors circumvented this asymptotic stability condition using a coupling technique. Exploring how to integrate their approach into the more general setting in this paper is a worthwhile endeavor.

\item(The Effect of Re-Solving) The numerical experiments in Section \ref{sec:numerical-experiments} indicate that more frequent re-solving always yields improved performance in this network utility maximization model. Additionally, the analysis in Appendix \hyperlink{subsec:Lipschitz-constant-of-V-rel-post}{EC.4.2} hints at the update policy's potential for increased noise immunity. Nonetheless, a comprehensive understanding of the effects of re-solving is far from trivial, since it requires comparing different approximation schemes for a stochastic problem. Note that this question should not be conflated with the concept of time-consistency in multi-stage stochastic programming \cite[Chapter 5]{pflug2014multistage}. Although it may seem counterintuitive, more re-solving is not always beneficial - a fact underscored by counter-examples for the inventory management model found in \cite{cooper2002asymptotic, secomandi2008analysis, bumpensanti2020re}. We plan to explore this issue further and report our findings in a follow-up paper, utilizing a broader framework than the one presented in this paper (convexity is not needed for studying the effect of re-solving). Our numerical experiments appear to suggest that the negative impact of re-solving is related to degeneracy, as we define in Appendix \hyperlink{appendix:policy-mapping}{EC.3}, in conjunction with the presence of stringent constraints.

\item(Discrete-Valuedness) One limitation in our model is the assumption of continuous-valued data. However, in many real-world applications, the model is intrinsically discrete-valued, such as when the state is represented by counting numbers as in the Markovian bandit model of Example \ref{example-1}, or when the action is binary as in some variants of the inventory management model of Example \ref{example-2}. In many cases, continuous-valuedness is more a hindrance than a help, as unless explicit analytical solutions are available, continuous functions must be discretized for numerical solving. In our context, however, it is the inverse: discrete-valuedness disrupts the desired smoothness and convexity properties essential for our CEC-based analysis. Certain situations allow us to overcome these difficulties by relaxing discrete values, e.g. by permitting fractional numbers of arms as discussed in Remark \ref{rem:scaling-under-affinity} and \cite{gast2021lp, gast2022lp}, or by treating binary action as a probability with a value between $0$ and $1$, as in the probabilistic allocation policy in \cite{bumpensanti2020re}.

\end{enumerate}

% Acknowledgments here
\ACKNOWLEDGMENT{ Chen Yan extends heartfelt gratitude to Nicolas Gast and Bruno Gaujal for their invaluable guidance and numerous enlightening discussions on the intriguing subjects related to this research. He also wishes to express his appreciation to IMT Atlantique (Campus de Brest) for its hospitality. The time spent there facilitated many productive conversations with Alexandre Reiffers-Masson, ultimately making this paper possible.  }

\bibliographystyle{informs2014} % outcomment this and next line in Case 1
\bibliography{reference}        % if more than one, comma separated

%% Here starts the e-companion (EC)
%%%%%%%%%%%%%%%%%%%%%%%%%%%%%%%%%%%%%%%%%%%%%%%%%%%%%%%%%%
\ECSwitch

%\ECDisclaimer
%%%%%%%%%%%%%%%%%%%%%%%%%%%%%%%%%%%%%%%%%%%%%%%%%%%%%%%%%%

%%% Main head for the e-companion
\ECHead{Certainty Equivalence Control-Based Heuristics in Multi-Stage Convex Stochastic Optimization Problems (Online Appendix)}

\paragraph*{Outline} Appendices \hyperlink{appendix:proof-well-defined}{EC.1} and \hyperlink{appendix:proof-of-main-results}{EC.2} contain the technical proofs of our theoretical results. In Appendix \hyperlink{appendix:policy-mapping}{EC.3}, we provide a unifying view on the various regularity conditions seen as sufficient for specific upper bounds. Appendix \hyperlink{appendix:policy-classes-extension}{EC.4} extends our proposed heuristic policies to several broader policy classes, followed by a discussion on the issue of the multiplicative exponentially growing constant.

\section{Existence of Optimal Solution}   \hypertarget{appendix:proof-well-defined}{}

In this appendix, we provide a proof for Proposition \ref{prop:existence-of-optimal-solution}, thereby validating the well-posedness of the optimization problem outlined in \eqref{eq:original-problem}. A key aspect within the proof process involves establishing the continuity of a maximization function with respect to a variable that is also part of the constraint set. To ensure this continuity, we leverage a technical result from Proposition 4.4 in \cite{bonnans2013perturbation}, which affirms this continuity provided certain constraint qualification is met on the constrained sets.

\begin{repeatproposition} [Proposition 1]  
  Under Assumptions \ref{assumption-1}-\ref{assumption-8}, the optimization problem \eqref{eq:original-problem} is well-defined and an optimal solution exists. Moreover, the mapping $\bx \mapsto V_{\mathrm{opt}} (\bx,T)$ is a continuous function of $\bx$. 
\end{repeatproposition}

\proof{Proof of Proposition \ref{prop:existence-of-optimal-solution}}
We use dynamic programming and proceed backward in time. Starting at the last time-step $T$, define for $(\bx,\bw)$ the optimal value function
\begin{equation*}
  \hat{V}_{\mathrm{opt}} (\bx,1,\bw) \defeq \max_{\bu \in \calU_T (\bx,\bw)} R_T (\bx, \bw, \bu)
\end{equation*}
By Weierstrass extreme value theorem, $\hat{V}_{\mathrm{opt}} (\bx,1,\bw)$ is well-defined and attained, since it is the maximum value of a continuous function $ \bu \mapsto R_T(\bx, \bw, \bu)$ over a non-empty compact set $\calU_T (\bx,\bw)$, where the compactness follows from our model assumptions. 

We next show that $\bx \mapsto \hat{V}_{\mathrm{opt}} (\bx,1,\bw)$ is continuous as a function of $\bx$, by using \cite[Proposition 4.4]{bonnans2013perturbation}. Note that this proposition is established under a more general setting of infinite dimensional Banach space, and for our finite dimensional case considered here, the only non-trivial condition that needs verification is Robinson's Constraint Qualification (Condition (iv) of this Proposition 4.4). In finite dimension, Robinson's CQ reduces to MFCQ (Mangasarian-Fromovitz Constraint Qualification, see also Appendix \hyperlink{appendix:policy-mapping}{EC.3} for more details). By \cite[Proposition 3.2.7]{facchinei2003finite}, the MFCQ holds at any point $\bu \in \calU_T(\bx,\bw)$ if and only if the Slater CQ holds for $\calU_T(\bx,\bw)$. By the first item of Assumption \ref{assumption-8}, this is true for any feasible set $\calU_T(\bx,\bw)$ parameterized by $(\bx,\bw)$. We conclude that $\bx \mapsto \hat{V}_{\mathrm{opt}} (\bx,1,\bw)$ is indeed a continuous function of $\bx$. Moreover, since $\bx$ takes values in a compact set, we deduce that $\bx \mapsto \hat{V}_{\mathrm{opt}} (\bx,1,\bw)$ is uniformly continuous. 

By definition, we have
\begin{equation*}
  V_{\mathrm{opt}} (\bx,1) = \int \hat{V}_{\mathrm{opt}} (\bx,1,\bw) \cdot f(\bw) d \bw
\end{equation*}
Since $\bx \mapsto \hat{V}_{\mathrm{opt}} (\bx,1,\bw)$ is uniformly continuous for any realisation of $\bw$, we deduce that $\bx \mapsto V_{\mathrm{opt}} (\bx,1)$ is also a continuous function of $\bx$. 

Next for time-step $T-1$, the dynamic programming equation writes
\begin{align*}
  \hat{V}_{\mathrm{opt}} (\bx,2,\bw) & = \max_{\bu \in \calU_{T-1} (\bx,\bw)} R_{T-1} (\bx, \bw, \bu) + \expect{V_{\mathrm{opt}} (\bX,1) \ \Big | \ \bx, \bw, \bu} \\
  & \ = \max_{\bu \in \calU_{T-1} (\bx,\bw)} R_{T-1} (\bx, \bw, \bu) + \int_{\bx' \in \R^{n_x}} V_{\mathrm{opt}} (\phi(\bx,\bw,\bu) + \bx' ,1) \ d \nu(\bx' \mid \bx, \bw, \bu)
\end{align*}
We argue that $\bu \mapsto R_{T-1} (\bx, \bw, \bu) + \int_{\bx' \in \R^{n_x}} V_{\mathrm{opt}} (\phi(\bx,\bw,\bu) + \bx' ,1) \ d \nu(\bx' \mid \bx, \bw, \bu)$ is a continuous function of $\bu$, by relying on the second item of Assumption \ref{assumption-8} for the continuity of the mapping $\bu \mapsto \nu(\cdot \mid \bx,\bw,\bu)$, as well as the fact that $\bu \mapsto R_{T-1} (\bx, \bw, \bu)$ and $\bu \mapsto V_{\mathrm{opt}} (\phi(\bx,\bw,\bu) + \bx' ,1)$ are all continuous, the later follows from our previous analysis for time-step $T$. Consequently, again by Weierstrass extreme value theorem, $\hat{V}_{\mathrm{opt}} (\bx,2,\bw)$ is well-defined and attained for every $(\bx,\bw)$. We show in a similar way as for time-step $T$ that, $\bx \mapsto R_{T-1} (\bx, \bw, \bu) + \int_{\bx' \in \R^{n_x}} V_{\mathrm{opt}} (\phi(\bx,\bw,\bu) + \bx' ,1) \ d \nu(\bx' \mid \bx, \bw, \bu)$ is a continuous function of $\bx$. Again by the first item of Assumption \ref{assumption-8} and \cite[Proposition 4.4]{bonnans2013perturbation}, we deduce that $\bx \mapsto \hat{V}_{\mathrm{opt}} (\bx,2,\bw)$ is a (uniformly) continuous function of $\bx$ for any realisation of $\bw$. Consequently $\bx \mapsto V_{\mathrm{opt}} (\bx,2) \defeq \int \hat{V}_{\mathrm{opt}} (\bx,2,\bw) \cdot f(\bw) d \bw$ is also continuous. 

More generally, suppose that for time-step $t$ the optimization problem $V_{\mathrm{opt}} (\bx,t)$ of \eqref{eq:original-problem} for time-span $[T-t+1,T]$ is well-defined and $\bx \mapsto V_{\mathrm{opt}} (\bx,t)$ is a continuous function of $\bx$. We remind the reader that since we are using backward induction, the "$t$" in $V_{\mathrm{opt}} (\bx,t)$ is indexing the time steps remaining until the end of horizon. We write the dynamic programming equation
\begin{equation*}
  \hat{V}_{\mathrm{opt}} (\bx,t+1,\bw) = \max_{\bu \in \calU_{T-t} (\bx,\bw)} R_{T-t} (\bx, \bw, \bu) + \int_{\bx' \in \R^{n_x}} V_{\mathrm{opt}} (\phi(\bx,\bw,\bu) + \bx', t) \ d \nu(\bx' \mid \bx, \bw, \bu)
\end{equation*}
We use the induction hypothesis and item 2 of Assumption \ref{assumption-8} to show that $\bu \mapsto R_{T-t} (\bx, \bw, \bu) + \int_{\bx' \in \R^{n_x}} V_{\mathrm{opt}} (\phi(\bx,\bw,\bu) + \bx', t) \ d \nu(\bx' \mid \bx, \bw, \bu)$ is continuous, so by Weierstrass extreme value theorem $\hat{V}_{\mathrm{opt}} (\bx,t+1,\bw)$ is wel-defined and attained. We use again the induction hypothesis and item 2 of Assumption \ref{assumption-8} to show that $\bx \mapsto R_{T-t} (\bx, \bw, \bu) + \int_{\bx' \in \R^{n_x}} V_{\mathrm{opt}} (\phi(\bx,\bw,\bu) + \bx', t) \ d \nu(\bx' \mid \bx, \bw, \bu)$ is continuous, and combine with item 1 of Assumption \ref{assumption-8} to show that $\bx \mapsto \hat{V}_{\mathrm{opt}} (\bx,t+1,\bw)$ is (uniformly) continuous for any realisation of $\bw$. Hence $\bx \mapsto V_{\mathrm{opt}} (\bx,t+1) \defeq \int \hat{V}_{\mathrm{opt}} (\bx,t+1,\bw) \cdot f(\bw) d \bw$ is continuous as well. This completes the induction step and concludes the proof.       \Halmos

\endproof

To establish the continuity result presented in Proposition \ref{prop:existence-of-optimal-solution}, we relied on the MFCQ. As we will elaborate in Appendix \hyperlink{appendix:policy-mapping}{EC.3}, stronger constraint qualifications are required to ensure stricter regularities, such as Lipschitz-continuity or $\calC^2$-smoothness.

\section{Proof of the Performance Bounds}     \hypertarget{appendix:proof-of-main-results}{}

In this appendix we prove the performance bounds in Section \ref{sec:main-performance-bound}. The following concentration inequality in vector form will be used in the sequence to bound the probability of the stochastic trajectory leaving outside an $\varepsilon$-neighbourhood.    %(the orange events in Figure \ref{fig:evolution-of-system})
\begin{lemma}[Vector Bernstein Inequality]  \label{lem:multivariate-hoeffding}
  Let $\bY$ be a random vector such that $\expect{\bY} = \by$, $\norme{\bY} < C$ and $\var{\bY} < \infty$. Then for $0 < \varepsilon < \var{\bY}/C$, we have 
  \begin{equation*}
    \proba{\norme{\bY - \by} \ge \varepsilon} \le \exp \left( - \frac{\varepsilon^2}{8 \var{\bY}} + \frac{1}{4} \right)
  \end{equation*}
\end{lemma}
For a proof of Lemma \ref{lem:multivariate-hoeffding} we refer to \cite[Lemma 18]{kohler2017sub}. 

\subsection{The Update Policy}  \hypertarget{subsec:proof-of-update}{}
For convenience, we repeat the assumptions and theorems below for ease of discussion.

\begin{repeatassumption}[Assumption I] 
For all $1 \le t \le T$, there exists $\varepsilon_t > 0$ such that for all $(\bx,\bw) \in \calB((\bx^*(t),\barw), \varepsilon_t)$, the set $S^*_t(\bx,\bw)$ is single-valued. Moreover, the (single-valued) function $S^*_t(\bx,\bw)$ defined in \eqref{eq:policy-mapping} is locally Lipschitz-continuous in $\calB((\bx^*(t),\barw), \varepsilon_t)$.  
\end{repeatassumption}

Combining the smoothness (hence the Lipschitz-continuity) of the reward functions $R_t(\cdot)$, Assumption \ref{assumption:lipschitz-continuity} implies that 
  \begin{itemize}
    \item The function $\hat{V}_{\mathrm{rel}+} (\cdot, T-t, \cdot): \R^{n_x} \times \R^{n_w} \rightarrow \R$ defined in \eqref{eq:original-problem-res} is Lipschitz-continuous with Lipschitz constant $K_t$ in the domain $\calB((\bx^*(t),\barw), \varepsilon_t)$:
        \begin{equation} \label{eq:lipschitz-approximation-1}
          \abs{\hat{V}_{\mathrm{rel}+} (\bx, T-t, \bw) - \hat{V}_{\mathrm{rel}+} (\bx^*(t), T-t, \barw)} \le K_t \cdot \norme{(\bx, \bw) - (\bx^*(t), \barw)}, \ \ \ \forall (\bx,\bw) \in \calB((\bx^*(t),\barw), \varepsilon_t)
        \end{equation}
  \end{itemize}
Since $V_{\mathrm{rel}-} (\bx,T+1-t) = \hat{V}_{\mathrm{rel}+} (\bx,T+1-t,\barw)$, this also implies that
  \begin{itemize}
    \item The function $V_{\mathrm{rel}-} (\cdot, T-t): \R^{n_x} \rightarrow \R$, defined in \eqref{eq:original-problem-rel} is Lipschitz-continuous with Lipschitz constant $L_t$ in the domain $\calB(\bx^*(t), \varepsilon_t)$:
        \begin{equation}  \label{eq:lipschitz-approximation-2}
         \abs{ V_{\mathrm{rel}-} (\bx,T+1-t) - V_{\mathrm{rel}-} (\bx^*(t),T+1-t) } \le L_t \cdot \norme{\bx - \bx^*(t)}, \ \ \ \forall \bx \in \calB(\bx^*(t), \varepsilon_t)
        \end{equation}
  \end{itemize}

\begin{repeattheorem}[Theorem 1 (restated with explicit constants)]  
  Let $V_{\mathrm{opt}} (\bx,T)$ be the value of the stochastic optimization problem \eqref{eq:original-problem} that satisfies Assumptions \ref{assumption-1}-\ref{assumption-8}, and let $V_{\mathrm{update}} (\bx,T)$ be the value of the update policy defined in Algorithm \ref{algo:update}. Under the additional Assumption \ref{assumption:lipschitz-continuity}, there exists constants $C_1, \mathfrak{P}, \bar{V} > 0$ such that 
  \begin{equation*}
    V_{\mathrm{opt}} (\bx,T) - V_{\mathrm{update}} (\bx,T) \le \mathfrak{P} C_1 +  (1 - \mathfrak{P}) \bar{V}
  \end{equation*}
 The constant $\mathfrak{P}$ is given in \eqref{eq:event-proba-formula} and converges to $1$ exponentially fast as both $\var{\bW}$ and $\var{\calE}$ converge to $0$. The constant $\bar{V}$ is a finite upper bound of $V_{\mathrm{opt}} (\bx,T)$.  
  The constant $C_1$ is
  \begin{equation}  \label{eq:choice-of-constant}
     C_1 \defeq \sum_{t=1}^{T} K_{t} \sqrt{\var{\bW}} + L_{t} \sqrt{\var{\calE}}
   \end{equation}
  where 
  \begin{itemize}
    \item $K_t$ is the Lipschitz constant of the function $\hat{V}_{\mathrm{rel}+} (\cdot, T-t, \cdot): \R^{n_x} \times \R^{n_w} \rightarrow \R$ defined in \eqref{eq:original-problem-res} for $(\bx,\bw) \in \calB((\bx^*(t),\barw), \varepsilon_t)$ 
    \item $L_t$ is the Lipschitz constant of the function $V_{\mathrm{rel}-} (\cdot, T-t): \R^{n_x} \rightarrow \R$ defined in \eqref{eq:original-problem-rel} for $\bx \in \calB(\bx^*(t), \varepsilon_t)$
    \item The distribution function of the random vector $\bW$ are given in Assumption \ref{assumption-3}
    \item $\var{\calE} < \infty$ is defined in \eqref{eq:law-of-stochastic-noise} of Assumption \ref{assumption-6}
  \end{itemize}

\end{repeattheorem}

For each time-step $1 \le t \le T$, we introduce the following value function
  \begin{equation} \label{eq:rel-bar-t}
  V_{\mathrm{rel}+} (\bx(t),T+1-t) \defeq \int \hat{V}_{\mathrm{rel}+} (\bx(t), T+1-t, \bw) \cdot f(\bw) d \bw
\end{equation} 
where recall that $f(\bw)$ is the distribution function of the random vector $\bW$. An interpretation of $V_{\mathrm{rel}+} (\bx(t),T+1-t)$ is that it is the expected value of the deterministic optimal control at time-step $t$, with the system being in $\bx(t)$, plus the additional \emph{hindsight} of the exact value of $\bw = \bW(t)$. In contrast, $V_{\mathrm{rel}-} (\bx(t),T+1-t)$ from \eqref{eq:original-problem-rel} can be seen as the expected value of the deterministic optimal control at time-step $t$, with the system being in $\bx(t)$, and zero knowledge of the realization of $\bW(t)$. So we have the bound $V_{\mathrm{rel}-} (\bx(t),T+1-t) \le V_{\mathrm{rel}+} (\bx(t),T+1-t)$. Their difference can be interpreted as the \emph{value of perfect information of knowing $\bW$}, see \cite{avriel1970value,huang1977sharp}. We are now ready to prove the theorem.

\proof{Proof of Theorem \ref{thm:general-theory-of-resolving-Lipschitz}}
We divide the proof into several steps.

\paragraph{Step One: Bounding the Probability of Leaving the $\varepsilon_t$-Neighbourhood} \

For $1 \le t \le T$, denote by $J_t$ the Lipschitz constant of the policy mapping $S^*_t$ inside $\calB((\bx^*(t),\barw), \varepsilon_t)$ claimed in Assumption \ref{assumption:lipschitz-continuity}, and write $J \defeq \max_t J_t$. We have
\begin{equation*}
  \bx^*(t+1) = \phi(\bx^*(t), \barw, S^*_t(\bx^*(t), \barw)) 
\end{equation*}
and 
\begin{equation*}
  \bX(t+1) = \phi(\bX(t), \bW(t), S^*_t(\bX(t), \bW(t))) + \calE(t)
\end{equation*}
where we have abbreviated $\calE(t) \defeq \calE(\bX(t), \bW(t), S^*_t(\bX(t), \bW(t)))$. Since $\phi(\cdot)$ is an affine function, we denote by $c_{\phi}$ its Lipschitz constant. Making the difference of the above two equations and re-arranging terms, we obtain
\begin{equation*}
  \norme{\bX(t+1) - \bx^*(t+1)} \le c_{\phi} (J+1) \norme{\bX(t) - \bx^*(t)} + c_{\phi} (J+1) \norme{\bW(t) - \barw} + \norme{\calE(t)}
\end{equation*}
Denote by $a \defeq c_{\phi} (J+1)$, $\varepsilon \defeq \min_t \varepsilon_t/2$, and 
\begin{equation*}
  b \defeq \begin{cases}
        \varepsilon/T, & \mbox{if } a = 1 \\
        \varepsilon(a-1)/a^T, & \mbox{if } a > 1 \\
        \varepsilon (1 - a), & \mbox{if } 0 < a < 1
      \end{cases}
\end{equation*}
An elementary calculation shows that 
\begin{equation*}
  \mbox{$a \norme{\bW(t) - \barw} + \norme{\calE(t)} \le b$ holds for all $t$} \Rightarrow \norme{\bX(t) - \bx^*(t)} \le \varepsilon_t/2 \mbox{ holds for all } t
\end{equation*}
Hence we deduce that
\begin{align} 
   & \norme{\calE(t)} \le b/2 \mbox{ and } \norme{\bW(t) - \barw} \le \min \{ b/(2a), \varepsilon \} \mbox{ hold for all time-step $t$}  \nonumber \\
  \Rightarrow & \qquad \mbox{the stochastic trajectory $\bX(t)$ together with all the realizations of $\bW(t)$ during the update} \nonumber \\
   & \qquad \mbox{policy remains inside the $\varepsilon_t$-neighbourhood required by Assumption \ref{assumption:lipschitz-continuity} for all time-step $t$}  \label{eq:event-for-time-t}
\end{align}
By Lemma \ref{lem:multivariate-hoeffding}, the event described on the left hand side of \eqref{eq:event-for-time-t} occurs with probability at least
\begin{equation} \label{eq:event-proba-formula}
 \mathfrak{P} \defeq \left[1 - \exp \left(- \frac{(b/2)^2}{8 \var{\calE}} + \frac{1}{4} \right) \right ]^T  \left[1 - \exp \left(- \frac{(\min \{ b/(2a), \varepsilon \})^2}{8 \var{\bW}} + \frac{1}{4} \right) \right ]^T
\end{equation}

\paragraph{Step Two: Computation Inside the $\varepsilon_t$-Neighbourhood} \

In this step we suppose that the event described in \eqref{eq:event-for-time-t} occurs. Fix $1 \le t \le T$ and suppose that arriving at time-step $t$, the system configuration is in $\bX(t) = \bx(t)$. We can bound the difference between $V_{\mathrm{rel}-} (\bx(t),T+1-t)$ and $V_{\mathrm{rel}+} (\bx(t),T+1-t)$ by
\begin{align} \label{eq:difference-of-prio-and-post}
   & \abs{V_{\mathrm{rel}-} (\bx(t),T+1-t) - V_{\mathrm{rel}+} (\bx(t),T+1-t)} \\
  \le &   \int \abs{ \hat{V}_{\mathrm{rel}+} (\bx(t), T+1-t, \expect{\bW}) - \hat{V}_{\mathrm{rel}+} (\bx(t), T+1-t, \bw) } \cdot f(\bw) d \bw \nonumber \\
  \le & \ K_{t-1} \int \norme{\expect{\bW} - \bw} \cdot f(\bw) d \bw \le K_{t-1} \sqrt{\var{\bW}}  \nonumber
\end{align}
where in the last line we have applied \eqref{eq:lipschitz-approximation-1} from Assumption \ref{assumption:lipschitz-continuity} with the implication from event \eqref{eq:event-for-time-t}. Next observe that according to the definition of the update policy, it is planned by using the full knowledge of $\bw = \bW(t)$, so we can write 
  \begin{equation} \label{eq:update-bar-t}
    V_{\mathrm{update}} (\bx(t),T+1-t) = \int \hat{V}_{\mathrm{update}} (\bx(t), T+1-t, \bw) \cdot f(\bw) d \bw
  \end{equation}
where
\begin{align}
  & \hat{V}_{\mathrm{update}} (\bx(t), T+1-t, \bw) \defeq R_{t} \left( \bx(t), \bw, \bu^*_{\bx,T+1-t,\bw}(t) \right)  \nonumber \\
  & + \expect{V_{\mathrm{update}} (\bX(t+1),T-t) \ \Big| \ \bX(t) = \bx(t), \bW(t) = \bw, \bU(t) = \bu^*_{\bx(t),T+1-t,\bw}(t)} \label{eq:principle-of-optimality-1}
\end{align}
In the above formula, $R_{t} \left( \bx(t), \bw, \bu^*_{\bx,T+1-t,\bw}(t) \right)$ is the instantaneous reward gained at time-step $t$, with the system being in configuration $\bx(t)$ and the control $\bu^*_{\bx,T+1-t,\bw}(t)$ obtained by solving \eqref{eq:original-problem-res} is applied. The system then evolves to the configuration $\bX(t+1)$ at time-step $t+1$, and the second term accounts for the expected reward onwards. 

On the other hand, by Bellman's principle of optimality, we claim that
\begin{equation}  \label{eq:principle-of-optimality}
  \hat{V}_{\mathrm{rel}+} (\bx(t),T+1-t,\bw) = R_{t} \big( \bx(t), \bw, \bu^*_{\bx(t),T+1-t,\bw}(t) \big) +  V_{\mathrm{rel}-} (\phi \big(\bx(t), \bw, \bu^*_{\bx(t),T+1-t,\bw}(t) \big),T-t)
\end{equation}
Indeed, $\bu^*_{\bx(t),T+1-t,\bw}(t)$ is an optimal control in the first-step from solving \eqref{eq:original-problem-res} for $\hat{V}_{\mathrm{rel}+} (\bx(t),T+1-t,\bw)$, and $\hat{\bx}^*(t+1) \defeq \phi \big(\bx(t), \bw, \bu^*_{\bx(t),T+1-t,\bw}(t) \big)$ is an optimal trajectory in the second-step. Since from the second-step onward, the optimal control and the optimal trajectory for $\hat{V}_{\mathrm{rel}+} (\bx(t),T+1-t,\bw)$ coincide with the ones for $V_{\mathrm{rel}-} (\hat{\bx}^*(t+1),T-t)$ in \eqref{eq:original-problem-rel}, we deduce \eqref{eq:principle-of-optimality}.

We remark that the reward at time-step $t$ in $\hat{V}_{\mathrm{update}} (\bx(t), T+1-t, \bw)$ and $\hat{V}_{\mathrm{rel}+} (\bx(t),T+1-t,\bw)$ are the same. Hence by taking the difference of \eqref{eq:update-bar-t} and \eqref{eq:rel-bar-t}, and using \eqref{eq:principle-of-optimality-1} and \eqref{eq:principle-of-optimality}, we obtain 
\begin{align}  \label{eq:difference}
 & V_{\mathrm{update}} (\bx(t),T+1-t) - V_{\mathrm{rel}+} (\bx(t),T+1-t) \nonumber  \\
 =  & \int \expect{V_{\mathrm{update}} (\bX(t+1),T-t) \ \Big| \ \bX(t) = \bx(t), \bW(t) = \bw, \bU(t) = \bu^*_{\bx(t),T+1-t,\bw}(t)} \cdot f(\bw) d \bw \nonumber \\ 
 & \qquad - \int V_{\mathrm{rel}-} (\phi \big(\bx(t), \bw, \bu^*_{\bx(t),T+1-t,\bw}(t) \big),T-t) \cdot f(\bw) d \bw 
\end{align}
We write out the following decomposition, using the abbreviation $\bx^*(t+1) \defeq \phi \big(\bx(t), \bw, \bu^*_{\bx(t),T+1-t,\bw}(t) \big)$ and $\bu^*(t) \defeq \bu^*_{\bx(t),T+1-t,\bw}(t)$:
\begin{align} \label{eq:decomposition-of-difference}
 & \expect{V_{\mathrm{update}} (\bX(t+1),T-t) \ \Big| \ \bX(t) = \bx(t), \bW(t) = \bw, \bU(t) = \bu^*(t)} - V_{\mathrm{rel}-} \big(\bx^*(t+1),T-t \big) \nonumber   \\ 
 = &  \ \underbrace{\expect{ V_{\mathrm{update}} (\bX(t+1),T-t) - V_{\mathrm{rel}+} (\bX(t+1),T-t) \ \Big| \ \bX(t) = \bx(t), \bW(t) = \bw, \bU(t) = \bu^*(t)}}_{\mytag{Term A}{term:A}} \nonumber \\
 & + \underbrace{\expect{ V_{\mathrm{rel}+} (\bX(t+1),T-t) - V_{\mathrm{rel}-} (\bX(t+1),T-t) \ \Big| \ \bX(t) = \bx(t), \bW(t) = \bw, \bU(t) = \bu^*(t)}}_{\mytag{Term B}{term:B}}  \\   \nonumber
 & + \underbrace{\expect{ V_{\mathrm{rel}-} (\bX(t+1),T-t) - V_{\mathrm{rel}-} \big(\bx^*(t+1),T-t \big) \ \Big| \ \bX(t) = \bx(t), \bW(t) = \bw, \bU(t) = \bu^*(t) }}_{\mytag{Term C}{term:C}}   \nonumber
\end{align}
We now analyse the three terms in \eqref{eq:decomposition-of-difference}. Denote by 
\begin{equation*}
  Z(t,\bX(t)) \defeq \expect{V_{\mathrm{update}} (\bX(t),T+1-t) - V_{\mathrm{rel}+} (\bX(t),T+1-t) \ \Big| \ \bX(t)}
\end{equation*} 
which is interpreted as follows: conditional on the value of $\bX(t)$, the difference between the expected value of the update policy performance on system state $\bX(t)$ for the time-span $[t , T]$ with the upper bound $V_{\mathrm{rel}+} (\bX(t),T+1-t)$. Applying $\expect{\ \cdot \ \big| \ \bX(t)}$ in \eqref{eq:difference} and use \eqref{eq:decomposition-of-difference}, we obtain
\begin{equation} \label{eq:decomposition-of-terms}
  Z(t,\bX(t)) = \expect{\int (\mbox{\ref{term:A}} + \mbox{\ref{term:B}} + \mbox{\ref{term:C}} ) \cdot f(\bw) d \bw \ \Big| \ \bX(t)}
\end{equation}
For $\mbox{\ref{term:A}}$, we have
\begin{align}  \label{eq:analysis-of-term-A}
  & \expect{\int (\mbox{\ref{term:A}} ) \cdot f(\bw) d \bw \ \Big| \ \bX(t)}  \nonumber \\
  = & \ \expect{\int \expect{ V_{\mathrm{update}} (\bX(t+1),T-t) - V_{\mathrm{rel}+} (\bX(t+1),T-t) \ \Big| \ \bW(t) = \bw, \bU(t) = \bu^*(t)}  \cdot f(\bw) d \bw \ \Big | \ \bX(t) }  \nonumber \\
  = & \ \expectup{ Z(t+1,\bX(t+1)) \ \Big| \ \bX(t)} 
\end{align}
where $\expectup{ Z(t+1,\bX(t+1)) \ \big| \ \bX(t)} $ means conditional on $\bX(t)$, the expected value of the quantity $Z(t+1,\bX(t+1))$ under the update policy.

For $\mbox{\ref{term:B}}$, by \eqref{eq:difference-of-prio-and-post} we have
\begin{equation*}
  \norme{ \expect{\int (\mbox{\ref{term:B}} ) \cdot f(\bw) d \bw} \ \Big| \ \bX(t) } \le K_{t} \sqrt{\var{\bW}}
\end{equation*}

For $\mbox{\ref{term:C}}$, we have
\begin{align*}
   & \norme{ \expect{\int (\mbox{\ref{term:C}} ) \cdot f(\bw) d \bw} \ \Big| \ \bX(t) } \\
  \le & \ L_{t} \int \expect{\norme{\bX(t+1) - \bx^*(t+1)} \ \big| \ \bX(t) = \bx(t), \bW(t) = \bw, \bU(t) = \bu^*(t) } \cdot f(\bw) d \bw  \\
  \le & \ L_{t} \sqrt{\var{\calE}}
\end{align*}
where in the first inequality we have applied \eqref{eq:lipschitz-approximation-2} from Assumption \ref{assumption:lipschitz-continuity} with the implication from event \eqref{eq:event-for-time-t}, and the second inequality follows from Assumption \ref{assumption-6}.

Now if we write 
\begin{equation} \label{eq:definition-of-small-z}
  z(t) \defeq \expectup{Z(t,\bX(t)) \mid \bX(1) = \bx}
\end{equation}
which is the expected value under the update policy for the quantity $Z(t,\bX(t))$ conditional merely on the initial system state $\bX(1)$, we have
\begin{align*}
  z(t) - z(t+1) &  = \expectup{Z(t,\bX(t)) - Z(t+1, \bX(t+1)) \ \big| \ \bX(1) = \bx } \\
    & = \expectup{Z(t,\bX(t)) - \expectup{Z(t+1, \bX(t+1)) \mid  \bX(t) }  \ \big| \ \bX(1) = \bx }
\end{align*}
So from the above calculations we obtain
\begin{equation*}
  \abs{z(t) - z(t+1)}  \le K_{t} \sqrt{\var{\bW}} + L_{t} \sqrt{\var{\calE}}
\end{equation*}

\paragraph{Step Three: Conclusion of the Proof} \

Since
\begin{equation*}
  z(1) = V_{\mathrm{update}} (\bx,T) - V_{\mathrm{rel}+} (\bx,T)
\end{equation*}
and $z(T+1) = 0$, we conclude that
\begin{align*}
   & \abs{V_{\mathrm{update}} (\bx,T) - V_{\mathrm{opt}} (\bx,T)} \le \abs{V_{\mathrm{update}} (\bx,T) - V_{\mathrm{rel}+} (\bx,T)}  \\
  = & \  \abs{z(1)} \le \sum_{t=1}^{T} \abs{z(t) - z(t+1)}  \\
  \le & \ \mathfrak{P} \left( \sum_{t=1}^{T} K_{t} \sqrt{\var{\bW}} + L_{t} \sqrt{\var{\calE}} \right) + (1 - \mathfrak{P}) \bar{V}
\end{align*}        \Halmos

\endproof

We next recall Assumption \ref{assumption:C2-smoothness}:

\begin{repeatassumption}[Assumption II] 
For all $1 \le t \le T$, there exists $\varepsilon_t > 0$ such that for all $(\bx,\bw) \in \calB((\bx^*(t),\barw), \varepsilon_t)$, the set $S^*_t(\bx,\bw)$ is single-valued. Moreover, the (single-valued) function $S^*_t(\bx,\bw)$ defined in \eqref{eq:policy-mapping} is locally $\calC^2$-smooth in $\calB((\bx^*(t),\barw), \varepsilon_t)$.  
\end{repeatassumption}

Combining the smoothness of the reward functions $R_t(\cdot)$, Assumption \ref{assumption:C2-smoothness} implies that 
\begin{itemize}
    \item The function $\hat{V}_{\mathrm{rel}+} (\cdot, T-t, \cdot): \R^{n_x} \times \R^{n_w} \rightarrow \R$ defined in \eqref{eq:original-problem-res} is $\calC^2$-smooth in the domain $\calB((\bx^*(t),\barw), \varepsilon_t)$, so that
        \begin{align} \label{eq:second-order-approximation-1}
          \hat{V}_{\mathrm{rel}+} (\bx, T-t, \bw) & = \hat{V}_{\mathrm{rel}+} (\bx^*(t), T-t, \barw) + \left[(\bx,\bw) - (\bx^*(t),\barw)\right] \cdot \nabla_{(\bx,\bw)}^\top \hat{V}_{\mathrm{rel}+} (\bx^*(t), T-t, \barw)  \nonumber \\
          & \qquad + \frac{1}{2} \left[(\bx,\bw) - (\bx^*(t),\barw)\right] \cdot \nabla^2_{(\bx,\bw)} \hat{V}_{\mathrm{rel}+} (\bx', T-t, \bw') \cdot \left[(\bx,\bw) - (\bx^*(t),\barw)\right]^\top 
        \end{align}
        for a certain $(\bx',\bw')$ along the line segment from $(\bx, \bw)$ to $(\bx^*(t),\barw)$ in $\R^{n_x} \times \R^{n_w}$. By our choice of norms, we have
        \begin{align*}
         & \abs{\left[(\bx,\bw) - (\bx^*(t),\barw)\right] \cdot \nabla^2_{(\bx,\bw)} \hat{V}_{\mathrm{rel}+} (\bx', T-t, \bw') \cdot \left[(\bx,\bw) - (\bx^*(t),\barw)\right]^\top} \\
          & \qquad \le \norme{(\bx,\bw) - (\bx^*(t),\barw)}^2 \norme{\nabla^2_{(\bx,\bw)} \hat{V}_{\mathrm{rel}+} (\bx', T-t, \bw')}
        \end{align*}
        Denote by $\mathfrak{K}_t (\bx',\bw') \defeq \norme{\nabla^2_{(\bx,\bw)} \hat{V}_{\mathrm{rel}+} (\bx', T-t, \bw')}$. We have
         \begin{equation*}
           \sup_{(\bx', \bw') \in \calB((\bx^*(t),\barw), \varepsilon_t)} \mathfrak{K}_t (\bx',\bw') \defeq \norminf{\mathfrak{K}_t} < \infty
         \end{equation*}
         
    \item The function $V_{\mathrm{rel}-} (\cdot, T-t): \R^{n_x} \rightarrow \R$, defined in \eqref{eq:original-problem-rel} is $\calC^2$-smooth in the domain $\calB(\bx^*(t), \varepsilon_t)$, so that
        \begin{align} \label{eq:second-order-approximation-2}
            V_{\mathrm{rel}-} (\bx, T-t) & = V_{\mathrm{rel}-} (\bx^*(t), T-t) + (\bx - \bx^*(t)) \cdot \nabla_{\bx}^\top V_{\mathrm{rel}-} (\bx^*(t), T-t)  \nonumber \\
           & \qquad + \frac{1}{2} (\bx - \bx^*(t)) \cdot \nabla^2_{\bx} V_{\mathrm{rel}-} (\bx', T-t) \cdot (\bx - \bx^*(t))^\top 
        \end{align}
        for a certain $\bx'$ along the line segment from $\bx$ to $\bx^*(t)$ in $\R^{n_x}$. We have
        \begin{equation*}
          \sup_{\bx' \in \calB(\bx^*(t), \varepsilon_t)} \norme{\nabla^2_{\bx} V_{\mathrm{rel}-} (\bx', T-t)} \defeq \norminf{\mathfrak{L}_t} < \infty
        \end{equation*} 
  \end{itemize}
  
\begin{repeattheorem}[Theorem 2 (restated with explicit constants)]  
  Let $V_{\mathrm{opt}} (\bx,T)$ be the value of the stochastic optimization problem \eqref{eq:original-problem} that satisfies Assumptions \ref{assumption-1}-\ref{assumption-8}, and let $V_{\mathrm{update}} (\bx,T)$ be the value of the update policy defined in Algorithm \ref{algo:update}. Under the additional Assumption \ref{assumption:C2-smoothness}, there exists constants $C_2, \mathfrak{P}, \bar{V} > 0$ such that 
  \begin{equation*}
    V_{\mathrm{opt}} (\bx,T) - V_{\mathrm{update}} (\bx,T) \le \mathfrak{P} C_2 +  (1 - \mathfrak{P}) \bar{V}
  \end{equation*}
  The constant $\mathfrak{P}$ is given in \eqref{eq:event-proba-formula} and converges to $1$ exponentially fast as both $\var{\bW}$ and $\var{\calE}$ converge to $0$. The constant $\bar{V}$ is a finite upper bound of $V_{\mathrm{opt}} (\bx,T)$. 
  The constant $C_2$ is
  \begin{equation}  \label{eq:choice-of-constant-refined}
     C_2 \defeq \frac{1}{2} \sum_{t=1}^{T} \norminf{\mathfrak{K}_t} \var{\bW} + \norminf{\mathfrak{L}_t} \var{\calE}
   \end{equation}
  where
  \begin{itemize}
    \item $\norminf{\mathfrak{K}_t} = \sup_{(\bx', \bw') \in \calB((\bx^*(t),\barw), \varepsilon_t)} \norme{\nabla^2_{(\bx,\bw)} \hat{V}_{\mathrm{rel}+} (\bx', T-t, \bw')}$ 
    \item $\norminf{\mathfrak{L}_t} = \sup_{\bx' \in \calB(\bx^*(t), \varepsilon_t)} \norme{\nabla^2_{\bx} V_{\mathrm{rel}-} (\bx', T-t)}$
    \item The distribution function of the random vector $\bW$ are given in Assumption \ref{assumption-3}
    \item $\var{\calE} < \infty$ is defined in \eqref{eq:law-of-stochastic-noise} of Assumption \ref{assumption-6}
  \end{itemize}
\end{repeattheorem}

\proof{Proof of Theorem \ref{thm:general-theory-of-resolving-C2}}
 The proof uses exactly the same lines of analysis as in Theorem \ref{thm:general-theory-of-resolving-Lipschitz}, up to arriving at Equation \eqref{eq:analysis-of-term-A}. From there, using \eqref{eq:second-order-approximation-1}, instead of \eqref{eq:difference-of-prio-and-post}, we have
 \begin{align} \label{eq:difference-of-prio-and-post-refine} 
   & V_{\mathrm{rel}-} (\bx(t),T+1-t) - V_{\mathrm{rel}+} (\bx(t),T+1-t)  \nonumber \\
  = &   \int  \left( \hat{V}_{\mathrm{rel}+} (\bx(t), T+1-t, \expect{\bW}) - \hat{V}_{\mathrm{rel}+} (\bx(t), T+1-t, \bw) \right) \cdot f(\bw) d \bw \nonumber \\
  = & \ \int (\expect{\bW}) - \bw) \cdot \nabla_{\bw}^\top \hat{V}_{\mathrm{rel}+} (\bx, T-t, \expect{\bW})) \cdot f(\bw) d \bw \nonumber \\
  & \qquad + \frac{1}{2} \int \left( (\expect{\bW}) - \bw) \cdot \nabla^2_{\bw} \hat{V}_{\mathrm{rel}+} (\bx, T-t, \bw') \cdot (\expect{\bW}) - \bw)^\top  \right) \cdot f(\bw) d \bw \nonumber \\
  = & \frac{1}{2} \int \left( (\expect{\bW}) - \bw) \cdot \nabla^2_{\bw} \hat{V}_{\mathrm{rel}+} (\bx, T-t, \bw') \cdot (\expect{\bW}) - \bw)^\top \right) \cdot f(\bw) d \bw
\end{align}
where in the above $\bw' = \bw'(\bw,\bw^*)$ depends on $\bw$ and $\bw^*$. The crucial part in \eqref{eq:difference-of-prio-and-post-refine} is that the first-order term cancels out, and only the second-order term is left.

Similarly, using \eqref{eq:second-order-approximation-2} and the abbreviation $\bx^*(t+1) \defeq \phi \big(\bx(t), \bw, \bu^*_{\bx(t),T+1-t,\bw}(t) \big)$, $\calE(t) \defeq \bX(t+1) - \bx^*(t+1)$, our calculation concerning \ref{term:C} becomes 
\begin{align*} 
   & V_{\mathrm{rel}-} (\bX(t+1),T-t) - V_{\mathrm{rel}-} \big(\bx^*(t+1),T-t \big)   \\
  = & \   \calE(t) \cdot \nabla_{\bx}^\top V_{\mathrm{rel}-} (\bx^*(t+1), T-t) + \frac{1}{2} \calE(t) \cdot \nabla^2_{\bx} V_{\mathrm{rel}-} (\bx', T-t) \cdot \calE(t)^\top 
\end{align*}
where in the above $\bx' = \bx'(\bX(t+1), \bx^*(t+1))$ depends on $\bX(t+1)$ and $\bx^*(t+1)$. Consequently 
\begin{equation}  \label{eq:analysis-of-term-C}
  \mbox{\ref{term:C}} = \frac{1}{2} \expectt{\calE(t) \cdot \nabla^2_{\bx} V_{\mathrm{rel}-} (\bx', T-t) \cdot \calE(t)^\top \ \big| \ \bX(t) = \bx(t), \bW(t) = \bw, \bU(t) = \bu^*_{\bx(t),T+1-t,\bw}(t)}
\end{equation}
Again the crucial part in \eqref{eq:analysis-of-term-C} is that the first-order term cancels out after taking expectation, and we are left with only the second-order term. Combining \eqref{eq:decomposition-of-terms}, \eqref{eq:analysis-of-term-A}, \eqref{eq:definition-of-small-z}, \eqref{eq:difference-of-prio-and-post-refine}, \eqref{eq:analysis-of-term-C} and using Assumption \ref{assumption:C2-smoothness}, we finally obtain
\begin{equation*}
  \abs{z(t) - z(t+1)} \le \frac{1}{2} \norminf{\mathfrak{K}_t} \var{\bW} + \frac{1}{2} \norminf{\mathfrak{L}_t} \var{\calE}
\end{equation*}
The rest of the analysis are similar to Theorem \ref{thm:general-theory-of-resolving-Lipschitz}.   \Halmos

\endproof

\subsection{The Projection Policy}   \hypertarget{subsec:proof-of-projection}{}
For convenience, we repeat the assumptions and theorems below for ease of discussion.

\begin{repeatassumption}[Assumption III]
For all $1 \le t \le T$, there exists $\varepsilon_t > 0$ such that $\proj_t(\bx,\bw)$ as defined in \eqref{eq:proj} is locally Lipschitz-continuous in $\calB((\bx^*(t),\barw), \varepsilon_t)$.
\end{repeatassumption}

\begin{repeattheorem}[Theorem 3 (restated with explicit constants)]  
Let $V_{\mathrm{opt}} (\bx,T)$ be the value of the stochastic optimization problem \eqref{eq:original-problem} that satisfies Assumptions \ref{assumption-1}-\ref{assumption-8}, and let $V_{\mathrm{proj}} (\bx,T)$ be the value of the projection policy defined in Algorithm \ref{algo:projection}. Under the additional Assumption \ref{assumption:lipschitz-continuity-projection}, there exists constants $C_3, \mathfrak{P}', \bar{V} > 0$ such that 
    \begin{equation*}
    V_{\mathrm{opt}} (\bx,T) - V_{\mathrm{proj}} (\bx,T) \le \mathfrak{P}' C_3 +  (1 - \mathfrak{P}') \bar{V}
  \end{equation*}
The constant $\mathfrak{P}'$ is given in \eqref{eq:event-proba-formula-1} and converges to $1$ exponentially fast as both $\var{\bW}$ and $\var{\calE}$ converge to $0$. The constant $\bar{V}$ is a finite upper bound of $V_{\mathrm{opt}} (\bx,T)$. The constant $C_3$ is 
\begin{equation}  \label{eq:choice-of-constant-projection-Lipschitz}
  C_3 \defeq \sum_{t=1}^{T} \left( c_{\bW}(t) \cdot \sqrt{\var{\bW}} + c_{\calE}(t) \cdot \sqrt{\var{\calE}} \right) 
\end{equation}  
where 
\begin{equation*}
  c_{\bW}(t) \defeq c_R (\cproj + 1) + \frac{c_R c_{\phi} (\cproj + 1)^2 ( c_{\phi}^{t-1} (\cproj + 1)^{t-1} - 1)}{ c_{\phi} (\cproj + 1) - 1}
\end{equation*}
\begin{equation*}
  c_{\calE}(t) \defeq \frac{c_R (\cproj + 1) ( c_{\phi}^{t-1} (\cproj + 1)^{t-1} - 1)}{ c_{\phi} (\cproj + 1) - 1}
\end{equation*}
with 
\begin{itemize}
  \item $c_{\phi}$ being the Lipschitz constant of $\phi(\cdot)$
  \item $c_R \defeq \max_t c_{R_t}$ and $c_{R_t}$ being the Lipschitz constant of the reward function $R_t(\cdot)$
  \item $c_{\mathrm{proj}} \defeq \max_t c_{\mathrm{proj},t}$ with $c_{\mathrm{proj},t}$ being the Lipschitz constant of the projection mapping $\proj(\bx,\bw,t)$ inside $\calB((\bx^*(t),\barw), \varepsilon_t)$ claimed in Assumption \ref{assumption:lipschitz-continuity-projection}
\end{itemize}
\end{repeattheorem}

\proof{Proof of Theorem \ref{thm:general-theory-of-projection-Lipschitz}}
The first step of the proof is very similar to the one of Theorem \ref{thm:general-theory-of-resolving-Lipschitz}. Namely, write 
\begin{equation} \label{eq:one-step-deterministic}
  \bx^*(t+1) = \phi(\bx^*(t), \barw, \proj_t(\bx^*(t), \barw)) 
\end{equation}
and 
\begin{equation}  \label{eq:one-step-stochastic}
  \bX(t+1) = \phi(\bX(t), \bW(t), \proj_t(\bX(t), \bW(t)) ) + \calE(t)
\end{equation}
where we have abbreviated $\calE(t) \defeq \calE(\bX(t), \bW(t), \proj_t(\bX(t), \bW(t)) )$. Taking the difference we obtain
\begin{equation} \label{eq:x-relation}
  \norme{\bX(t+1) - \bx^*(t+1)} \le c_{\phi} (c_{\mathrm{proj}}+1) \norme{\bX(t) - \bx^*(t)} + c_{\phi} (c_{\mathrm{proj}}+1) \norme{\bW(t) - \barw} + \norme{\calE(t)}
\end{equation}
where $c_{\phi}$ is the Lipschitz constant of $\phi(\cdot)$, and $c_{\mathrm{proj}} \defeq \max_t c_{\mathrm{proj},t}$, where for $1 \le t \le T$, $c_{\mathrm{proj},t}$ is the Lipschitz constant of the projection mapping $\proj(\bx,\bw,t)$ inside $\calB((\bx^*(t),\barw), \varepsilon_t)$ claimed in Assumption \ref{assumption:lipschitz-continuity-projection}. 

Denote by $\alpha \defeq c_{\phi} (c_{\mathrm{proj}} + 1)$, $\varepsilon \defeq \min_t \varepsilon_t/2$, and 
\begin{equation*}
  \beta \defeq \begin{cases}
        \varepsilon/T, & \mbox{if } \alpha = 1 \\
        \varepsilon(\alpha-1)/\alpha^T, & \mbox{if } \alpha > 1 \\
        \varepsilon (1 - \alpha), & \mbox{if } 0 < \alpha < 1
      \end{cases}
\end{equation*}
An elementary calculation shows that 
\begin{equation*}
  \mbox{$\alpha \norme{\bW(t) - \barw} + \norme{\calE(t)} \le \beta$ holds for all $t$} \Rightarrow \norme{\bX(t) - \bx^*(t)} \le \varepsilon_t/2 \mbox{ holds for all } t
\end{equation*}
Hence we deduce that
\begin{align} 
   & \norme{\calE(t)} \le \beta/2 \mbox{ and } \norme{\bW(t) - \barw} \le \min \{ \beta/(2 \alpha), \varepsilon \} \mbox{ hold for all time-step $t$}  \nonumber \\
  \Rightarrow & \qquad \mbox{the stochastic trajectory $\bX(t)$ together with all the realizations of $\bW(t)$ during the projection} \nonumber \\
   & \qquad \mbox{policy remains inside the $\varepsilon_t$-neighbourhood required by Assumption \ref{assumption:lipschitz-continuity-projection} for all time-step $t$}  \label{eq:event-for-time-t-1}
\end{align}
By Lemma \ref{lem:multivariate-hoeffding}, the event described on the left hand side of \eqref{eq:event-for-time-t-1} occurs with probability at least
\begin{equation} \label{eq:event-proba-formula-1}
 \mathfrak{P}' \defeq \left[1 - \exp \left(- \frac{(\beta/2)^2}{8 \var{\calE}} + \frac{1}{4} \right) \right ]^T  \left[1 - \exp \left(- \frac{(\min \{ \beta/(2\alpha), \varepsilon \})^2}{8 \var{\bW}} + \frac{1}{4} \right) \right ]^T
\end{equation}

We subsequently assume that the right hand side of \eqref{eq:event-for-time-t-1} holds. Taking expectation in \eqref{eq:x-relation} we obtain
\begin{equation}  \label{eq:illustration-3}
  \expectproj{\norme{\bX(t+1) - \bx^*(t+1)}} \le \sqrt{\var{\calE}} + c_{\phi} (\cproj + 1) \sqrt{\var{\bW}} + c_{\phi} (\cproj + 1) \expectproj{\norme{\bX(t) - \bx^*(t)}}
\end{equation}
where $\expectproj{\cdot} = \expectproj{ \ \cdot \mid \bX(1) = \bx}$ is the expectation taken under the projection policy. Using the abbreviation $\mathfrak{h} \defeq \sqrt{\var{\calE}} + c_{\phi} (\cproj + 1) \sqrt{\var{\bW}}$, an elementary calculation implies that for all $t$:
\begin{equation*}
  \expectproj{\norme{\bX(t) - \bx^*(t)}} \le \mathfrak{h} \cdot \frac{( c_{\phi}^{t-1} (\cproj + 1)^{t-1} - 1)}{ c_{\phi} (\cproj + 1) - 1} 
\end{equation*}
where the right hand side is interpreted as $\mathfrak{h}(t-1)$ if $ c_{\phi} (\cproj + 1) = 1$. 

By definition, we have
\begin{equation*}
  V_{\mathrm{rel}-} (\bx,T) = \sum_{t=1}^{T} R_t (\bx^*(t),\barw,\bu^*(t))
\end{equation*}
and
\begin{equation*}
  V_{\mathrm{proj}} (\bx,T) =  \expectproj{\sum_{t=1}^{T} R_t \left( \bX(t), \bW(t), \proj(\bX(t),\bW(t),t) \right)}
\end{equation*}
Let $c_{R_t}$ be the Lipschitz constant of the reward function $R_t(\cdot)$ and write $c_R \defeq \max_t c_{R_t}$. We deduce that 
\begin{align}
  & \abs{V_{\mathrm{rel}-} (\bx,T) - V_{\mathrm{proj}} (\bx,T)}  \nonumber \\
  \le & \ \sum_{t=1}^{T} c_R (\cproj + 1) \expectproj{\norme{\bX(t) - \bx^*(t)}} + c_R (\cproj + 1) \sqrt{\var{\bW}} \label{eq:illustration-4} \\
  \le & \ \sum_{t=1}^{T} c_R (\cproj + 1) \left[ \mathfrak{h} \cdot \frac{( c_{\phi}^{t-1} (\cproj + 1)^{t-1} - 1)}{ c_{\phi} (\cproj + 1) - 1}  +  \sqrt{\var{\bW}} \right] \nonumber
\end{align}
By rearranging terms we achieve the optimality gap bound claimed in the theorem.   \Halmos

\endproof

It is important to note that in Theorem \ref{thm:general-theory-of-projection-Lipschitz}, we were satisfied with a first-order result. Therefore, when establishing upper bounds of the absolute value $\abs{V_{\mathrm{rel}-} (\bx,T) - V_{\mathrm{proj}} (\bx,T)}$, we could conveniently place the norm operator inside the expectation as in \eqref{eq:illustration-4}, and proceed in \eqref{eq:illustration-3} to handle $\expect{\norme{\bX(t) - \bx^*(t)}}$. In contrast, in order to obtain a second-order result, we need to directly estimate $\norme{\bx^*(t) - \expect{\bX(t)}}$, which is smaller than $\expect{\norme{\bX(t) - \bx^*(t)}}$. A similar situation arises in the context of refining a mean field approximation, as noted in \cite[Equations (1) and (2)]{gast2017expected}.

Let us recall Assumption \ref{assumption:C2-smoothness-projection} and restate Theorem \ref{thm:general-theory-of-projection-C2} with explicit constants:

\begin{repeatassumption}[Assumption IV] 
For all $1 \le t \le T$, there exists $\varepsilon_t > 0$ such that $\proj_t(\bx,\bw)$ as defined in \eqref{eq:proj} is locally $\calC^2$-smooth in $\calB((\bx^*(t),\barw), \varepsilon_t)$.
\end{repeatassumption}

\begin{repeattheorem}[Theorem 4 (restated with explicit constants)]  
Let $V_{\mathrm{opt}} (\bx,T)$ be the value of the stochastic optimization problem \eqref{eq:original-problem} that satisfies Assumptions \ref{assumption-1}-\ref{assumption-8}, and let $V_{\mathrm{proj}} (\bx,T)$ be the value of the projection policy defined in Algorithm \ref{algo:projection}. Under the additional Assumption \ref{assumption:C2-smoothness-projection}, there exists constants $C_4, \mathfrak{P}', \bar{V} > 0$ such that 
    \begin{equation*}
    V_{\mathrm{opt}} (\bx,T) - V_{\mathrm{proj}} (\bx,T) \le \mathfrak{P}' C_4 +  (1 - \mathfrak{P}') \bar{V}
  \end{equation*}
The constant $\mathfrak{P}'$ is given in \eqref{eq:event-proba-formula-1} and converges to $1$ exponentially fast as both $\var{\bW}$ and $\var{\calE}$ converge to $0$. The constant $\bar{V}$ is a finite upper bound of $V_{\mathrm{opt}} (\bx,T)$. The constant $C_4$ is 
\begin{equation}  \label{eq:choice-of-constant-projection-C2}
  C_4 \defeq \sum_{t=1}^{T} c_9(t) \cdot \var{\bW} + c_{10}(t) \cdot \var{\calE}
\end{equation}
with the constants $c_9(t)$ and $c_{10}(t)$ given respectively in \eqref{eq:choice-of-constant-projection-C2-1} and \eqref{eq:choice-of-constant-projection-C2-2}, which do not depend on the stochastic part of the system. 

\end{repeattheorem}

The essence of the proof given below can be summarized as follows: For each time-step $t$, when we apply the norm in a second-order Taylor expansion of $\bX(t)$ around $\bx^*(t)$, a term of order $\expect{ \norme{\bx^*(t) - \bX(t)}^2 }$ appears in the Hessian. We first need to bound this term with second-order terms. Upon doing this, it is then used to bind $\norme{\bx^*(t) - \expect{\bX(t)}}$ with second-order terms. These second-order approximations are recursively applied to attain second-order approximations of other quantities in the following sequence: The control $\bU(t) = \proj_t(\cdot)$ is a $\calC^2$-smooth function of $\bX(t)$ and $\bW(t)$ from Assumption \ref{assumption:C2-smoothness-projection}; this subsequently determines the next system configuration $\bX(t+1)$ as a $\calC^2$-smooth function of $\bU(t), \bX(t)$, and $\bW(t)$; overall, the reward function $R_t(\cdot)$ depends $\calC^2$-smoothly on all three entries $\bX(t)$, $\bU(t)$, and $\bW(t)$.

\proof{Proof of Theorem \ref{thm:general-theory-of-projection-C2}}
We divide the proof into multiple steps, and keep an effort to make every constant explicit. To ease the notation, unless otherwise specified, in the following proof we omit the \emph{proj} in the notation of expectations: $\expect{\cdot} = \expectproj{ \ \cdot \mid \bX(1) = \bx}$, so all expectations are understood to be taken under the projection policy.

\paragraph{Step One: Bound $\expect{ \norme{\bx^*(t) - \bX(t)}^2 }$ by Second-Order Terms} \

We construct two sequences $c_1(t) \ge 0$ and $c_2(t) \ge 0$ by induction, such that for all time-step $t$ we have
\begin{equation*}
  \expect{ \norme{\bx^*(t) - \bX(t)}^2 } \le c_1(t) \cdot \var{\bW} + c_2(t) \cdot \var{\calE}
\end{equation*}
For $t=1$, we simply take $c_1(1) = c_2(1) = 0$. Write
\begin{equation} \label{eq:x-relation-1}
  \bX(t+1) = \Phi_t(\bX(t),\bW(t)) + \calE(t)
\end{equation} 
where $\Phi_t(\bx,\bw) \defeq \phi(\bx,\bw, \proj_t(\bx,\bw))$ is $\calC^2$-smooth in $\calB((\bx^*(t),\barw), \varepsilon_t)$ from Assumption \ref{assumption:C2-smoothness-projection}, by shrinking $\varepsilon_t$ if necessary. Write out the second-order Taylor expansion of \eqref{eq:x-relation-1} for each coordinate $1 \le j \le n_x$, by recalling that $\bx^*(t+1) = \Phi_t(\bx^*(t),\barw)$, we obtain
\begin{align}
  X_j(t+1) - x^*_j(t+1) =  & \ \underbrace{\calE_j(t) + \left((\bX(t),\bW(t)) - (\bx^*(t),\barw) \right) \cdot \nabla^\top_{(\bx,\bw)} \Phi_t(\bx^*(t),\barw)_j}_{\mbox{first-order terms}} \nonumber \\
  + & \ \underbrace{ \frac{1}{2} \left((\bX(t),\bW(t)) - (\bx^*(t),\barw) \right) \cdot \nabla^2 \Phi_t(\bx'_j,\bw'_j)_j \cdot \left((\bX(t),\bW(t)) - (\bx^*(t),\barw) \right)^\top}_{\mbox{second-order terms}} \label{eq:diff-of-x}
\end{align}
for a certain vector $(\bx'_j,\bw'_j)$ along the line segment from $(\bx^*(t),\barw)$ to $(\bX(t),\bW(t))$. Denote by $\bGPhi_{t,j,\bx} \defeq \left( \nabla^\top_{(\bx,\bw)} \Phi_t(\bx^*(t),\barw)_j \right)_{\bx}$ the $\bx$-part of the Jacobian $\nabla^\top_{(\bx,\bw)} \Phi_t(\bx^*(t),\barw)_j$, with a similar notation for $\bGPhi_{t,j,\bw}$. Apply $\expect{\norme{\cdot}^2}$ on both sides of \eqref{eq:diff-of-x} and use Cauchy-Schwartz: 
\begin{align*}
  & \expect{\norme{\bX(t+1) - \bx^*(t+1)}^2} = \sum_{j=1}^{n_x} \expect{\left( X_j(t+1) - x^*_j(t+1) \right)^2}  \\
  \le & \ 3 \cdot \sum_{j=1}^{n_x} \Bigg( \expect{\left( \calE_j(t) \right)^2} + \expect{\left( \left((\bX(t),\bW(t)) - (\bx^*(t),\barw) \right) \cdot \nabla^\top_{(\bx,\bw)} \Phi_t(\bx^*(t),\barw)_j \right)^2}  \\
  & \ + \underbrace{\expect{\frac{1}{4} \left( \left((\bX(t),\bW(t)) - (\bx^*(t),\barw) \right) \cdot \nabla^2 \Phi_t(\bx'_j,\bw'_j)_j \cdot \left((\bX(t),\bW(t)) - (\bx^*(t),\barw) \right)^\top  \right)^2}}_{\mbox{fourth-order terms}} \Bigg) \\
  \le & \ 4 \cdot \var{\calE} + 4 \left( \sum_{j=1}^{n_x} \norme{\bGPhi_{t,j,\bx}}^2 \right) \cdot \expect{\norme{\bX(t) - \bx^*(t)}^2} + 4 \left( \sum_{j=1}^{n_x} \norme{\bGPhi_{t,j,\bw}}^2 \right) \cdot \var{\bW}
\end{align*}
where in the last step we have absorbed the fourth-order terms into the second-order terms, hence the coefficient "4". From the induction hypothesis we have subsequently 
\begin{align*}
  \expect{\norme{\bX(t+1) - \bx^*(t+1)}^2} & \le \left( 4 \left( \sum_{j=1}^{n_x} \norme{\bGPhi_{t,j,\bx}}^2 \right) c_1(t) + 4 \left( \sum_{j=1}^{n_x} \norme{\bGPhi_{t,j,\bw}}^2 \right)  \right) \cdot \var{\bW} \\
   & \qquad + \left( 4 + 4 \left( \sum_{j=1}^{n_x} \norme{\bGPhi_{t,j,\bx}}^2 \right) c_2(t)  \right) \cdot \var{\calE}
\end{align*}
So we can define the sequences $c_1(t)$ and $c_2(t)$ recurrently by
\begin{equation}  \label{eq:c-def-1}
  c_1(t+1) = 4 \left( \sum_{j=1}^{n_x} \norme{\bGPhi_{t,j,\bx}}^2 \right) c_1(t) + 4 \left( \sum_{j=1}^{n_x} \norme{\bGPhi_{t,j,\bw}}^2 \right)
\end{equation}
\begin{equation} \label{eq:c-def-2}
  c_2(t+1) = 4 + 4 \left( \sum_{j=1}^{n_x} \norme{\bGPhi_{t,j,\bx}}^2 \right) c_2(t)
\end{equation}
with $c_1(1) = c_2(1) = 0$. 

\paragraph{Step Two: Bound $\expect{\norme{\bu^*(t) - \proj_t(\bX(t),\bW(t))}^2}$ by Second-Order Terms} \

Recall that $\proj_t(\bx^*(t),\barw) = \bu^*(t)$. Using the local $\calC^2$-smoothness of the $\proj_t(\cdot)$ function, for each coordinate $1 \le i \le n_w$, we have the Taylor expansion
\begin{align}
  & u^*_i(t) - \proj_t(\bX(t),\bW(t))_i  = \left[(\bx^*(t), \barw) - (\bX(t), \bW(t)) \right] \cdot \nabla^\top_{(\bx,\bw)} \proj_t(\bx^*(t), \barw)_i  \nonumber \\
  & + \frac{1}{2} \left[(\bx^*(t), \barw) - (\bX(t), \bW(t)) \right] \cdot \nabla^2_{(\bx,\bw)} \proj_t(\bx_j, \bw_j)_i \cdot \left[(\bx^*(t), \barw) - (\bX(t), \bW(t)) \right]^\top  \label{eq:diff-of-u}
\end{align}
for a certain vector $(\bx_i,\bw_i)$ along the line segment from $(\bx^*(t),\barw)$ to $(\bX(t),\bW(t))$, where we have written $\proj_t(\bX(t),\bW(t))_i$ to denote the $i$-th coordinate of $\proj_t(\bX(t),\bW(t))$. We use the abbreviation $ \bGProj_{t,i,\bx} \defeq \left( \nabla^\top_{(\bx,\bw)} \proj(\bx^*(t), \barw,t)_i \right)_{\bx}$ to denote the $\bx$-part of the gradient $\nabla^\top_{(\bx,\bw)} \proj_t(\bx^*(t), \barw)_i$, with a similar notation for $\bGProj_{t,i,\bw}$.

Similar to what we did for the estimation of $\expect{ \norme{\bx^*(t) - \bX(t)}^2 }$, apply $\expect{\norme{\cdot}^2}$ in \eqref{eq:diff-of-u} on both sides and use Cauchy-Schwartz, we obtain:
\begin{equation*}
  \expect{\norme{\bu^*(t) - \proj_t(\bX(t),\bW(t))}^2} \le 3 \left(\sum_{i=1}^{n_u} \norme{\bGProj_{t,i,\bx}}^2 \right) \cdot \expect{\norme{\bx^*(t) - \bX(t)}^2} + 3 \left(\sum_{i=1}^{n_u} \norme{\bGProj_{t,i,\bw}}^2 \right) \cdot \var{\bW}
\end{equation*}
So we can construct two sequences $c_3(t),c_4(t) \ge 0$ such that for all time-step $t$ we have
\begin{equation*}
  \expect{\norme{\bu^*(t) - \proj_t(\bX(t),\bW(t))}^2} \le c_3(t) \cdot \var{\bW} + c_4(t) \cdot \var{\calE}
\end{equation*}
where 
\begin{equation}  \label{eq:c-def-3}
  c_3(t) = 3 \left(\sum_{i=1}^{n_u} \norme{\bGProj_{t,i,\bx}}^2 \right) c_1(t) + 3 \left(\sum_{i=1}^{n_u} \norme{\bGProj_{t,i,\bw}}^2 \right) \mbox{ and } c_4(t) = 3 \left(\sum_{i=1}^{n_u} \norme{\bGProj_{t,i,\bx}}^2 \right) c_2(t)
\end{equation}

\paragraph{Step Three: Bound $\norme{\bu^*(t) - \expect{\proj_t(\bX(t),\bW(t))}}$ and $\norme{\bx^*(t) - \expect{\bX(t)}}$ by Second-Order Terms} \

We use the notation 
\begin{equation*}
  \bHProj_{t,i,\bx \bw} \defeq \left( \nabla^2_{(\bx,\bw)} \proj_t(\bx_i, \bw_i)_i \right)_{\bx \bw}
\end{equation*}
to denote the $\bx\bw$-part for the Hessian matrix of the function $\proj(\cdot)$. And we use 
\begin{equation*}
  \calHProj_{t,i,\bx\bw} \defeq \sup_{(\bx_i,\bw_i) \in \calB((\bx^*(t),\barw), \varepsilon_t)} \norme{ \left( \nabla^2_{(\bx,\bw)} \proj_t(\bx_i, \bw_i)_i \right)_{\bx\bw} }
\end{equation*}
to denote an upper-bound to the norm of $\bHProj_{t,i,\bx\bw}$ in the neighbourhood $\calB((\bx^*(t),\barw), \varepsilon_t)$. The Hessian part of \eqref{eq:diff-of-u} can be written more explicitly as 
\begin{align*}
  & \left[(\bx^*(t), \barw) - (\bX(t), \bW(t)) \right] \cdot \nabla^2_{(\bx,\bw)} \proj_t(\bx_i, \bw_i)_i \cdot \left[(\bx^*(t), \barw) - (\bX(t), \bW(t)) \right]^\top    \\
  = & \ (\bx^*(t) - \bX(t)) \cdot \bHProj_{t,i, \bx \bx} \cdot (\bx^*(t) - \bX(t))^\top + (\barw - \bW(t)) \cdot \bHProj_{t,i, \bw \bx} \cdot (\bx^*(t) - \bX(t))^\top \\
  & \qquad + (\bx^*(t) - \bX(t)) \cdot \bHProj_{t,i, \bx \bw} \cdot (\barw - \bW(t))^\top + (\barw - \bW(t)) \cdot \bHProj_{t,i, \bw \bw} \cdot (\barw - \bW(t))^\top
\end{align*}
Notice that the random vectors $\bW(t)$ and $\bX(t)$ are independent, hence
\begin{equation*}
  \expect{(\barw - \bW(t)) \cdot \bHProj_{t,i, \bw \bx} \cdot (\bx^*(t) - \bX(t))^\top} = \expect{(\barw - \bW(t))} \cdot \bHProj_{t,i, \bw \bx} \cdot \expect{(\bx^*(t) - \bX(t))^\top} = \bzero
\end{equation*}
So finally if we take expectation in \eqref{eq:diff-of-u} and then take the norm, we obtain the bound
\begin{align*}
  & \abs{u^*_i(t) - \expect{\proj_t(\bX(t),\bW(t))}_i} \\
  & \le \abs{\left(\bx^*(t) - \expect{\bX(t)} \right) \cdot \bGProj_{t,i,\bx}} + \frac{1}{2} \left( \calHProj_{t,i, \bx \bx} \cdot \expect{ \norme{\bx^*(t) - \bX(t)}^2 }  + \calHProj_{t,i, \bw \bw} \cdot \var{\bW} \right)
\end{align*}
From which we deduce that
\begin{align}
  & \norme{\bu^*(t) - \expect{\proj_t(\bX(t),\bW(t))}} \le \sum_{i=1}^{n_u} \abs{u^*_i(t) - \expect{\proj_t(\bX(t),\bW(t))}_i} \ \mbox{ ($\calL^2$-norm is smaller than $\calL^1$-norm) } \nonumber \\
  \le & \ \sum_{i=1}^{n_u} \bigg( \norme{\bGProj_{t,i,\bx}} \cdot \norme{\bx^*(t) - \expect{\bX(t)}} + \frac{1}{2} \left( \calHProj_{t,i,\bx \bx} \cdot \expect{ \norme{\bx^*(t) - \bX(t)}^2 } + \calHProj_{t,i, \bw \bw} \cdot \var{\bW} \right) \bigg) \label{eq:etimation-of-u}
\end{align}
Apply the same procedure for $\bX(t)$, starting from \eqref{eq:diff-of-x}, we obtain 
\begin{align}
  & \norme{\bx^*(t+1) - \expect{\bX(t+1)}} \le \sum_{j=1}^{n_x} \abs{x^*_j(t+1) - \expect{X_j(t+1)}}  \nonumber \\
  \le & \sum_{j=1}^{n_x} \bigg( \norme{\bGPhi_{t,j,\bx}} \cdot \norme{\bx^*(t) - \expect{\bX(t)}} + \frac{1}{2} \left( \calHPhi_{t,j,\bx \bx} \cdot \expect{ \norme{\bx^*(t) - \bX(t)}^2 } + \calHPhi_{t,j, \bw \bw} \cdot \var{\bW} \right) \bigg)  \label{eq:estimation-of-x}
\end{align}
We remark that previously we have applied $\expect{\norme{\cdot}^2}$ in \eqref{eq:diff-of-x} to obtain a second-order estimation for $\norme{\bx^*(t) - \expect{\bX(t)}}^2$, whereas now we apply $\norme{\expect{\cdot}}$ in both \eqref{eq:diff-of-x} and \eqref{eq:diff-of-u} to obtain second-order estimation of $\norme{\bu^*(t) - \expect{\proj_t(\bX(t),\bW(t))}}$ and $\norme{\bx^*(t) - \expect{\bX(t)}}$. This can not be done without a second-order estimation of $\norme{\bx^*(t) - \expect{\bX(t)}}^2$ at first hand.

We next apply $\expect{ \norme{\bx^*(t) - \bX(t)}^2 } \le c_1(t) \cdot \var{\bW} + c_2(t) \cdot \var{\calE}$ in both \eqref{eq:etimation-of-u} and \eqref{eq:estimation-of-x}, rearrange terms to obtain
\begin{align}
  & \norme{\bu^*(t) - \expect{\proj_t(\bX(t),\bW(t))}} \le \left( \sum_{i=1}^{n_u} \norme{\bGProj_{t,i,\bx}} \right) \cdot \norme{\bx^*(t) - \expect{\bX(t)}} \nonumber \\
  & + \left( \frac{1}{2} \sum_{i=1}^{n_u} \left( \calHProj_{t,i,\bw\bw} + c_1(t) \cdot \calHProj_{t,i,\bx\bx} \right) \right) \cdot \var{\bW}  + \left( \frac{1}{2} \sum_{i=1}^{n_w} c_2(t) \cdot \calHProj_{t,i,\bx\bx} \right) \cdot \var{\calE}    \label{eq:estimation-of-u-2}
\end{align}
\begin{align}
  & \norme{\bx^*(t+1) - \expect{\bX(t+1)}} \le \left( \sum_{j=1}^{n_x} \norme{\bGPhi_{t,j,\bx}} \right) \cdot \norme{\bx^*(t) - \expect{\bX(t)}} \nonumber\\
  & \qquad + \left( \frac{1}{2} \sum_{j=1}^{n_x} \left( \calHPhi_{t,j,\bw\bw} + c_1(t) \cdot \calHPhi_{t,j,\bx\bx} \right) \right) \cdot \var{\bW}  + \left( \frac{1}{2} \sum_{j=1}^{n_x} c_2(t) \cdot \calHPhi_{t,j,\bx\bx} \right) \cdot \var{\calE}   \label{eq:estimation-of-x-2}
\end{align}
We now apply a similar technique to construct four sequences $c_5(t), c_6(t), c_7(t), c_8(t) \ge 0$ by induction, such that for all time-step $t$ we have
\begin{equation*}
  \norme{\bx^*(t) - \expect{\bX(t)}} \le c_5(t) \cdot \var{\bW} + c_6(t) \cdot \var{\calE}
\end{equation*}
\begin{equation*}
  \norme{\bu^*(t) - \expect{\proj_t(\bX(t),\bW(t))}} \le  c_7(t) \cdot \var{\bW} + c_8(t) \cdot \var{\calE}
\end{equation*}
Their initial values at $t=1$ are 
\begin{equation*}
  c_5(1) = c_6(1) = c_8(1) = 0, \ c_7(1) = \frac{1}{2} \sum_{i=1}^{n_u} \calHProj_{1,i,\bw\bw} 
\end{equation*}
From \eqref{eq:estimation-of-x-2} the recurrent relations to uniquely define $c_5(t)$ and $c_6(t)$ are 
\begin{equation} \label{eq:c-def-4}
  c_5(t+1) = c_5(t) \cdot \left( \sum_{j=1}^{n_x} \norme{\bGPhi_{t,j,\bx}} \right) + \frac{1}{2} \sum_{j=1}^{n_x} \left( \calHPhi_{t,j,\bw\bw} + c_1(t) \cdot \calHPhi_{t,j,\bx\bx} \right) 
\end{equation}
\begin{equation}  \label{eq:c-def-5}
  c_6(t+1) = c_6(t) \cdot \left( \sum_{j=1}^{n_x} \norme{\bGPhi_{t,j,\bx}} \right) + \frac{1}{2} \sum_{j=1}^{n_x} c_2(t) \cdot \calHPhi_{t,j,\bx\bx}
\end{equation}
And from \eqref{eq:estimation-of-u-2}, $c_7(t)$ and $c_8(t)$ are given by 
\begin{equation}  \label{eq:c-def-6}
  c_7(t) = c_5(t) \cdot \left( \sum_{i=1}^{n_u} \norme{\bGProj_{t,i,\bx}} \right) + \frac{1}{2} \sum_{i=1}^{n_u} \left( \calHProj_{t,i,\bw\bw} + c_1(t) \cdot \calHProj_{t,i,\bx\bx} \right)
\end{equation}
\begin{equation}  \label{eq:c-def-7}
  c_8(t) = c_6(t) \cdot \left( \sum_{i=1}^{n_u} \norme{\bGProj_{t,i,\bx}} \right) + \frac{1}{2} \sum_{i=1}^{n_w} c_2(t) \cdot \calHProj_{t,i,\bx\bx}
\end{equation}

\paragraph{Step Four: Conclusion of the Proof} \

We begin by writing
\begin{equation}  \label{eq:basic-relation-1}
  V_{\mathrm{rel}-} (\bx,T) - V_{\mathrm{proj}} (\bx,T) = \expect{\sum_{t=1}^{T} \expect{R_t (\bv^*(t)) - R_t(\bV(t)) \ \Big| \ \bX(t) } }
\end{equation}
where we have abbreviated $\bv^*(t) \defeq (\bx^*(t),\barw,\bu^*(t))$ and $\bV(t) \defeq (\bX(t), \bW(t), \proj_t(\bX(t),\bW(t)) )$. Using the smoothness of the reward function $R_t(\cdot)$, we have
\begin{equation} \label{eq:basic-relation-2}
    R_t (\bv^*(t)) - R_t (\bV(t) ) =   \left[ \bv^*(t) - \bV(t) \right] \cdot \nabla^\top_{(\bx,\bw,\bu)} R_t (\bv^*(t)) + \frac{1}{2} \left[ \bv^*(t) - \bV(t) \right] \cdot \nabla^2_{(\bx,\bw,\bu)} R_t (\bv(t)) \cdot \left[ \bv^*(t) - \bV(t) \right]^\top
\end{equation}
with $\bv(t)$ being a certain vector along the line segment from $\bv^*(t)$ to $\bV(t)$. In the following we use the abbreviation $\bU(t) \defeq \proj_t(\bX(t),\bW(t))$. The meaning for $\bGR_{t,\bx}$, $\bHR_{t,\bx\bx}$, $\calHR_{t,\bx\bx}$ is analogue to what we have used before for $\Phi(\cdot)$ and $\proj(\cdot)$. Combine everything we have obtained so far, we deduce that  
\begin{align*}
 & \abs{V_{\mathrm{rel}-} (\bx,T) - V_{\mathrm{proj}} (\bx,T)} \le \sum_{t=1}^{T} \bigg( \norme{\bx^*(t) - \expect{\bX(t)}} \norme{\bGR_{t,\bx}} + \norme{\bu^*(t) - \expect{\bU(t)}} \norme{\bGR_{t,\bu}} \\
 & \qquad + \frac{1}{2} \calHR_{t,\bx\bx} \cdot \expect{\norme{\bx^*(t) - \bX(t)}^2}  + \frac{1}{2} \calHR_{t,\bw\bw} \cdot \var{\bW}  + \frac{1}{2} \calHR_{t,\bu\bu} \cdot \expect{\norme{\bu^*(t) - \bU(t)}^2}  \\
 & \qquad + \expect{\norme{(\bx^*(t) - \bX(t)) \cdot \bHR_{t,\bx\bw} \cdot (\barw - \bW(t))^\top }} + \expect{\norme{(\bx^*(t) - \bX(t)) \cdot \bHR_{t,\bx\bu} \cdot (\bu^*(t) - \bU(t))^\top }} \\
 & \qquad + \expect{\norme{(\bu^*(t) - \bU(t)) \cdot \bHR_{t,\bu\bw} \cdot (\barw - \bW(t))^\top }} \bigg) \\
 & \ \le \sum_{t=1}^{T} \bigg( \norme{\bx^*(t) - \expect{\bX(t)}} \norme{\bGR_{t,\bx}} + \norme{\bu^*(t) - \expect{\bU(t)}} \norme{\bGR_{t,\bu}} \\
 & \qquad + \frac{1}{2} (\calHR_{t,\bx\bx}+\calHR_{t,\bx\bw}+\calHR_{t,\bx\bu}) \cdot \expect{\norme{\bx^*(t) - \bX(t)}^2}  + \frac{1}{2} (\calHR_{t,\bw\bx}+\calHR_{t,\bw\bw}+\calHR_{t,\bw\bu}) \cdot \var{\bW} \\
 & \qquad + \frac{1}{2} (\calHR_{t,\bu\bx}+\calHR_{t,\bu\bw}+\calHR_{t,\bu\bu}) \cdot \expect{\norme{\bu^*(t) - \bU(t)}^2}  \bigg) 
\end{align*}
\begin{align*} 
 & \ \le \sum_{t=1}^{T} \bigg( \Big( c_5(t) \norme{\bGR_{t,\bx}} + c_7(t) \norme{\bGR_{t,\bu}} + \frac{1}{2} c_1(t) (\calHR_{t,\bx\bx}+\calHR_{t,\bx\bw}+\calHR_{t,\bx\bu}) + \frac{1}{2} (\calHR_{t,\bw\bx}+\calHR_{t,\bw\bw}+\calHR_{t,\bw\bu}) \\
 & \qquad + \frac{1}{2} c_3(t) (\calHR_{t,\bu\bx}+\calHR_{t,\bu\bw}+\calHR_{t,\bu\bu})   \Big) \cdot \var{\bW}  \\
 & +  \left( c_6(t) \norme{\bGR_{t,\bx}} +  c_8(t) \norme{\bGR_{t,\bu}} + \frac{1}{2} c_2(t) (\calHR_{t,\bx\bx}+\calHR_{t,\bx\bw}+\calHR_{t,\bx\bu}) + \frac{1}{2} c_4(t) (\calHR_{t,\bu\bx}+\calHR_{t,\bu\bw}+\calHR_{t,\bu\bu})  \right)  \cdot \var{\calE}  \bigg)  \\
 & \ = \sum_{t=1}^{T} c_9(t) \cdot \var{\bW} + c_{10}(t) \cdot \var{\calE}
\end{align*}
where we have applied Cauchy-Schwartz in the second step. The constants $c_9(t)$ and $c_{10}(t)$ are given by 
\begin{align} 
  c_9(t) & \defeq c_5(t) \norme{\bGR_{t,\bx}} + c_7(t) \norme{\bGR_{t,\bu}} \nonumber \\
   & \qquad + \frac{1}{2} \left( c_1(t) (\calHR_{t,\bx\bx}+\calHR_{t,\bx\bw}+\calHR_{t,\bx\bu}) + (\calHR_{t,\bw\bx}+\calHR_{t,\bw\bw}+\calHR_{t,\bw\bu}) + c_3(t) (\calHR_{t,\bu\bx}+\calHR_{t,\bu\bw}+\calHR_{t,\bu\bu}) \right) \label{eq:choice-of-constant-projection-C2-1}
\end{align}  
\begin{equation} \label{eq:choice-of-constant-projection-C2-2}
  c_{10}(t) \defeq c_6(t) \norme{\bGR_{t,\bx}} +  c_8(t) \norme{\bGR_{t,\bu}} + \frac{1}{2} \left( c_2(t) (\calHR_{t,\bx\bx}+\calHR_{t,\bx\bw}+\calHR_{t,\bx\bu}) + c_4(t) (\calHR_{t,\bu\bx}+\calHR_{t,\bu\bw}+\calHR_{t,\bu\bu}) \right)
\end{equation} 
where the sequences $c_1(t)$ and $c_2(t)$ are defined in \eqref{eq:c-def-1} and \eqref{eq:c-def-2}; $c_3(t)$ and $c_4(t)$ are defined in \eqref{eq:c-def-3}; $c_5(t)$ and $c_6(t)$ are defined in \eqref{eq:c-def-4} and \eqref{eq:c-def-5}; $c_7(t)$ and $c_8(t)$ are defined in \eqref{eq:c-def-6} and \eqref{eq:c-def-7}. This concludes the proof of the theorem.   \Halmos

\endproof

\section{Policy Mapping and Sensitivity Analysis}     \hypertarget{appendix:policy-mapping}{}

In this appendix we provide sufficient conditions for satisfying Assumptions \ref{assumption:lipschitz-continuity}, \ref{assumption:C2-smoothness} for the update policy, and Assumptions \ref{assumption:lipschitz-continuity-projection}, \ref{assumption:C2-smoothness-projection} for the projection policy, that transform into studying the parameterized solution mapping to a mathematical program and their sensitivity analysis. To give a unifying and general discussion, we first formulate some background definition that pertain to variational inequalities (VI) in Section \hyperlink{subsec:policy-mapping-preliminary}{EC.3.1}. The sufficient conditions for satisfying the local $\calC^2$-smoothness Assumptions \ref{assumption:C2-smoothness} and \ref{assumption:C2-smoothness-projection} are discussed in Section \hyperlink{subsec:non-degeneracy-and-strict-complementarity}{EC.3.2}. The sufficient conditions for satisfying the local Lipschitz continuity Assumptions \ref{assumption:lipschitz-continuity} and \ref{assumption:lipschitz-continuity-projection} are discussed in Section \hyperlink{subsec:lipschitz-property}{EC.3.3}. Further properties related to the Euclidean projector and non-degeneracy are given in Section \hyperlink{subsec:euclidean-projector}{EC.3.4}. We refer to \cite{facchinei2003finite} for omitted arguments in our discussion. 

\subsection{Preliminaries}  \hypertarget{subsec:policy-mapping-preliminary}{}

Given a closed, convex and non-empty subset $K$ of the Euclidean $n$-dimensional space $\R^n$ and a mapping $F: K \rightarrow \R^n$, the \emph{variational inequality}, denoted as VI($K,F$), is the problem of finding the set of vectors $\bx \in K$ such that $(\by - \bx) \cdot F(\bx)^\top \ge 0$ hold true for all $\by \in K$. We write SOL($K,F$) for the solution set to the corresponding VI($K,F$). Denote by $\calN_{K}(\bx')$ the \emph{normal cone} to $K$ at $\bx'$, i.e. $\calN_{K}(\bx') \defeq \left\{ \bd \in \R^n \mid \bd \cdot (\by - \bx')^\top \le 0, \ \forall \by \in K \right\}$, then $\bx$ solves VI($K,F$) if and only if $\mathbf{0} \in F(\bx) + \calN_{K}(\bx)$ \cite[Section 1.1]{facchinei2003finite}. Suppose that $K$ is finitely representable and is given by the following set of equations and inequalities:
\begin{equation*}
  K \defeq \left\{ \bx \in \R^n \ \big | \ h_j(\bx) = 0, \forall j = 1,\dots,J \mbox{ and } g_i(\bx) \le 0, \forall i = 1,\dots,I \right\}
\end{equation*}
where $h_j: \R^n \rightarrow \R$ are affine and $g_i: \R^n \rightarrow \R$ are convex and $\calC^2$-smooth. Consider a convex function $f: \R^n \rightarrow \R$ that is defined and at least $\calC^1$ on an open set that contains $K$. Then minimize $f(\cdot)$ in $K$ is equivalent to finding $\bx$ in SOL($K,F$) with $F = \nabla f$. A function $F$ that can be written as the gradient of another function $f$ is called \emph{integrable}. The variational inequality VI($K,F$) generalizes the constrained nonlinear programming by considering a general function $F$ not necessarily being integrable \cite[Section 1.3.1]{facchinei2003finite}.

We next introduce an additional perturbation parameter space $\calP$ for the VI, so that the parameterized problem is formulated with $F: \calD \times \calP \subset \R^n \times \R^p \rightarrow \R^n$ being a function of two arguments $(\bx,\bp)$, and $K: \R^p \rightarrow \calD$ is a multi-function with values in $\calD$:
\begin{equation}  \label{eq:K-definition}
  K(\bp) \defeq \left\{ \bx \in \R^n \ \big | \ h_j(\bx, \bp) = 0, \forall j = 1,\dots,J \mbox{ and } g_i(\bx, \bp) \le 0, \forall i = 1,\dots,I \right\}
\end{equation}
where $h_j: \R^n \times \R^p \rightarrow \R$ are affine and $g_i: \R^n \times \R^p \rightarrow \R$ are convex and $\calC^2$-smooth, jointly in $(\bx,\bp)$. We are given a solution $\bx^*$ of the VI($K(\bp^*),F(\cdot, \bp^*)$), and we aims at understanding the properties of the solution mapping $\bx(\bp) \in \mbox{SOL}(K(\bp),F(\cdot, \bp))$ for $\bp$ near $\bp^*$. 

Under the problem setting of this paper, the VI's of interest for the update policy are the convex programs \eqref{eq:original-problem-res} for $1 \le t \le T$, parameterized by $\bp = (\bx,\bw) \in \R^{n_x} \times \R^{n_w}$, and the solution mapping $\bx(\bp)$ is the policy mapping 
$S^*_t(\bx,\bw)$ defined in \eqref{eq:policy-mapping} for $(\bx,\bw)$ in a neighbourhood of $(\bx^*(t),\barw)$. The VI's of interest for the projection policy are those prescribed by \eqref{eq:feasible-set-for-projection} and \eqref{eq:convex-program-for-projection} for $1 \le t \le T$, parameterized by $\bp = (\bx,\bw) \in \R^{n_x} \times \R^{n_w}$, and the solution mapping $\bx(\bp)$ is the policy mapping $\proj_t(\bx,\bw)$.

\subsection{The Non-Degenerate Property and LICQ}    \hypertarget{subsec:non-degeneracy-and-strict-complementarity}{}

\begin{figure}[htbp]
  \centering
  \input{figure/degeneracy1.tex}
  \caption{Illustration of various concepts: normal cone, tangent cone, critical cone. Also, note that $\bx_1$ is non-degenerate while $\bx_2$ is degenerate.} \label{fig:degeneracy1}
\end{figure}
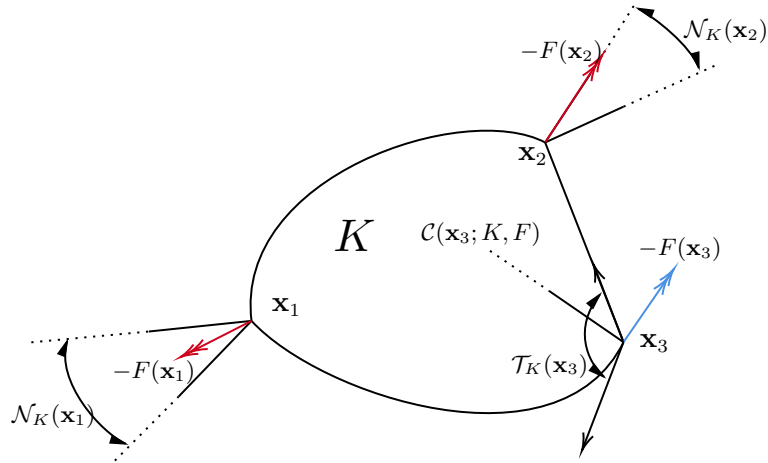

Consider a non-parameterized VI($K,F$). If $\bx \in \mbox{SOL}(K,F)$, then $-F(\bx) \in \calN_K(\bx)$. We call a solution $\bx \in \mbox{SOL}(K,F)$ \emph{non-degenerate} if $-F(\bx)$ is in the \emph{relative interior} of $\calN_K(\bx)$. Otherwise it is called \emph{degenerate}. See Figure \ref{fig:degeneracy1} for an illustration. By definition, the LICQ (Linear Independent Constraint Qualification) holds at $\bx \in K$ if 
\begin{equation*}
  \left\{ \nabla g_i(\bx) \mid i \in \calI(\bx) \right\} \cup \left\{ \nabla h_j(\bx) \mid j = 1,\dots,J \right\}
\end{equation*}
are linearly independent, where $\calI(\bx) \defeq \left\{ i \mid g_i(\bx)=0 \right\}$ is the set of active constraints at $\bx$. 

\begin{theorem} [Sufficient Conditions for Satisfying Assumptions \ref{assumption:C2-smoothness} and \ref{assumption:C2-smoothness-projection}] \hypertarget{thm:sufficient-conditions-C2}{} 

Let $\bx^* \in \emph{SOL}(K(\bp^*),F(\cdot, \bp^*))$ satisfies the LICQ and is non-degenerate, where $K(\bp)$ is given by \eqref{eq:K-definition} and $F(\bx, \bp) = \nabla_{\bx} f(\bx,\bp)$ is $\calC^2$-smooth jointly in $(\bx,\bp)$. Then there exists a neighbourhood $\calX$ of $\bx^*$, a neighbourhood $\calW$ of $\bp^*$ and a $\calC^2$-smooth function $\bx: \calW \rightarrow \calX$, such that $\bx(\bp^*) = \bx^*$. And for every $\bp \in \calW$, $\bx(\bp)$ is the only solution in \emph{SOL}($K(\bp),F(\cdot,\bp)$)
\end{theorem}
We give two references for a complete proof of the above Theorem \hyperlink{thm:sufficient-conditions-C2}{EC.1}. The first is \cite[Theorem 5.4.15]{facchinei2003finite}. Note that in this proof $F(\cdot)$ is only assumed to be $\calC^1$-smooth, and no assumption of integrability is needed. Subsequently the solution mapping $\bx(\cdot)$ is only shown to be locally $\calC^1$-smooth. Another proof that is more suited for our purpose is \cite[Theorem 3.3]{robinson1987local}, for which we discuss below. In this proof, $F(\cdot)$ is assumed to be $\calC^r$-smooth for $r \ge 1$. As a consequence we deduce that $\bx(\cdot)$ is locally $\calC^r$-smooth. In other words, the solution mapping $\bx(\cdot)$ inherits the same level of smoothness as the function $F(\cdot)$ under the setting of the theorem. 

To prepare the proof, we first give several equivalent characterization of the non-degenerate property. The \emph{tangent cone} $\calT_K(\bx)$ of $K$ at $\bx$ consists of all vectors $\bd \in \R^n$ for which there exists a sequence of vectors $\{\by_k\} \subset K$ and a sequence of positive scalars $\{\tau_k\}$ such that 
\begin{equation*}
  \lim_{k \rightarrow \infty} \by_k = \bx, \ \lim_{k \rightarrow \infty} \tau_k = 0, \mbox{ and } \lim_{k \rightarrow \infty} \frac{\by_k - \bx}{\tau_k} = \bd
\end{equation*}
By definition, the Abadie's CQ postulates that $\calT_K(\bx)$ is equal to the \emph{linearization cone} of $K$ at $\bx$, defined as 
\begin{equation*}
  \calL_K (\bx) \defeq \{ \bv \in \R^n \mid \nabla h_j(\bx) \cdot \bv^\top = 0 \ \forall j=1,\dots,J \mbox{ and } \nabla g_i(\bx) \cdot \bv^\top \le 0 \ \forall i \in \calI(\bx) \}
\end{equation*}
where $\calI(\bx) \defeq \left\{ i \mid g_i(\bx)=0 \right\}$ is the set of active constraints at $\bx$. The \emph{critical cone} of the pair $(K,F)$ at $\bx \in K$ is defined as $\calC(\bx; K,F) \defeq \calT_K(\bx) \cap F(\bx)^{\perp}$. The \emph{Karush-Kuhn-Tucker} (KKT) system associated with the VI($K,F$) is the problem of finding $(\bx,\mu,\lambda)$ that satisfies the following:
\begin{equation*}
  0 = F(\bx) + \sum_{j=1}^{J} \mu_j \nabla h_j(\bx) + \sum_{i=1}^{I} \lambda_i \nabla g_i(\bx) 
\end{equation*}
\begin{equation}  \label{eq:KKT}
  \bzero = h(\bx) \mbox{ and } \bzero \le \lambda \perp g(\bx) \le \bzero
\end{equation}
A solution $(\bx,\mu,\lambda)$ to the KKT system is called a KKT triple, and $(\mu,\lambda)$ is called the KKT multiplier. Under Abadie's CQ, there exists KKT multiplier $(\mu,\lambda)$ so that $(\bx,\mu,\lambda)$ is a KKT triple if and only if $\bx \in \mbox{SOL}(K,F)$ \cite[Proposition 1.3.4]{facchinei2003finite}. The equivalent characterizations of the non-degenerate property are:
\begin{enumerate} [label=(\roman*)]
  \item $\bx \in \mbox{SOL}(K,F)$ is non-degenerate if and only if the critical cone $\calC(\bx;K,F)$ is a linear subspace, i.e. $\calC(\bx;K,F) = \calC(\bx;K,F) \cap - \calC(\bx;K,F)$. And it must equal the linearity space of the tangent cone $\calT_K(\bx)$, i.e. the largest linear subspace contained in $\calT_K(\bx)$ \cite[Proposition 3.4.2]{facchinei2003finite}.
  \item Suppose that $\bx \in \mbox{SOL}(K,F)$ at which the Abadie's CQ holds. Then $\bx$ is non-degenerate if and only if there exists KKT triple $(\bx,\mu,\lambda)$ such that $\lambda - g(\bx) > \bzero$. The later is referred to as \emph{strict complementarity} \cite[Corollary 3.4.3]{facchinei2003finite}.
\end{enumerate}
In Figure \ref{fig:degeneracy1}, $\bx_1 \in \mbox{SOL}(K,F)$ is non-degenerate since the critical cone $\calC(\bx_1;K,F)$ is the singleton $\bx_1$, a $0$-dimensional linear subspace; $\bx_2 \in \mbox{SOL}(K,F)$ is degenerate since the critical cone $\calC(\bx_2;K,F)$ is a half-line. To avoid confusion, we mention that there is an unfortunate conflict of terminology in the literature, as Robinson in \cite{robinson1982generalized, robinson1987local} has attributed a property that generalizes LICQ the name "non-degeneracy", while the non-degeneracy we consider here has also been studied in \cite{robinson1982generalized, robinson1987local} without giving an explicit name. The more common terminology of non-degeneracy for which we follow in this paper is due to \cite{dunn1987convergence}. See \cite[Section 3.8]{facchinei2003finite} for more discussion on this issue.

\proof{Proof of Theorem \hyperlink{thm:sufficient-conditions-C2}{EC.1}, sketched from \cite{robinson1987local} }
  Let $\bx^* \in \mbox{SOL}(K(\bp^*),F(\cdot, \bp^*))$ satisfies the LICQ and is non-degenerate. From LICQ, for each $\bp \in \calW$, we can construct a $\calC^r$-diffeomorphism between VI($K(\bp),F(\cdot,\bp)$) and VI($\calC(\bx^*;K,F),\tilde{F}(\cdot,\bp)$), where $\tilde{F}(\cdot,\bp)$ is a $\calC^r$-smooth function defined in a neighbourhood of the origin of $\R^n$. The noticeable feature of the latter is that the constraint set $\calC(\bx^*;K,F)$ is no longer parameterized by $\bp$. Furthermore, from non-degeneracy, $\calC(\bx^*;K,F)$ is a linear subspace. So that VI($\calC(\bx^*;K,F),\tilde{F}(\cdot,\bp)$) is equivalent to optimizing $\tilde{F}(\cdot,\bp)$ in an unconstrained set $\R^m$, with $m$ being the dimension of the subspace $\calC(\bx^*;K,F)$. We may now apply the classical Implicit Function Theorem (IFT) to show that the solution mapping $\tilde{\bx}(\bp)$ of VI($\calC(\bx^*;K,F),\tilde{F}(\cdot,\bp)$) is $\calC^r$-smooth in $\calW$. Using the $\calC^r$-diffeomorphism, we transform $\tilde{\bx}(\bp)$ back into the solution mapping $\bx(\bp) \in \mbox{SOL}(K(\bp),F(\cdot,\bp))$, and conclude that $\bx(\bp)$ is locally $\calC^r$-smooth in $\calW$.     \Halmos
\endproof
To summarize the ideas in this proof: 
\begin{quote}
  Under LICQ and non-degeneracy, a parameterized constrained optimization problem behaves like a non-parametric (from LICQ) and unconstrained (from non-degeneracy) problem in a neighbourhood of the stationary point in question.
\end{quote}
This observation has played a key role in our refined analysis of Theorem \ref{thm:general-theory-of-resolving-C2} and \ref{thm:general-theory-of-projection-C2}, which guarantees respectively that the policy mappings $S^*_t(\bx,\bw)$ and $\mathrm{Proj}_t(\bx,\bw)$ 
inherit the same level of smoothness of the corresponding objective functions. In general, explicit computation of a such solution mapping is not possible, as already can be seen from the classical IFT, which guarantees the existence but not a construction of the implicit function. Nevertheless, for the particular case of WCMDP in Example \ref{example-2}, where the relaxed problems are linear programs, such solution mapping can be efficiently computed via inverting some coefficient matrix. This has played an important role in the implementation of an improved update policy for WCMDP depicted in Algorithm $2$ of \cite{gast2022lp}.

It is important to note that both the LICQ and non-degeneracy are seen as generic properties for a \emph{single} problem. While a rigorous formulation of this claim is not provided in this paper, we refer to a result due to Morse in \cite[Theorem 1.48]{ioffe2017variational} and related concepts therein for a discussion. However, the conditions required for \eqref{eq:original-problem-rel} to meet Assumptions \ref{assumption:C2-smoothness} and \ref{assumption:C2-smoothness-projection} are notably more stringent: they necessitate that the LICQ and non-degeneracy conditions be simultaneously satisfied in a total of $T$ sub-programs, all deduced from an optimal solution to the original program \eqref{eq:original-problem-rel}. In instances where these conditions are not met, we look to establish weaker conditions ensuring Lipschitz continuity, as encapsulated in Assumptions \ref{assumption:lipschitz-continuity} and \ref{assumption:lipschitz-continuity-projection}, which we will explore next.

\subsection{The Lipschitz Property}  \hypertarget{subsec:lipschitz-property}{}

In the literature, the Lipschitz property to the solution mapping is often referred to as a sort of solution stability \cite[Definitions 5.2.3 and 5.2.6]{facchinei2003finite} or regularity \cite[Theorem 5.4.12]{facchinei2003finite}, with the Lipschitz constant measuring the speed of change to SOL($K,F$) when the pair $(K,F)$ undergoes small perturbations. We first state two well known constraint qualifications that are weakened from LICQ. 
\begin{enumerate}[label=\alph*)]
\item The Mangasarian-Fromovitz Constraint Qualification (MFCQ) holds at $\bx \in K$, if
\begin{enumerate}[label=(\roman*)]
  \item The gradients $\left\{ \nabla h_j(\bx), \forall j = 1,\dots,J  \right\}$ are linearly independent.
  \item There exists $\bv \in \R^n$ such that $\nabla h_j(\bx) \cdot \bv^\top = 0$, for all $j=1,\dots,J$, and $\nabla g_i(\bx) \cdot \bv^\top < 0$, for all $i \in \calI(\bx)$. 
\end{enumerate}
\item The Constant Rank Constraint Qualification (CRCQ) holds at $\bx \in K$ if there exists a neighbourhood $\calB(\bx,\varepsilon)$ of $\bx$ such that for every pair of index subsets $\calI' \subset \calI(\bx)$ and $\calJ' \subset \{1,\dots,J \}$, the family of gradient vectors 
\begin{equation*}
  \left\{ \nabla g_i(\bx') \mid i \in \calI' \right\} \cup  \left\{ \nabla h_j(\bx') \mid j \in \calJ' \right\}
\end{equation*}
has the same rank for all $\bx' \in \calB(\bx,\varepsilon) \cap K$ (which depends on $(\calI',\calJ')$). 
\end{enumerate}
It can be shown that LICQ implies both MFCQ and CRCQ, while MFCQ and CRCQ can not be compared with each other. Also any one of these three CQ's implies Abadie's CQ \cite[Proposition 3.2.1]{facchinei2003finite}. 

It turns out that to ensure the Lipschitz property, with MFCQ and CRCQ together are not sufficient, and we need to find a third condition in the absence of non-degeneracy. This technical condition is called Strong Coherent Orientation Condition (SCOC), to which we introduce now. For simplicity, we remove the equality constraints defining the constraint set (by removing redundant constraints and using change of variables), and set
\begin{equation*}  
  K(\bp) \defeq \left\{ \bx \in \R^n \ \big| \ g(\bx,\bp) \le \bzero \in \R^I \right\}
\end{equation*}
where for each $i = 1,\dots,I$ and $\bp \in \calP$, $g_i (\cdot, \bp):\R^n \rightarrow \R$ are convex and $\calC^2$-smooth. For each pair $(\bx,\bp)$ with $\bx \in \mbox{SOL}(K(\bp),F(\cdot, \bp))$, and for each multiplier $\lambda \in \R^I$, denote the Lagrangian function 
\begin{equation*}
  L(\bx,\lambda,\bp) \defeq F(\bx,\bp) + \sum_{i=1}^{I} \lambda_i \nabla_{\bx} g_i(\bx,\bp)
\end{equation*} 
Let $\calM(\bx^*,\bp^*)$ denote the set of multipliers $\lambda \in \R^I$ satisfying the KKT system \eqref{eq:KKT}. A basic result states that $\calM(\bx^*,\bp^*)$ is non-empty if $\bx^* \in \mbox{SOL}(K(\bp^*),F(\cdot, \bp^*))$ and MFCQ holds at $\bx^*$ \cite[Proposition 3.2.1]{facchinei2003finite}. 

We introduce a partition of the index set $\left\{1,\dots,I \right\}$ with respect to the triple $(\bx^*,\lambda,\bp^*)$ such that $\lambda \in \calM(\bx^*,\bp^*)$:
\begin{equation*}
  \mbox{support}(\lambda) \defeq \{ i \mid \lambda_i > 0 = g_i(\bx^*,\bp^*) \}
\end{equation*}
\begin{equation*}
  \mbox{degenerate}(\lambda) \defeq \{ i \mid \lambda_i = 0 = g_i(\bx^*,\bp^*) \}
\end{equation*}
\begin{equation*}
  \mbox{inactive}(\lambda) \defeq \{ i \mid \lambda_i = 0 > g_i(\bx^*,\bp^*) \}
\end{equation*}
In addition, denote by $\mbox{active}(\lambda) \defeq \mbox{support}(\lambda) \cup \mbox{degenerate}(\lambda)$. We remark that by our second characterization of non-degeneracy, $\bx^* \in \mbox{SOL}(K(\bp^*),F(\cdot, \bp^*))$ is non-degenerate if and only if $\mbox{degenerate}(\lambda) = \emptyset$. Let $\calM^e(\bx^*,\bp^*)$ ("e" for \emph{extremal}) be the subset of $\calM(\bx^*,\bp^*)$ consisting of $\lambda \in \calM(\bx^*,\bp^*)$ for which $\{ \nabla_{\bx} g_i(\bx^*,\bp^*) \mid i \in \mbox{active}(\lambda) \}$ are linearly independent.

For each $d \bp \in \R^p$ and each $\lambda \in \calM(\bx^*,\bp^*)$, define the \emph{directional critical set}
\begin{align*}
  \calC (\bx^*, \lambda; \bp^*, d \bp) & \defeq \big\{ \bv \in \R^n: \\
   & \nabla_{\bx} g_i(\bx^*, \bp^*) \cdot \bv^\top + \nabla_{\bp} g_i(\bx^*,\bp^*) \cdot d \bp = 0 \ \ \forall i \in \mbox{support}(\lambda) \\
   & \nabla_{\bx} g_i(\bx^*, \bp^*) \cdot \bv^\top + \nabla_{\bp} g_i(\bx^*,\bp^*) \cdot d \bp \le 0 \ \ \forall i \in \mbox{degenerate}(\lambda)  \big\}
\end{align*}
Proposition 5.4.7 of \cite{facchinei2003finite} states that for each $\lambda' \in \calM(\bx^*,\bp^*)$, the set $\calC (\bx^*, \lambda'; \bp^*, d \bp)$ is non-empty if and only if $\lambda'$ solves the linear program
\begin{equation}  \label{eq:lp-1}
  \mbox{maximize } \ \sum_{i=1}^{I} \lambda_i \nabla_{\bp} g_i(\bx^*,\bp^*) d \bp^\top \mbox{ subject to } \lambda \in \calM(\bx^*,\bp^*)
\end{equation}
The dual of \eqref{eq:lp-1} is the linear program 
\begin{equation}\label{eq:lp-2}
  \mbox{minimize } \ F(\bx^*,\bp^*) \cdot d \bx^\top \mbox{ subject to } \nabla_{\bx} g_i(\bx^*,\bp^*) d \bx^\top + \nabla_{\bp} g_i(\bx^*, \bp^*) \cdot d \bp^\top 
\end{equation}
Let $\calM^c(\bx^*,\bp^*;d\bp)$ and $\calD^c(\bx^*,\bp^*;d\bp)$ be respectively the set of solutions to \eqref{eq:lp-1} and \eqref{eq:lp-2} ("c" for critical). We have 
\begin{equation*}
  \calC (\bx^*, \lambda; \bp^*, d \bp) = \begin{cases}
                                           \calD^c(\bx^*,\bp^*;d\bp), & \mbox{if } \lambda \in \calM^c(\bx^*,\bp^*;d\bp) \\
                                           \emptyset, & \mbox{if } \lambda \in \calM(\bx^*,\bp^*;d\bp) \backslash \calM^c(\bx^*,\bp^*;d\bp)
                                         \end{cases}
\end{equation*}
so it is either empty or a polyhedron set defined by $\calD^c(\bx^*,\bp^*;d\bp)$. We use the abbreviation $\nabla_{\bx} g_{\calI}$ for the set of vectors $\{\nabla_{\bx} g_{i} \mid i \in \calI \}$. Define the \emph{SCOC family of index sets} as 
\begin{align*}
  \calB(\bx^*,\bp^*) & \defeq \big\{ \calI \subset \{1,\dots,I\} \ \big| \ \exists \lambda \in \calM(\bx^*,\bp^*) \ \mbox{such that } \\ 
  & \qquad \mbox{support}(\lambda) \subset \calI \subset \mbox{active}(\lambda) \mbox{ and } \nabla_{\bx} g_{\calI}(\bx^*,\bp^*) \mbox{ are linearly independent} \big\}
\end{align*}
with "$\calB$" for basis. $\calB(\bx^*,\bp^*)$ is non-empty and finite: non-empty since $\calM^e(\bx,\bp) \subset \calB(\bx^*,\bp^*)$, finite since there are only finitely many constraints. Let $B$ be the cardinal of $\calB(\bx^*,\bp^*)$ and enumerate the elements as $\calB(\bx^*,\bp^*) = \{\calI^1,\dots,\calI^B \}$. Each $\calI^b$ gives rise to a unique multiplier $\lambda^b \in \calM^e(\bx,\bp)$; while an element of $\calM^e(\bx,\bp)$ may correspond to multiple index sets in $\calB(\bx^*,\bp^*)$. Define the matrix 
\begin{equation*}
 \Lambda^b \defeq \begin{pmatrix}
  J_{\bx} L(\bx^*,\lambda^b,\bp^*) & \nabla_{\bx} g_{\calI^b} (\bx^*,\bp^*)^\top \\
  \nabla_{\bx} g_{\calI^b} (\bx^*,\bp^*) & \bzero 
\end{pmatrix}
\end{equation*}
where "$J_{\bx}$" means the Jacobian with derivation on $\bx$. By definition, the SCOC holds at $\bx^* \in K(\bp^*)$ if all the $B$ matrices $\Lambda^b$ for $b=1,\dots,B$ have the same non-zero determinant sign. Again, note that if $\bx^* \in \mbox{SOL}(K(\bp^*),F(\cdot, \bp^*))$ is non-degenerate, then $B=1$ and SCOC is trivial. 

The following result taken from \cite[Theorem 5.4.12]{facchinei2003finite} provides sufficient conditions we desired for satisfying Assumptions \ref{assumption:lipschitz-continuity} and \ref{assumption:lipschitz-continuity-projection}. The proof is technical and we refer to the aforementioned reference for details. Note that the conclusion is stronger than Lipschitz continuity, as it shows that the solution mapping $\bx(\bp)$ is actually piecewise-$\calC^1$. In addition, it also provides explicit formula to compute the directional derivatives. 

\begin{theorem} [Sufficient Conditions for Satisfying Assumptions \ref{assumption:lipschitz-continuity} and \ref{assumption:lipschitz-continuity-projection}] \hypertarget{thm:sufficient-conditions-lipschitz}{}
  Let $\bx^* \in \emph{SOL}(K(\bp^*),F(\cdot, \bp^*))$ satisfies the MFCQ,CRCQ and SCOC, where $K(\cdot)$ is given by \eqref{eq:K-definition} and $F(\cdot)$ is $\calC^2$-smooth. Then 
  \begin{enumerate}[label=\alph*)]
    \item There exists a neighbourhood $\calX$ of $\bx^*$, a neighbourhood $\calW$ of $\bp^*$ and a piecewise-$\calC^1$ function $\bx: \calW \rightarrow \calX$ such that $\bx(\bp^*) = \bx^*$. And for every $\bp \in \calW$, $\bx(\bp)$ is the only solution in \emph{SOL}$(K(\bp),F(\cdot,\bp))$ 
    \item If $F(\bx) = \bq + \bx \cdot \bM $ for some matrix $\bM \in \R^{n \times n}$ and some vector $\bq \in \R^n$, then we write VI($K,F$) as AVI$(K,\bq,\bM)$ ("A" for \emph{affine}). For all $d \bp \in \R^p$ and $\lambda \in \calM^e (\bx^*,\bp^*) \cap \calM^c (\bx^*,\bp^*,d\bp)$, the following
    \begin{equation*}
      \emph{AVI}\left( \calD^c(\bx^*,\bp^*;d\bp), \ d\bp \cdot J_{\bp} L(\bx^*,\lambda,\bp^*), \ J_{\bx} L(\bx^*,\lambda,\bp^*) \right)
    \end{equation*}
    has a unique solution and is equal to the directional derivative $\bx^{'}(\bx^*;d\bp)$  
    \item If in addition LICQ holds, then $\calM(\bx^*,\bp^*)$ reduced to a singleton $\{\lambda^*\}$. Moreover, for each $\bp \in \calW$, the pair $(\bx(\bp),\lambda(\bp))$ is the unique solution in $\calX \times \R^p$ for the KKT system \eqref{eq:KKT}, and the function $\bp \mapsto (\bx(\bp),\lambda(\bp))$ is piecewise-$\calC^1$  
  \end{enumerate}
\end{theorem}

If we specific to the case where the feasible regions are polyhedral (i.e. all the functions $g_{t,i}$ are affine), some refinements and simplifications of the above theorem are possible, as already CRCQ holds trivially under this situation. By exploring further this additional affinity, in \cite[Theorem 5.2]{lu2008variational}, it is shown that under a similar determinantal condition as SCOC alone is enough to establish the piecewise-$\calC^1$ property, so that we do not need MFCQ. Furthermore, in \cite[Theorem 4.2]{lu2008variational}, it is shown that this determinantal condition is also a necessary condition for the solution mapping being locally single-valued and Lipschitz-continuous. In addition, in \cite[Equation (6)]{robinson2003constraint} a collection of equivalent conditions are given for the local Lipschitz-continuity under the polyhedral convex feasible set situation.  

We mention that in the end of \cite{robinson1982generalized}, an example of projecting the origin onto a convex polytope is shown to be continuous but not Lipschitz under small perturbations in the constraints defining the polytope. Despite in the linear program case where the Lipschitz property can often be guaranteed via the Hoffman's error bound \cite[Lemma 3.2.3]{facchinei2003finite}, verification of this property is in general non-trivial and involves deep results, see for instance \cite{ioffe2017variational}.

\subsection{More on the Euclidean Projector and the Non-Degeneracy}  \hypertarget{subsec:euclidean-projector}{}

Let $K(\bp)$ be defined as in \eqref{eq:K-definition}. Recall that the Euclidean projection of $\bx \in \R^n$ onto $K(\bp)$, denoted as $\Pi_{K(\bp)}(\bx)$, is the unique solution $\by$ to the convex program 
\begin{equation*}
  \mbox{minimize } \frac{1}{2} (\by-\bx) \cdot (\by-\bx)^\top \mbox{ subject to } \by \in K(\bp)
\end{equation*}
This in turn can be written as VI($K(\bp), \bI - \bx$), where $\bI$ is the identity map. Hence the parametric analysis on the Euclidean projector is a special case of Theorems \hyperlink{thm:sufficient-conditions-C2}{EC.1} and \hyperlink{thm:sufficient-conditions-lipschitz}{EC.2}. This allows us to obtain simplifications on the assumptions. In particular, the technical assumption of SCOC can be dropped.

\begin{theorem} [Sufficient Conditions related to the Euclidean Projector]   \hypertarget{thm:sufficient-conditions-projector}{}
  Let $K(\bp)$ be given by \eqref{eq:K-definition} with $g$ being $\calC^2$-smooth in a neighbourhood of the pair $(\bx_{\pi}^*,\bp^*)$, where $\bp^* \in \calP$ and $\bx_{\pi}^* \defeq \Pi_{K(\bp^*)}(\bx^*)$. If the MFCQ and CRCQ hold at $\bx_{\pi}^* \in K(\bp^*)$, then the function
  \begin{equation} \label{eq:projection-function}
    (\bx,\bp) \mapsto \Pi_{K(\bp)}(\bx)
  \end{equation}
  is piecewise-$\calC^1$ near $(\bx^*,\bp^*)$. Moreover, if $\bx_{\pi}^* \in \emph{VI}(K(\bp^*), \bI - \bx^*)$ is non-degenerate, then \eqref{eq:projection-function} is locally $\calC^2$-smooth. 
\end{theorem}

\proof{Proof}
  The piecewise-$\calC^1$ part is Theorem 4.7.5 of \cite{facchinei2003finite}. The $\calC^2$-smooth part is obtained by combing Corollary 4.1.2 of the same reference.   \Halmos
\endproof

It is instructive to visualize the effect of non-degeneracy via the Euclidean projection. For this purpose, let us suppose that the functions $g(\cdot) $ defining the feasible region are all affine, so that $K(\bp)$ are polyhedrons for all $\bp \in \calP$.

\begin{figure}[htbp]
  \centering
  \input{figure/degeneracy2.tex}
  \caption{Illustration of the non-degeneracy: On the left, $\bx_{\pi}^*$ is non-degenerate. Under small perturbation of $\bp^*$, the projection $\bx_{\pi} = \Pi_{K(\bp)}(\bx^*)$ remains on the same "sticky face". On the right, $\bx_{\pi}^*$ is degenerate. After a small perturbation, $\bx_{\pi}$ "jumps" to a different face, in this case, from a vertex (a $0$-dimensional face) $\bx_{\pi}^*$ to an edge (a $1$-dimensional face) of the polygon on which $\bx_{\pi}$ belongs to. The green lines are used to illustrate the normal manifold induced by $K(\bp^*)$.} \label{fig:degeneracy2}
\end{figure}

Define a collection of index sets
\begin{equation*}
  \mathfrak{F}(\bp) \defeq \left\{ \calI \subset \{1,\dots,I \} \ \big| \ \exists \bx \in \R^n \mbox{ such that } g_i(\bx,\bp) = 0 \ \forall i \in \calI \mbox{ and }  g_i(\bx,\bp) < 0 \ \forall i \notin \calI \right\}
\end{equation*}
Each element $\calI$ in $\mathfrak{F}(\bp)$ is in one-one correspondence to a non-empty face $\mathfrak{F}_{\calI}(\bp)$ of $K(\bp)$: 
\begin{equation*}
  \mathfrak{F}_{\calI}(\bp) \defeq \left\{ \bx \in K(\bp) \ \big| \ g_i(\bx,\bp) = 0 \ \forall i \in \calI \right\}
\end{equation*}
The relative interior of a face $\mathfrak{F}_{\calI}(\bp)$ is given by
\begin{equation*}
  \mathrm{ri}\mathfrak{F}_{\calI}(\bp) \defeq \left\{ \bx \in K(\bp) \ \big| \ g_i(\bx,\bp) = 0 \ \forall i \in \calI \mbox{ and }  g_i(\bx,\bp) < 0 \ \forall i \notin \calI \right\}
\end{equation*}
Then, $\bx_{\pi}^* = \Pi_{K(\bp^*)}(\bx^*)$ is non-degenerate if and only if there exists $\calI^* \in \mathfrak{F}(\bp)$ such that $\bx_{\pi}^* \in \mathrm{ri}\mathfrak{F}_{\calI^*}(\bp^*)$. Moreover, upon small perturbation of $\bp^*$, the projection $\bx_{\pi} = \Pi_{K(\bp)}(\bx^*)$ is confined to the face defined by the \emph{same} index set $\calI^*$. For this reason $\mathfrak{F}_{\calI^*}(\bp^*)$ is called a "sticky-face" to $\bx^*$. This observation has important implications in computation: If $\bx_{\pi}^* = \Pi_{K(\bp^*)}(\bx^*)$ happens to be non-degenerate, then to compute $\bx_{\pi} = \Pi_{K(\bp)}(\bx^*)$ for $\bp$ near $\bp^*$, instead of projecting onto $K(\bp)$, we only need to project onto $\mathfrak{F}_{\calI^*}(\bp)$, which is itself a polyhedron but with simpler structure. In other words, under non-degeneracy
\begin{equation*}
  \Pi_{K(\bp)}(\bx^*) = \Pi_{\mathfrak{F}_{\calI^*}(\bp)}(\bx^*) \mbox{ for $\bp$ near $\bp^*$}
\end{equation*}

\section{The Policy Classes $\slip$, $\slisse$}  \hypertarget{appendix:policy-classes-extension}{}

In this appendix, we build upon the findings in the main text, extending the discussion on the update and projection policies to a wider spectrum of policy classes. We then investigate the implications associated with policies belonging to these classes. Furthermore, we delve into the intricacies of the established optimality gap bounds, focusing on the exponential growth of a multiplicative constant with respect to the stage number. From a computational complexity perspective, we argue that such exponential growth is generally unavoidable in a multi-stage optimization problem (\cite{dyer2006computational, shapiro2005complexity, reaiche2016note}).

\subsection{Extension to $\slip$ and $\slisse$}

Upon examining the proofs of Theorems \ref{thm:general-theory-of-projection-Lipschitz} and \ref{thm:general-theory-of-projection-C2} pertaining to the projection policy, it becomes apparent that aside from Lipschitz-continuity and $\calC^2$-smoothness, we did not rely on any additional attributes of the Euclidean projection mapping to reach the conclusions of these theorems. This motivates a more general formulation that encapsulates only the essential conditions needed to support the same claim.

To set the stage, we fix an optimal solution $\bu^*[1,T]$ and the corresponding $\bx^*[1,T]$ by solving \eqref{eq:original-problem-rel} for $V_{\mathrm{rel}-} (\bx,T)$. We also fix a positive sequence $\varepsilon_t$ for $1 \le t \le T$ in line with Assumptions \ref{assumption:lipschitz-continuity-projection} or \ref{assumption:C2-smoothness-projection}, depending on the context. Recall the definition of the feasible region at time-step $t$, represented as $\calU_t(\bx,\bw)$ in \eqref{eq:feasible-region-time-t}. By a \emph{feasible policy (mapping)} $S[1 , T]$, we imply T single-valued functions $S_t(\cdot)$ that map $(\bx,\bw)$ to a control $\bu(t) \in \calU_t(\bx,\bw)$, for $1 \le t \le T$. Following that, we define the following two policy classes
\begin{equation*}
  \slip \defeq \left\{ S[1,T]  \ \big| \ S_t(\bx,\bw) \in \calU_t(\bx,\bw) \mbox{ and $S_t(\cdot)$ is Lipschitz-continuous in $\calB((\bx^*(t),\barw), \varepsilon_t) \mbox{ for } 1 \le t \le T $} \right\}
\end{equation*}
\begin{equation*}
  \slisse \defeq \left\{ S[1,T]  \ \big| \ S_t(\bx,\bw) \in \calU_t(\bx,\bw) \mbox{ and $S_t(\cdot)$ is $\calC^2$-smooth in $\calB((\bx^*(t),\barw), \varepsilon_t) \mbox{ for } 1 \le t \le T $} \right\}
\end{equation*}
In the above definitions, we only insist on the single-valuedness and Lipschitz-continuity (resp. $\calC^2$-smoothness) of $S_t(\cdot)$ locally in $\calB((\bx^*(t),\barw), \varepsilon_t)$. Outside these neighborhoods, we only demand the feasibility $S_t(\cdot) \subseteq \calU_t(\cdot)$. This approach stems from our concern for policy mappings to be well-behaved in the $\varepsilon_t$-neighborhood of $(\bx^*(t),\barw)$.

For each $1 \le t \le T$, recall that $\bu^*_{\bx,T+1-t,\bw}(t)$ is the first control from the control sequence $\bu^*_{\bx,T+1-t,\bw}[t,T]$ that stands as an optimal solution of \eqref{eq:original-problem-res} for $\hat{V}_{\mathrm{rel}+} (\bx,T+1-t,\bw)$. For $(\bx(t),\bw) \in \calB((\bx^*(t),\barw), \varepsilon_t)$, we define
\begin{equation}  \label{eq:policy-mapping-one-point-1}
  \slip^{\bx(t),\bw} \defeq \left\{ S[1,T] \in \slip \ \big| \  S_t(\bx(t),\bw) = \bu^*_{\bx(t), T+1-t, \bw}(t)  \right\}
\end{equation}
\begin{equation}   \label{eq:policy-mapping-one-point-2}
  \slisse^{\bx(t),\bw} \defeq \left\{ S[1,T] \in \slisse \ \big| \ S_t(\bx(t),\bw) = \bu^*_{\bx(t), T+1-t, \bw}(t) \right\}
\end{equation}
These are subsets of $\slip$ and $\slisse$ that match on a specific point $(\bx(t), \bw)$ for $S_t(\cdot)$ with the policy mapping $S_t^*(\cdot)$ defined in \eqref{eq:policy-mapping}. With our primary focus on the deterministic trajectory in the asymptotic limit when all variances converge to zero, we denote by
\begin{equation} \label{eq:intersection-of-deterministic}
  \slip^* \defeq \bigcap_{t=1}^T \slip^{\bx^*(t),\barw} \mbox{ and } \slisse^* \defeq \bigcap_{t=1}^T \slisse^{\bx^*(t),\barw}
\end{equation}
As per our earlier analysis, the projection policy belongs to $\slip^*$ (resp. $\slisse^*$) under Assumption \ref{assumption:lipschitz-continuity-projection} (resp. Assumption \ref{assumption:C2-smoothness-projection}). The subsequent result emerges directly from the proofs of Theorems \ref{thm:general-theory-of-projection-Lipschitz} and \ref{thm:general-theory-of-projection-C2}.

\begin{corollary} [Optimality Gap Bounds for Policy Mappings in $\slip^*$ and $\slisse^*$]   \hypertarget{cor:optimality-gap-bounds-general}{}

Recall the setup in Theorems \ref{thm:general-theory-of-projection-Lipschitz} and \ref{thm:general-theory-of-projection-C2}. For any policy mapping $S[1,T] \in \slip^*$, we have
  \begin{equation}  \label{eq:general-conclusion-1}
    V_{\mathrm{opt}} (\bx,T) - V_{S[1,T]} (\bx,T) \le \mathfrak{P}' C'_3 +  (1 - \mathfrak{P}') \bar{V}
  \end{equation}
 And for any policy mapping $S[1,T] \in \slisse^*$, we have
  \begin{equation}  \label{eq:general-conclusion-2}
    V_{\mathrm{opt}} (\bx,T) - V_{S[1,T]} (\bx,T) \le \mathfrak{P}' C'_4 +  (1 - \mathfrak{P}') \bar{V}
  \end{equation}
  where $C'_3$ and $C'_4$ are given respectively in \eqref{eq:choice-of-constant-projection-Lipschitz} and \eqref{eq:choice-of-constant-projection-C2}, with the constants related to \emph{proj} replaced by the same quantities from policy mapping $S[1,T]$.
\end{corollary}

\begin{remark}[Asymptotic Optimality and Robust Optimality]  \hypertarget{rem:asymptotic-vs-robustness}{} 

The motivation to consider policy mappings in $\slip^*$ and $\slisse^*$, as opposed to the broader classes $\slip$ and $\slisse$, stems from the pursuit of asymptotic optimality. If the variances $\var{\bW}$ and $\var{\calE}$ can be manipulated and reduced, then as variances approach zero, the performance of the policies in these classes tends towards the optimal value. This has been illustrated in Section \ref{subsec:interpretation-of-the-optimality-gap-bounds}. However, under circumstances where variances cannot be controlled by the decision-maker, adopting a robustness perspective may be more appropriate. Specifically, by selecting a suitable constant $\delta_t > 0$ for each time-step $1 \le t \le T$, and finding a policy mapping in the larger class $\slip$ and satisfies
  \begin{equation}  \label{eq:reboust-formulation}
    \inf_{(\bx,\bw) \in \calB((\bx^*(t),\barw), \varepsilon_t)} R_t \left( \bx, \bw, S_t(\bx, \bw) \right) \ge R_t \left( \bx^*(t), \barw, S^*_t(\bx^*(t),\barw) \right) - \delta_t \ \ \mbox{ for all } 1 \le t \le T
  \end{equation}
it will lead to the same conclusion of \eqref{eq:general-conclusion-1} for this policy mapping $S[1,T]$, assuming that the sum $\sum_{t=1}^{T} \delta_t$ does not exceed the constant $C'_3$, and by treating the situation outside the $\varepsilon$-neighbourhood in the same manner. 

It is worth noting that solving \eqref{eq:reboust-formulation} essentially amounts to controlling the worst-case performance within an $\varepsilon$-neighborhood of a deterministically optimal trajectory (the nominal one), which is a typical objective in a robust optimization problem. While stochastic optimization operates under the presumption that uncertainty has a probabilistic description, robust optimization posits a deterministic, set-based model for uncertainty (\cite{bertsimas2011theory}). This approach aims to devise a solution that remains feasible for \emph{any} realization of uncertainty within a specified set. Interestingly, our CEC-based heuristics already exhibit a distributional robustness, as they focus solely on the first and second moments of these distributions. Therefore, as indicated in Remark \hyperlink{rem:asymptotic-vs-robustness}{EC.1}, by choosing the uncertainty set as an $\varepsilon$-neighbourhood of the deterministically optimal trajectory, solving the robust counterpart \eqref{eq:reboust-formulation} results in the same optimality gap bounds. Robust finite-horizon Markov decision processes with finite state and action spaces have been explored in \cite[Chapter 13]{ben2009robust}, and a robust multi-stage optimization framework is discussed in \cite[Chapter 14]{ben2009robust}. Nevertheless, solving \eqref{eq:reboust-formulation} remains challenging due to its infinite dimensionality. A potential solution could be considering affine solution mappings, as proposed in \cite[Chapter 14]{ben2009robust}.     \Halmos
  
\end{remark}

The above corollary essentially posits that for the establishment of a first-order or second-order optimality gap bound, one needs only to extrapolate a \emph{feasible} policy mapping as a Lipschitz-continuous, or a $\calC^2$-smooth function around a deterministically optimal trajectory. The major complexities stem from the feasibility condition, imposed by the stringent constraints. The update and projection policies proffer two universal methods for such a construction. However, an unanswered query remains: do they represent a sound choice within the corresponding policy class in terms of performance?

\subsection{The Lipschitz Constant of $\hat{V}_{\mathrm{rel}+} (\cdot, T-t, \cdot)$}  \hypertarget{subsec:Lipschitz-constant-of-V-rel-post}{}

To investigate this question, consider the non-stochastic part of the constant $C_1$ from \eqref{eq:choice-of-constant} in the update policy's optimality gap bound, it is governed by the sum of $T$ Lipschitz constants $K_t$ and $L_t$ of the optimal value functions $\hat{V}_{\mathrm{rel}+} (\cdot, T-t, \cdot)$. While in the corresponding constant $C_3$ from \eqref{eq:choice-of-constant-projection-Lipschitz} for the projection policy, it is governed by $c_{\phi}^{T} (\cproj + 1)^{T}$, a term that may explode exponentially with $T$. We need to understand the growth rate of the constants $K_t$ with $t$.

First observe that under Assumption \ref{assumption:lipschitz-continuity}, the policy mapping $S^*[1,T]$ obtained from $S_t^*(\cdot)$ in \eqref{eq:policy-mapping} with $1 \le t \le T$ is such that 
\begin{equation} \label{eq:intersection-of-policy-class}
  S^*[1,T] \in \bigcap_{t=1}^T \bigcap_{\substack{(\bx(t),\bw) \in \\ \calB((\bx^*(t),\barw), \varepsilon_t)}} \slip^{\bx(t),\bw}
\end{equation}
In particular, the intersection on the right-hand-side of \eqref{eq:intersection-of-policy-class} is non-empty. The following result then characterizes an upper bound to the Lipschitz constants $K_t$ of the mappings $(\bx,\bw) \mapsto \hat{V}_{\mathrm{rel}+} (\bx, T-t, \bw)$ for $1 \le t \le T$.

\begin{corollary} [Explicit Bound with Lipschitz-Continuity in Update Policy]     \hypertarget{cor:optimality-gap-bound-more-explicit}{}
  Under Assumption \ref{assumption:lipschitz-continuity}, denote by $c_{\mathrm{lip}, S'_t(\cdot)}^{\bx(t), \bw}$ the Lipschitz constant of \emph{any} policy mapping $S'_t(\cdot) \in \slip^{\bx(t),\bw}$, where by Assumption \ref{assumption:lipschitz-continuity} $\slip^{\bx(t),\bw}$ is non-empty. Define 
\begin{equation}  \label{eq:choice-of-constant-c-lip}
     c_{\mathrm{lip},t} \defeq \sup_{\substack{(\bx(t),\bw) \in \\ \calB((\bx^*(t),\barw), \varepsilon_t)}} \inf_{S'_t(\cdot) \in \slip^{\bx(t),\bw}} c_{\mathrm{lip}, S'_t(\cdot)}^{\bx(t), \bw}
\end{equation}
  for $1 \le t \le T$. Then an upper bound $\overline{K}_t$ for the Lipschitz constant $K_t$ of the value function $\hat{V}_{\mathrm{rel}+} (\cdot, T-t, \cdot): \R^{n_x} \times \R^{n_w} \rightarrow \R$ can be constructed inductively backward on $t$ via $\overline{K}_{T+1} = 0$, and 
  \begin{equation}  \label{eq:choice-of-K}
    \overline{K}_{t-1} = \overline{K}_t c_{\phi} \cdot (c_{\mathrm{lip},t} + 1) + (c_R+1) \cdot c_{\mathrm{lip},t},
  \end{equation}
  where $c_{\phi}$ is the Lipschitz constant of the deterministic state transition function $\phi(\cdot)$, and $c_R \defeq \max_t c_{R_t}$ with $c_{R_t}$ being the Lipschitz constant of the reward function $R_t(\cdot)$.   
   
\end{corollary}

\proof{Proof of Corollary \hyperlink{cor:optimality-gap-bound-more-explicit}{EC.2}}

Recall that we denote by $c_{\phi}$ the Lipschitz constant of the deterministic state transition function $\phi(\cdot)$, and $c_R \defeq \max_t c_{R_t}$ with $c_{R_t}$ being the Lipschitz constant of the reward function $R_t(\cdot)$.

In the following, we construct an upper bound $\overline{K}_t$ for the Lipschitz constant $K_t$ of the value function $\hat{V}_{\mathrm{rel}+} (\cdot, T+1-t, \cdot): \R^{n_x} \times \R^{n_w} \rightarrow \R$ using induction backward on $t$, starting with $\overline{K}_{T+1}=0$. Suppose at time-step $t$ we have an upper bound $\overline{K}_t$ for the Lipschitz constant $K_t$ such that 
\begin{align}
  & \mbox{For any } (\bx(t+1),\bw) \in \calB \left( (\bx^*(t+1),\barw), \varepsilon_{t+1} \right) \mbox{ we have:} \nonumber \\
  & \qquad \abs{\hat{V}_{\mathrm{rel}+} (\bx(t+1),T-t,\bw) - \hat{V}_{\mathrm{rel}+} (\bx^*(t+1),T-t,\barw)} \le \overline{K}_t \norme{(\bx(t+1),\bw) - (\bx^*(t+1),\barw)} \label{eq:induction-on-t}
\end{align}
We proceed to construct an upper bound $\overline{K}_{t-1}$ for $K_{t-1}$. 

Choose $S_t(\cdot) \in \slip^{\bx^*(t),\barw}$ and denote by  $c_{\mathrm{lip}, S_t(\cdot)}^{\bx^*(t), \barw}$ the Lipschitz constant of the policy mapping $S_t(\cdot)$ inside $\calB((\bx^*(t),\barw), \varepsilon_t)$. Fix any $(\bx(t), \bw) \in \calB \left( (\bx^*(t),\barw), \varepsilon_{t} \right)$. From \eqref{eq:principle-of-optimality} and by the optimality of the control $\bu^*_{\bx(t),T+1-t,\bw}(t)$, we deduce that 
\begin{align*}
  \hat{V}_{\mathrm{rel}+} (\bx(t),T+1-t,\bw) = & \ R_{t} \big( \bx(t), \bw, \bu^*_{\bx(t),T+1-t,\bw}(t) \big) +  V_{\mathrm{rel}-} (\phi \big(\bx(t), \bw, \bu^*_{\bx(t),T+1-t,\bw}(t) \big),T-t)  \\
  \ge & \ R_{t} \big( \bx(t), \bw, S_t(\bx(t),\bw) \big) + V_{\mathrm{rel}-} (\phi \big(\bx(t), \bw, S_t(\bx(t),\bw) \big),T-t)
\end{align*}
While by construction $S_t(\bx^*(t),\barw) = \bu^*_{\bx^*(t),T+1-t,\barw}(t)$ is an optimal control, hence
\begin{equation*}
  \hat{V}_{\mathrm{rel}+} (\bx^*(t),T+1-t,\barw) = R_{t} \big( \bx^*(t), \barw, S_t(\bx^*(t),\barw) \big) + V_{\mathrm{rel}-} (\phi \big(\bx^*(t), \barw, S_t(\bx^*(t),\barw) \big),T-t)
\end{equation*}
Consequently by taking the difference, and recall that $V_{\mathrm{rel}-} (\bx,T-t) = \hat{V}_{\mathrm{rel}+} (\bx,T-t,\barw)$, we obtain
\begin{align}
  & \hat{V}_{\mathrm{rel}+} (\bx^*(t),T+1-t,\barw) - \hat{V}_{\mathrm{rel}+} (\bx(t),T+1-t,\bw) \nonumber \\
  & \le \ (c_R +1) c_{\mathrm{lip}, S_t(\cdot)}^{\bx^*(t), \barw} \cdot \norme{(\bx^*(t),\barw) - (\bx(t),\bw)} \nonumber \\
  & \qquad + \hat{V}_{\mathrm{rel}+} (\phi \big(\bx^*(t), \barw, S_t(\bx^*(t),\barw) \big),T-t, \barw) - \hat{V}_{\mathrm{rel}+} (\phi \big(\bx(t), \bw, S_t(\bx(t),\bw) \big),T-t,\barw) \nonumber \\
  & \le \ \left( \overline{K}_t c_{\phi} \cdot (c_{\mathrm{lip}, S_t(\cdot)}^{\bx^*(t), \barw} +1) + (c_R+1) \cdot c_{\mathrm{lip}, S_t(\cdot)}^{\bx^*(t), \barw} \right) \cdot \norme{(\bx^*(t),\barw) - (\bx(t),\bw)} \label{eq:bound-direction-1-fix}
\end{align}
where the last step follows from our induction hypothesis on $\overline{K}_t$ in \eqref{eq:induction-on-t}, and by shrinking $\varepsilon_t$ if necessary, so that $\phi \big(\bx(t), \bw, S_t(\bx(t),\bw) \big)$ remains in the $\varepsilon_{t+1}$-neighbourhood of $\bx^*(t+1)$, for any $(\bx(t), \bw) \in \calB \left( (\bx^*(t),\barw), \varepsilon_{t} \right)$. 

Since \eqref{eq:bound-direction-1-fix} holds for \emph{any} $S_t(\cdot) \in \slip^{\bx^*(t),\barw}$, we deduce that
\begin{align} 
  & \hat{V}_{\mathrm{rel}+} (\bx^*(t),T+1-t,\barw) - \hat{V}_{\mathrm{rel}+} (\bx(t),T+1-t,\bw)   \nonumber \\
  & \le \qquad \inf_{S_t(\cdot) \in \slip^{\bx^*(t),\barw}} \left( \overline{K}_t c_{\phi} \cdot (c_{\mathrm{lip}, S_t(\cdot)}^{\bx^*(t), \barw} +1) + (c_R+1) \cdot c_{\mathrm{lip}, S_t(\cdot)}^{\bx^*(t), \barw} \right) \cdot \norme{(\bx^*(t),\barw) - (\bx(t),\bw)} \label{eq:bound-direction-1}
\end{align}

We now interchange the roles of $(\bx^*(t), \barw)$ and $(\bx(t),\bw)$ above, by choosing a policy mapping $S'_t(\cdot) \in \slip^{\bx(t),\bw}$ with Lipschitz constant $c_{\mathrm{lip}, S'_t(\cdot)}^{\bx(t), \bw}$. This time we obtain  
\begin{equation*}
  \hat{V}_{\mathrm{rel}+} (\bx(t),T+1-t,\bw) = R_{t} \big( \bx(t), \bw, S'_t(\bx(t),\bw) \big) + V_{\mathrm{rel}-} (\phi \big(\bx(t), \bw, S'_t(\bx(t),\bw) \big),T-t)
\end{equation*}
while 
\begin{equation*}
  \hat{V}_{\mathrm{rel}+} (\bx^*(t),T+1-t,\barw) \ge R_{t} \big( \bx^*(t), \barw, S'_t(\bx^*(t),\barw) \big) + V_{\mathrm{rel}-} (\phi \big(\bx^*(t), \barw, S'_t(\bx^*(t),\barw) \big),T-t) 
\end{equation*}
Hence we deduce the following inequality in the other direction
\begin{align} 
  & \hat{V}_{\mathrm{rel}+} (\bx(t),T+1-t,\bw) - \hat{V}_{\mathrm{rel}+} (\bx^*(t),T+1-t,\barw) \nonumber \\
  & \qquad \le \inf_{S'_t(\cdot) \in \slip^{\bx(t),\bw}}  \left( \overline{K}_t c_{\phi} \cdot (c_{\mathrm{lip}, S'_t(\cdot)}^{\bx(t), \bw} +1) + (c_R+1) \cdot c_{\mathrm{lip}, S'_t(\cdot)}^{\bx(t), \bw} \right) \cdot \norme{(\bx^*(t),\barw) - (\bx(t),\bw)} \label{eq:bound-direction-2}
\end{align}
Combining \eqref{eq:bound-direction-1} and \eqref{eq:bound-direction-2}, we can take the absolute value and deduce that
\begin{align}
  & \abs{\hat{V}_{\mathrm{rel}+} (\bx^*(t),T+1-t,\barw) - \hat{V}_{\mathrm{rel}+} (\bx(t),T+1-t,\bw)} \le \norme{(\bx^*(t),\barw) - (\bx(t),\bw)} \cdot  \nonumber \\
  & \max \Bigg\{ \inf_{S_t(\cdot) \in \slip^{\bx^*(t),\barw}} \left( \overline{K}_t c_{\phi} \cdot (c_{\mathrm{lip}, S_t(\cdot)}^{\bx^*(t), \barw} +1) + (c_R+1) \cdot c_{\mathrm{lip}, S_t(\cdot)}^{\bx^*(t), \barw} \right),  \nonumber \\
  & \inf_{S'_t(\cdot) \in \slip^{\bx(t),\bw}}  \left( \overline{K}_t c_{\phi} \cdot (c_{\mathrm{lip}, S'_t(\cdot)}^{\bx(t), \bw} +1) + (c_R+1) \cdot c_{\mathrm{lip}, S'_t(\cdot)}^{\bx(t), \bw} \right)  \Bigg\}  \label{eq:bound-both-direction}
\end{align}
Define $c_{\mathrm{lip},t}$ as claimed in \eqref{eq:choice-of-constant-c-lip} of the corollary. From \eqref{eq:bound-both-direction} we deduce that for all $(\bx(t),\bw) \in \calB((\bx^*(t),\barw), \varepsilon_t)$, we have
\begin{align}
  & \abs{\hat{V}_{\mathrm{rel}+} (\bx^*(t),T+1-t,\barw) - \hat{V}_{\mathrm{rel}+} (\bx(t),T+1-t,\bw)} \nonumber \\
  & \le \left( \overline{K}_t c_{\phi} \cdot (c_{\mathrm{lip},t} + 1) + (c_R+1) \cdot c_{\mathrm{lip},t} \right) \cdot \norme{(\bx^*(t),\barw) - (\bx(t),\bw)} \label{eq:bound-of-lip}
\end{align}
So that we can choose $\overline{K}_{t-1} \defeq \overline{K}_t c_{\phi} \cdot (c_{\mathrm{lip},t} + 1) + (c_R+1) \cdot c_{\mathrm{lip},t}$ to complete the induction step.   \Halmos

\endproof

Denote by $c_{\mathrm{update},t}$ the Lipschitz constant of $S^*_t(\cdot)$, then from \eqref{eq:intersection-of-policy-class} and \eqref{eq:choice-of-constant-c-lip} we have $c_{\mathrm{update},t} \ge c_{\mathrm{lip},t}$. By inspecting the proof of Corollary \hyperlink{cor:optimality-gap-bound-more-explicit}{EC.2}, we see that a looser, and non-explicit upper bound of $K_t$ can be constructed by replacing $c_{\mathrm{lip},t}$ in \eqref{eq:choice-of-K} by $c_{\mathrm{update},t}$. On the other hand, with $c_{\mathrm{proj},t}$ denoting the Lipschitz constant of the projection mapping $\proj(\bx,\bw,t)$, guaranteed from Assumption \ref{assumption:lipschitz-continuity-projection}, we have $c_{\mathrm{proj},t} \ge \inf_{S_t(\cdot) \in \slip^{\bx^*(t),\barw}} c_{\mathrm{lip}, S_t(\cdot)}^{\bx^*(t), \barw}$. While it is in general not possible to directly compare $ c_{\mathrm{lip},t}$ with $c_{\mathrm{proj},t}$ due to the 
"$\sup_{\substack{(\bx(t),\bw) \in \\ \calB((\bx^*(t),\barw), \varepsilon_t)}}$" part in \eqref{eq:choice-of-constant-c-lip}, the point is that as long as there \emph{exists} one Lipschitz extrapolation to a policy mapping $S'_t(\cdot)$ satisfying $S'_t(\bx(t),\bw) = \bu^*_{\bx(t), T+1-t, \bw}(t)$ and having a small Lipschitz constant, uniformly for all $(\bx(t),\bw) \in \calB((\bx^*(t),\barw), \varepsilon_t)$, then $c_{\mathrm{lip},t}$ would be small. This consideration is \emph{independent} of the policy that we use; while the constant $c_{\mathrm{proj},t}$ depends on the specific policy (the projection policy here) taken. 

As an illustration, consider the weakly coupled Markov decision processes (WCMDPs) discussed in Example \ref{example-2}. From \eqref{eq:WCMDP-2}, the function $\phi(\cdot)$ depends merely on the control part $\bU$, and not on $\bX$ nor $\bW$. This implies that \eqref{eq:choice-of-K} can be refined to $\overline{K}_{t-1} = \overline{K}_t c_{\phi} \cdot c_{\mathrm{lip},t} + (c_R+1) \cdot c_{\mathrm{lip},t}$. One can then show that $c_{\phi} \le 1$, and, by exhibiting one Lipschitz extrapolation as in \cite[Proposition 4.1]{brown2022fluid}, that $c_{\mathrm{lip},t} \le 1$ for $1 \le t \le T$. It implies that the constant $C_1$ in \eqref{eq:choice-of-constant} of Theorem \ref{thm:general-theory-of-resolving-Lipschitz} can be chosen to grow at most quadratically with $T$ (and not exponentially). In contrast, it is observed numerically in \cite{yan2022close} that, for a large class of policies that do not involve any re-solving (a so-called "one-pass policy"), there exists some WCMDP on which the constant $c_{\phi} \cdot c_{*,t} $ is strictly larger than $1$, where $c_{*,t}$ is the Lipschitz constant depending on the specific one-pass policy under consideration.   

One may argue that from \eqref{eq:choice-of-K} the constant $K_t$ can still explode exponentially with $t$. We believe that this is unavoidable, unless the problem admits additional structure as in Example \ref{example-2}. This is based on the computational complexity results already mentioned in Remark \ref{rem:model-assumptions}, that a $T$-stage stochastic optimization problem is drastically more difficult than a $2$-stage problem. Additionally, there are evidences from the so-called \emph{sample complexity} on $T$-stage problems with only the i.i.d. noises $\bW$, that supports the exponential growth, see e.g. \cite[Section 3.2]{shapiro2005complexity}. A lower bound in \cite{reaiche2016note} for a certain class of $T$-stage problems even gives a growth rate of the sample complexity with a multiplicative constant $T^T$.

\bibliographystyle{informs2014} 
\bibliography{reference}

% Appendix here
% Options are (1) APPENDIX (with or without general title) or
%             (2) APPENDICES (if it has more than one unrelated sections)
% Outcomment the appropriate case if necessary
%
% \begin{APPENDIX}{<Title of the Appendix>}
% \end{APPENDIX}
%
%   or
%
% \begin{APPENDICES}
% \section{<Title of Section A>}
% \section{<Title of Section B>}
% etc
% \end{APPENDICES}

% References here (outcomment the appropriate case)

% CASE 1: BiBTeX used to constantly update the references
%   (while the paper is being written).

% CASE 2: BiBTeX used to generate mypaper.bbl (to be further fine tuned)
%\input{mypaper.bbl} % outcomment this line in Case 2

%If you don't use BiBTex, you can manually itemize references as shown below.

%%%%%%%%%%%%%%%%%
\end{document}

%% file: figure/network4.tex
\tikzset{every picture/.style={line width=0.75pt}} %set default line width to 0.75pt        

\begin{tikzpicture}[x=0.75pt,y=0.75pt,yscale=-1,xscale=1]
%uncomment if require: \path (0,300); %set diagram left start at 0, and has height of 300

%Flowchart: Sort [id:dp7225100030156695] 
\draw   (417.05,150.99) -- (301.87,234.24) -- (185.36,152.85) -- (300.54,69.6) -- cycle ; \draw   (300.54,69.6) -- (301.87,234.24) ;
%Straight Lines [id:da6981791395270927] 
\draw [color={rgb, 255:red, 208; green, 2; blue, 27 }  ,draw opacity=1 ][line width=1.5]    (183.43,147.79) -- (301.99,62.54) ;
\draw [shift={(304.43,60.79)}, rotate = 144.28] [color={rgb, 255:red, 208; green, 2; blue, 27 }  ,draw opacity=1 ][line width=1.5]    (14.21,-4.28) .. controls (9.04,-1.82) and (4.3,-0.39) .. (0,0) .. controls (4.3,0.39) and (9.04,1.82) .. (14.21,4.28)   ;
%Straight Lines [id:da09003210883438051] 
\draw [color={rgb, 255:red, 208; green, 2; blue, 27 }  ,draw opacity=1 ][line width=1.5]    (298.05,63.99) -- (413.94,142.11) ;
\draw [shift={(416.43,143.79)}, rotate = 213.98] [color={rgb, 255:red, 208; green, 2; blue, 27 }  ,draw opacity=1 ][line width=1.5]    (14.21,-4.28) .. controls (9.04,-1.82) and (4.3,-0.39) .. (0,0) .. controls (4.3,0.39) and (9.04,1.82) .. (14.21,4.28)   ;
%Straight Lines [id:da8988965015378902] 
\draw [color={rgb, 255:red, 126; green, 211; blue, 33 }  ,draw opacity=1 ][line width=1.5]    (190.43,156.79) -- (299.01,77.55) ;
\draw [shift={(301.43,75.79)}, rotate = 143.88] [color={rgb, 255:red, 126; green, 211; blue, 33 }  ,draw opacity=1 ][line width=1.5]    (14.21,-4.28) .. controls (9.04,-1.82) and (4.3,-0.39) .. (0,0) .. controls (4.3,0.39) and (9.04,1.82) .. (14.21,4.28)   ;
%Straight Lines [id:da8281125818687816] 
\draw [color={rgb, 255:red, 126; green, 211; blue, 33 }  ,draw opacity=1 ][line width=1.5]    (304.43,80.79) -- (306.39,223.79) ;
\draw [shift={(306.43,226.79)}, rotate = 269.22] [color={rgb, 255:red, 126; green, 211; blue, 33 }  ,draw opacity=1 ][line width=1.5]    (14.21,-4.28) .. controls (9.04,-1.82) and (4.3,-0.39) .. (0,0) .. controls (4.3,0.39) and (9.04,1.82) .. (14.21,4.28)   ;
%Straight Lines [id:da38565612515775727] 
\draw [color={rgb, 255:red, 126; green, 211; blue, 33 }  ,draw opacity=1 ][line width=1.5]    (306.43,226.79) -- (407.01,152.57) ;
\draw [shift={(409.43,150.79)}, rotate = 143.58] [color={rgb, 255:red, 126; green, 211; blue, 33 }  ,draw opacity=1 ][line width=1.5]    (14.21,-4.28) .. controls (9.04,-1.82) and (4.3,-0.39) .. (0,0) .. controls (4.3,0.39) and (9.04,1.82) .. (14.21,4.28)   ;
%Straight Lines [id:da6667371691228992] 
\draw [color={rgb, 255:red, 74; green, 144; blue, 226 }  ,draw opacity=1 ][line width=1.5]    (184.43,158.79) -- (299.96,239.07) ;
\draw [shift={(302.43,240.79)}, rotate = 214.8] [color={rgb, 255:red, 74; green, 144; blue, 226 }  ,draw opacity=1 ][line width=1.5]    (14.21,-4.28) .. controls (9.04,-1.82) and (4.3,-0.39) .. (0,0) .. controls (4.3,0.39) and (9.04,1.82) .. (14.21,4.28)   ;
%Straight Lines [id:da8330134149212509] 
\draw [color={rgb, 255:red, 74; green, 144; blue, 226 }  ,draw opacity=1 ][line width=1.5]    (302.43,240.79) -- (416,158.55) ;
\draw [shift={(418.43,156.79)}, rotate = 144.09] [color={rgb, 255:red, 74; green, 144; blue, 226 }  ,draw opacity=1 ][line width=1.5]    (14.21,-4.28) .. controls (9.04,-1.82) and (4.3,-0.39) .. (0,0) .. controls (4.3,0.39) and (9.04,1.82) .. (14.21,4.28)   ;
%Straight Lines [id:da8449394481472501] 
\draw [color={rgb, 255:red, 208; green, 2; blue, 27 }  ,draw opacity=1 ][line width=1.5]    (57.43,37.79) -- (107.43,37.79) ;
\draw [shift={(110.43,37.79)}, rotate = 180] [color={rgb, 255:red, 208; green, 2; blue, 27 }  ,draw opacity=1 ][line width=1.5]    (14.21,-4.28) .. controls (9.04,-1.82) and (4.3,-0.39) .. (0,0) .. controls (4.3,0.39) and (9.04,1.82) .. (14.21,4.28)   ;

%Straight Lines [id:da30415552903211895] 
\draw [color={rgb, 255:red, 126; green, 211; blue, 33 }  ,draw opacity=1 ][line width=1.5]    (55.43,66.79) -- (105.43,66.79) ;
\draw [shift={(108.43,66.79)}, rotate = 180] [color={rgb, 255:red, 126; green, 211; blue, 33 }  ,draw opacity=1 ][line width=1.5]    (14.21,-4.28) .. controls (9.04,-1.82) and (4.3,-0.39) .. (0,0) .. controls (4.3,0.39) and (9.04,1.82) .. (14.21,4.28)   ;

%Straight Lines [id:da11157743819959332] 
\draw [color={rgb, 255:red, 74; green, 144; blue, 226 }  ,draw opacity=1 ][line width=1.5]    (58.43,98.79) -- (108.43,98.79) ;
\draw [shift={(111.43,98.79)}, rotate = 180] [color={rgb, 255:red, 74; green, 144; blue, 226 }  ,draw opacity=1 ][line width=1.5]    (14.21,-4.28) .. controls (9.04,-1.82) and (4.3,-0.39) .. (0,0) .. controls (4.3,0.39) and (9.04,1.82) .. (14.21,4.28)   ;

%Straight Lines [id:da7894236624633972] 
\draw    (53.43,147.79) -- (113.43,147.79) ;
\draw [shift={(115.43,147.79)}, rotate = 180] [color={rgb, 255:red, 0; green, 0; blue, 0 }  ][line width=0.75]    (10.93,-3.29) .. controls (6.95,-1.4) and (3.31,-0.3) .. (0,0) .. controls (3.31,0.3) and (6.95,1.4) .. (10.93,3.29)   ;
%Straight Lines [id:da052252246619302944] 
\draw    (54.43,154.79) -- (114.43,154.79) ;
\draw [shift={(116.43,154.79)}, rotate = 180] [color={rgb, 255:red, 0; green, 0; blue, 0 }  ][line width=0.75]    (10.93,-3.29) .. controls (6.95,-1.4) and (3.31,-0.3) .. (0,0) .. controls (3.31,0.3) and (6.95,1.4) .. (10.93,3.29)   ;
%Straight Lines [id:da01737174532506658] 
\draw    (54.43,161.79) -- (114.43,161.79) ;
\draw [shift={(116.43,161.79)}, rotate = 180] [color={rgb, 255:red, 0; green, 0; blue, 0 }  ][line width=0.75]    (10.93,-3.29) .. controls (6.95,-1.4) and (3.31,-0.3) .. (0,0) .. controls (3.31,0.3) and (6.95,1.4) .. (10.93,3.29)   ;

%Bend Arrow [id:dp4168661738184214] 
\draw   (155.32,165.78) -- (162.41,172.42) .. controls (168.28,177.91) and (168.58,187.13) .. (163.09,193) -- (123.67,235.1) -- (125.64,236.94) -- (111.51,242.18) -- (115.8,227.73) -- (117.77,229.57) -- (157.18,187.47) .. controls (159.62,184.86) and (159.49,180.76) .. (156.88,178.32) -- (149.79,171.69) -- cycle ;

% Text Node
\draw (204.15,101.12) node [anchor=north west][inner sep=0.75pt]  [font=\large,rotate=-325.11] [align=left] {\textit{Link 1}};
% Text Node
\draw (219.26,191.7) node [anchor=north west][inner sep=0.75pt]  [font=\normalsize,rotate=-35.23] [align=left] {{\large \textit{Link 2}}};
% Text Node
\draw (298.14,125.36) node [anchor=north west][inner sep=0.75pt]  [rotate=-89.16] [align=left] {{\large \textit{Link 3}}};
% Text Node
\draw (358.07,73.22) node [anchor=north west][inner sep=0.75pt]  [font=\large,rotate=-36.43] [align=left] {\textit{Link 4}};
% Text Node
\draw (349.66,216.71) node [anchor=north west][inner sep=0.75pt]  [rotate=-325.33] [align=left] {{\large \textit{Link 5}}};
% Text Node
\draw (119,142) node [anchor=north west][inner sep=0.75pt]   [align=left] {{\large \textit{\textbf{Source}}}};
% Text Node
\draw (428,143) node [anchor=north west][inner sep=0.75pt]   [align=left] {{\large \textit{\textbf{Destination}}}};
% Text Node
\draw (114,27) node [anchor=north west][inner sep=0.75pt]   [align=left] {\textcolor[rgb]{0.82,0.01,0.11}{{\large \textit{Path 1}}}};
% Text Node
\draw (114,57) node [anchor=north west][inner sep=0.75pt]   [align=left] {{\large \textcolor[rgb]{0.49,0.83,0.13}{\textit{Path 2}}}};
% Text Node
\draw (116,89) node [anchor=north west][inner sep=0.75pt]   [align=left] {{\large \textcolor[rgb]{0.29,0.56,0.89}{\textit{Path 3}}}};
% Text Node
\draw (57,129) node [anchor=north west][inner sep=0.75pt]   [align=left] {Arrivals};
% Text Node
\draw (93.29,226.29) node [anchor=north west][inner sep=0.75pt]  [rotate=-310.58] [align=left] {Rejections};

\end{tikzpicture}

%% file: figure/degeneracy1.tex
\tikzset{every picture/.style={line width=0.75pt}} %set default line width to 0.75pt        

\begin{tikzpicture}[x=0.6pt,y=0.6pt,yscale=-1,xscale=1]
%uncomment if require: \path (0,300); %set diagram left start at 0, and has height of 300

%Curve Lines [id:da9484608892207627] 
\draw    (212.81,207.14) .. controls (201.68,115.45) and (347.12,67.74) .. (397.92,94.58) ;
%Curve Lines [id:da7625294224491299] 
\draw    (212.81,207.14) .. controls (256.66,254.85) and (406.28,303.3) .. (447.33,220.56) ;
%Straight Lines [id:da10654821206509424] 
\draw    (212.81,207.14) -- (148.1,213.85) ;
%Straight Lines [id:da6410309130431535] 
\draw    (212.81,207.14) -- (166.88,254.85) ;
%Straight Lines [id:da19107976084054368] 
\draw  [dash pattern={on 0.84pt off 2.51pt}]  (74,220.67) -- (148.1,213.85) ;
%Straight Lines [id:da7471648577327883] 
\draw  [dash pattern={on 0.84pt off 2.51pt}]  (126.67,295.33) -- (166.88,254.85) ;
%Straight Lines [id:da43963330821163393] 
\draw    (397.92,94.58) -- (447.33,72) ;
%Straight Lines [id:da46649906484861337] 
\draw    (397.92,94.58) -- (426,52) ;
%Straight Lines [id:da9379951969627214] 
\draw  [dash pattern={on 0.84pt off 2.51pt}]  (426,52) -- (455.33,6) ;
%Straight Lines [id:da6538678254011328] 
\draw  [dash pattern={on 0.84pt off 2.51pt}]  (447.33,72) -- (506.67,45.33) ;
%Curve Lines [id:da8666186144846504] 
\draw    (96.62,220.02) .. controls (88.69,245.39) and (118.48,278.05) .. (133.53,284.65) ;
\draw [shift={(135.33,285.33)}, rotate = 197.65] [fill={rgb, 255:red, 0; green, 0; blue, 0 }  ][line width=0.08]  [draw opacity=0] (12,-3) -- (0,0) -- (12,3) -- cycle    ;
\draw [shift={(97.33,218)}, rotate = 111.8] [fill={rgb, 255:red, 0; green, 0; blue, 0 }  ][line width=0.08]  [draw opacity=0] (12,-3) -- (0,0) -- (12,3) -- cycle    ;
%Curve Lines [id:da48450215812432207] 
\draw    (456.59,10.2) .. controls (466.46,15.21) and (489.64,33.45) .. (494.79,47.59) ;
\draw [shift={(495.33,49.33)}, rotate = 255.38] [fill={rgb, 255:red, 0; green, 0; blue, 0 }  ][line width=0.08]  [draw opacity=0] (12,-3) -- (0,0) -- (12,3) -- cycle    ;
\draw [shift={(454.67,9.33)}, rotate = 21.04] [fill={rgb, 255:red, 0; green, 0; blue, 0 }  ][line width=0.08]  [draw opacity=0] (12,-3) -- (0,0) -- (12,3) -- cycle    ;
%Straight Lines [id:da5240203724945007] 
\draw [color={rgb, 255:red, 208; green, 2; blue, 27 }  ,draw opacity=1 ][fill={rgb, 255:red, 208; green, 2; blue, 27 }  ,fill opacity=1 ]   (212.81,207.14) -- (168,230) ;
\draw [shift={(168,230)}, rotate = 332.98] [color={rgb, 255:red, 208; green, 2; blue, 27 }  ,draw opacity=1 ][line width=0.75]    (17.64,-3.29) .. controls (13.66,-1.4) and (10.02,-0.3) .. (6.71,0) .. controls (10.02,0.3) and (13.66,1.4) .. (17.64,3.29)(10.93,-3.29) .. controls (6.95,-1.4) and (3.31,-0.3) .. (0,0) .. controls (3.31,0.3) and (6.95,1.4) .. (10.93,3.29)   ;
%Straight Lines [id:da9054732726736885] 
\draw [color={rgb, 255:red, 208; green, 2; blue, 27 }  ,draw opacity=1 ]   (397.92,94.58) -- (434,39.33) ;
\draw [shift={(434,39.33)}, rotate = 123.14] [color={rgb, 255:red, 208; green, 2; blue, 27 }  ,draw opacity=1 ][line width=0.75]    (17.64,-3.29) .. controls (13.66,-1.4) and (10.02,-0.3) .. (6.71,0) .. controls (10.02,0.3) and (13.66,1.4) .. (17.64,3.29)(10.93,-3.29) .. controls (6.95,-1.4) and (3.31,-0.3) .. (0,0) .. controls (3.31,0.3) and (6.95,1.4) .. (10.93,3.29)   ;
%Straight Lines [id:da707609788296532] 
\draw    (397.92,94.58) -- (447.33,220.56) ;
%Straight Lines [id:da6636840309759928] 
\draw    (447.33,220.56) -- (421.86,287.08) ;
\draw [shift={(421.14,288.95)}, rotate = 290.95] [color={rgb, 255:red, 0; green, 0; blue, 0 }  ][line width=0.75]    (10.93,-3.29) .. controls (6.95,-1.4) and (3.31,-0.3) .. (0,0) .. controls (3.31,0.3) and (6.95,1.4) .. (10.93,3.29)   ;
%Straight Lines [id:da07543086320798165] 
\draw    (447.33,220.56) -- (429.2,173.49) ;
\draw [shift={(428.48,171.62)}, rotate = 68.93] [color={rgb, 255:red, 0; green, 0; blue, 0 }  ][line width=0.75]    (10.93,-3.29) .. controls (6.95,-1.4) and (3.31,-0.3) .. (0,0) .. controls (3.31,0.3) and (6.95,1.4) .. (10.93,3.29)   ;
%Curve Lines [id:da20259408161659498] 
\draw    (432.41,190.6) .. controls (415.67,211.25) and (424.04,232.73) .. (435.08,242.46) ;
\draw [shift={(436.48,243.62)}, rotate = 217.87] [fill={rgb, 255:red, 0; green, 0; blue, 0 }  ][line width=0.08]  [draw opacity=0] (12,-3) -- (0,0) -- (12,3) -- cycle    ;
\draw [shift={(433.81,188.95)}, rotate = 131.31] [fill={rgb, 255:red, 0; green, 0; blue, 0 }  ][line width=0.08]  [draw opacity=0] (12,-3) -- (0,0) -- (12,3) -- cycle    ;
%Straight Lines [id:da657036584432648] 
\draw [color={rgb, 255:red, 74; green, 144; blue, 226 }  ,draw opacity=1 ][fill={rgb, 255:red, 208; green, 2; blue, 27 }  ,fill opacity=1 ]   (447.33,220.56) -- (479.14,174.29) ;
\draw [shift={(479.14,174.29)}, rotate = 124.51] [color={rgb, 255:red, 74; green, 144; blue, 226 }  ,draw opacity=1 ][line width=0.75]    (17.64,-3.29) .. controls (13.66,-1.4) and (10.02,-0.3) .. (6.71,0) .. controls (10.02,0.3) and (13.66,1.4) .. (17.64,3.29)(10.93,-3.29) .. controls (6.95,-1.4) and (3.31,-0.3) .. (0,0) .. controls (3.31,0.3) and (6.95,1.4) .. (10.93,3.29)   ;
%Straight Lines [id:da5487376524222474] 
\draw    (447.33,220.56) -- (400.48,188.29) ;
%Straight Lines [id:da5581736787393039] 
\draw  [dash pattern={on 0.84pt off 2.51pt}]  (363.14,162.95) -- (400.48,188.29) ;

% Text Node
\draw (223.67,191.07) node [anchor=north west][inner sep=0.75pt]    {$\mathbf{x}_{1}$};
% Text Node
\draw (262.67,141.4) node [anchor=north west][inner sep=0.75pt]  [font=\Large]  {$K$};
% Text Node
\draw (379,95.73) node [anchor=north west][inner sep=0.75pt]    {$\mathbf{x}_{2}$};
% Text Node
\draw (60,257.73) node [anchor=north west][inner sep=0.75pt]  [font=\footnotesize]  {$\mathcal{N}_{K}(\mathbf{x}_{1})$};
% Text Node
\draw (484,15.07) node [anchor=north west][inner sep=0.75pt]  [font=\footnotesize]  {$\mathcal{N}_{K}(\mathbf{x}_{2})$};
% Text Node
\draw (123.33,230.07) node [anchor=north west][inner sep=0.75pt]  [font=\footnotesize]  {$-F(\mathbf{x}_{1})$};
% Text Node
\draw (380.33,28.4) node [anchor=north west][inner sep=0.75pt]  [font=\footnotesize]  {$-F(\mathbf{x}_{2})$};
% Text Node
\draw (374.67,226.07) node [anchor=north west][inner sep=0.75pt]  [font=\footnotesize]  {$\mathcal{T}_{K}(\mathbf{x}_{3})$};
% Text Node
\draw (455,152.4) node [anchor=north west][inner sep=0.75pt]  [font=\footnotesize]  {$-F(\mathbf{x}_{3})$};
% Text Node
\draw (455.67,213.73) node [anchor=north west][inner sep=0.75pt]    {$\mathbf{x}_{3}$};
% Text Node
\draw (318,143.07) node [anchor=north west][inner sep=0.75pt]  [font=\footnotesize]  {$\mathcal{C}(\mathbf{x}_{3} ;K,F)$};

\end{tikzpicture}

%% file: figure/degeneracy2.tex
\tikzset{every picture/.style={line width=0.75pt}} %set default line width to 0.75pt        

\begin{tikzpicture}[x=0.75pt,y=0.75pt,yscale=-1,xscale=1]
%uncomment if require: \path (0,300); %set diagram left start at 0, and has height of 300

%Shape: Polygon [id:ds5436398967614768] 
\draw  [color={rgb, 255:red, 208; green, 2; blue, 27 }  ,draw opacity=1 ][line width=1.5]  (113.76,86.6) -- (168.73,119.34) -- (113.76,206.67) -- (76.7,182.77) -- (58.79,119.34) -- cycle ;
%Shape: Polygon [id:ds9235946213210156] 
\draw  [color={rgb, 255:red, 74; green, 144; blue, 226 }  ,draw opacity=1 ][dash pattern={on 1.69pt off 2.76pt}][line width=1.5]  (116.39,85.29) -- (160.8,112.36) -- (128.37,210.61) -- (70.1,175.13) -- (63.5,123.83) -- cycle ;
%Shape: Free Drawing [id:dp6343064848926974] 
\draw  [line width=3] [line join = round][line cap = round] (195.43,193.14) .. controls (195.62,193.14) and (195.8,193.14) .. (195.98,193.14) ;
%Straight Lines [id:da5312333165506258] 
\draw [color={rgb, 255:red, 208; green, 2; blue, 27 }  ,draw opacity=1 ][line width=1.5]    (142.11,159.85) -- (195.43,193.14) ;
%Straight Lines [id:da4688802965622785] 
\draw [color={rgb, 255:red, 74; green, 144; blue, 226 }  ,draw opacity=1 ][line width=1.5]  [dash pattern={on 1.69pt off 2.76pt}]  (139.36,176.22) -- (195.43,193.14) ;
%Shape: Free Drawing [id:dp005186514068028725] 
\draw  [color={rgb, 255:red, 208; green, 2; blue, 27 }  ,draw opacity=1 ][line width=3] [line join = round][line cap = round] (142.66,160.39) .. controls (142.66,160.39) and (142.66,160.39) .. (142.66,160.39) ;
%Shape: Free Drawing [id:dp5790381087469854] 
\draw  [color={rgb, 255:red, 74; green, 144; blue, 226 }  ,draw opacity=1 ][line width=3] [line join = round][line cap = round] (139.36,176.22) .. controls (139.36,176.22) and (139.36,176.22) .. (139.36,176.22) ;
%Curve Lines [id:da8084225305153967] 
\draw [line width=0.75]  [dash pattern={on 0.75pt off 0.75pt}]  (194.29,150.11) .. controls (191.9,150.48) and (190.48,149.53) .. (190.05,147.27) .. controls (189.39,144.95) and (187.9,144.13) .. (185.58,144.82) .. controls (183.47,145.75) and (181.96,145.2) .. (181.03,143.17) .. controls (179.48,141.23) and (177.82,141.11) .. (176.07,142.8) .. controls (174.88,144.74) and (173.26,145.13) .. (171.23,143.98) .. controls (169.62,142.54) and (167.96,142.48) .. (166.26,143.8) -- (162.5,141.48) -- (156.95,136) ;
\draw [shift={(155.85,134.74)}, rotate = 48.71] [color={rgb, 255:red, 0; green, 0; blue, 0 }  ][line width=0.75]    (6.56,-1.97) .. controls (4.17,-0.84) and (1.99,-0.18) .. (0,0) .. controls (1.99,0.18) and (4.17,0.84) .. (6.56,1.97)   ;
%Shape: Polygon [id:ds26815657559008454] 
\draw  [color={rgb, 255:red, 208; green, 2; blue, 27 }  ,draw opacity=1 ][line width=1.5]  (368.39,100.48) -- (413.54,126.86) -- (368.39,197.19) -- (341.3,179.61) -- (323.24,126.86) -- cycle ;
%Shape: Free Drawing [id:dp7202141240369466] 
\draw  [color={rgb, 255:red, 0; green, 0; blue, 0 }  ,draw opacity=1 ][line width=3] [line join = round][line cap = round] (465.73,157.28) .. controls (465.73,157.28) and (465.73,157.28) .. (465.73,157.28) ;
%Straight Lines [id:da04588033336646058] 
\draw [color={rgb, 255:red, 208; green, 2; blue, 27 }  ,draw opacity=1 ][line width=1.5]    (413.54,126.86) -- (465.73,156.84) ;
%Shape: Free Drawing [id:dp37315345430251146] 
\draw  [color={rgb, 255:red, 208; green, 2; blue, 27 }  ,draw opacity=1 ][line width=3] [line join = round][line cap = round] (413.8,126.51) .. controls (413.8,126.66) and (413.8,126.8) .. (413.8,126.95) ;
%Straight Lines [id:da5988038385003105] 
\draw [color={rgb, 255:red, 126; green, 211; blue, 33 }  ,draw opacity=1 ][line width=1.5]  [dash pattern={on 5.63pt off 4.5pt}]  (453.53,65.41) -- (413.54,126.86) ;
%Straight Lines [id:da14838515525409135] 
\draw [color={rgb, 255:red, 126; green, 211; blue, 33 }  ,draw opacity=1 ][line width=1.5]  [dash pattern={on 5.63pt off 4.5pt}]  (404.77,42.11) -- (368.39,100.48) ;
%Straight Lines [id:da12249529267557757] 
\draw [color={rgb, 255:red, 126; green, 211; blue, 33 }  ,draw opacity=1 ][line width=1.5]  [dash pattern={on 5.63pt off 4.5pt}]  (333.43,45.62) -- (368.39,100.48) ;
%Straight Lines [id:da08385366359793478] 
\draw [color={rgb, 255:red, 126; green, 211; blue, 33 }  ,draw opacity=1 ][line width=1.5]  [dash pattern={on 5.63pt off 4.5pt}]  (287.38,73.32) -- (323.24,126.86) ;
%Straight Lines [id:da25877654383796145] 
\draw [color={rgb, 255:red, 126; green, 211; blue, 33 }  ,draw opacity=1 ][line width=1.5]  [dash pattern={on 5.63pt off 4.5pt}]  (260.29,143.21) -- (323.24,126.86) ;
%Straight Lines [id:da46331853452045557] 
\draw [color={rgb, 255:red, 126; green, 211; blue, 33 }  ,draw opacity=1 ][line width=1.5]  [dash pattern={on 5.63pt off 4.5pt}]  (274.67,196.67) -- (341.3,179.61) ;
%Straight Lines [id:da33595265458983326] 
\draw [color={rgb, 255:red, 126; green, 211; blue, 33 }  ,draw opacity=1 ][line width=1.5]  [dash pattern={on 5.63pt off 4.5pt}]  (309.05,219.7) -- (341.3,179.61) ;
%Straight Lines [id:da2716729991316227] 
\draw [color={rgb, 255:red, 126; green, 211; blue, 33 }  ,draw opacity=1 ][line width=1.5]  [dash pattern={on 5.63pt off 4.5pt}]  (340.2,238.61) -- (368.39,197.19) ;
%Straight Lines [id:da197933950851632] 
\draw [color={rgb, 255:red, 126; green, 211; blue, 33 }  ,draw opacity=1 ][line width=1.5]  [dash pattern={on 5.63pt off 4.5pt}]  (421.03,229.38) -- (368.39,197.19) ;
%Straight Lines [id:da8565963570282511] 
\draw [color={rgb, 255:red, 126; green, 211; blue, 33 }  ,draw opacity=1 ][line width=1.5]  [dash pattern={on 5.63pt off 4.5pt}]  (483.79,166.95) -- (413.54,126.86) ;
%Shape: Polygon [id:ds08435384668280554] 
\draw  [color={rgb, 255:red, 74; green, 144; blue, 226 }  ,draw opacity=1 ][dash pattern={on 1.69pt off 2.76pt}][line width=1.5]  (375.4,99.2) -- (406.1,114.14) -- (385.33,199.42) -- (333.86,173.93) -- (328.18,125.04) -- cycle ;
%Straight Lines [id:da4768218463623777] 
\draw [color={rgb, 255:red, 74; green, 144; blue, 226 }  ,draw opacity=1 ][line width=1.5]  [dash pattern={on 1.69pt off 2.76pt}]  (398.9,143.65) -- (465.73,156.84) ;
%Shape: Free Drawing [id:dp802040955279782] 
\draw  [color={rgb, 255:red, 74; green, 144; blue, 226 }  ,draw opacity=1 ][line width=3] [line join = round][line cap = round] (399.8,143.21) .. controls (399.8,143.21) and (399.8,143.21) .. (399.8,143.21) ;

% Text Node
\draw (202.31,191.43) node [anchor=north west][inner sep=0.75pt]    {$\mathbf{x}^{*}$};
% Text Node
\draw (47.24,88.04) node [anchor=north west][inner sep=0.75pt]  [font=\footnotesize,color={rgb, 255:red, 208; green, 2; blue, 27 }  ,opacity=1 ]  {$K(\mathbf{p} *)$};
% Text Node
\draw (138.74,82.59) node [anchor=north west][inner sep=0.75pt]  [font=\footnotesize,color={rgb, 255:red, 74; green, 144; blue, 226 }  ,opacity=1 ]  {$K(\mathbf{p})$};
% Text Node
\draw (124.45,141.72) node [anchor=north west][inner sep=0.75pt]  [font=\footnotesize,color={rgb, 255:red, 208; green, 2; blue, 27 }  ,opacity=1 ]  {$\mathbf{x}_{\pi }^{*}$};
% Text Node
\draw (139.61,178.4) node [anchor=north west][inner sep=0.75pt]  [font=\footnotesize,color={rgb, 255:red, 74; green, 144; blue, 226 }  ,opacity=1 ]  {$\mathbf{x}_{\pi }$};
% Text Node
\draw (163.3,151.86) node [anchor=north west][inner sep=0.75pt]  [font=\scriptsize] [align=left] {"sticky face"};
% Text Node
\draw (465.03,137.26) node [anchor=north west][inner sep=0.75pt]    {$\mathbf{x}^{*}$};
% Text Node
\draw (418.19,115.8) node [anchor=north west][inner sep=0.75pt]  [font=\footnotesize,color={rgb, 255:red, 208; green, 2; blue, 27 }  ,opacity=1 ]  {$\mathbf{x}_{\pi }^{*}$};
% Text Node
\draw (379.01,132.16) node [anchor=north west][inner sep=0.75pt]  [font=\footnotesize,color={rgb, 255:red, 74; green, 144; blue, 226 }  ,opacity=1 ]  {$\mathbf{x}_{\pi }$};
% Text Node
\draw (299.06,148.57) node [anchor=north west][inner sep=0.75pt]  [font=\footnotesize,color={rgb, 255:red, 208; green, 2; blue, 27 }  ,opacity=1 ]  {$K(\mathbf{p} *)$};
% Text Node
\draw (321.93,98.02) node [anchor=north west][inner sep=0.75pt]  [font=\footnotesize,color={rgb, 255:red, 74; green, 144; blue, 226 }  ,opacity=1 ]  {$K(\mathbf{p})$};

\end{tikzpicture}